\documentclass[english, twoside, 12pt]{smfart}
\usepackage{fullpage}
\usepackage{amsfonts}
\usepackage{amssymb}
\input psfig
\def\IMSmarkvadjust{0 pt}
\def\IMSmarkhadjust{0 pt}
\def\SBIMSMark#1#2#3{
 \font\SBF=cmss10 at 10 true pt
 \font\SBI=cmssi10 at 10 true pt
 \setbox0=\hbox{\SBF Stony Brook IMS Preprint \##1}
 \setbox2=\hbox to \wd0{\hfil \SBI #2}
 \setbox4=\hbox to \wd0{\hfil \SBI #3}
 \setbox6=\hbox to \wd0{\hss
             \vbox{\hsize=\wd0 \parskip=0pt \baselineskip=10 true pt
                   \copy0 \break%
                   \copy2 \break%
                   \copy4 \break}}
 \dimen0=\ht6   \advance\dimen0 by \vsize \advance\dimen0 by 8 true pt
                \advance\dimen0 by -\pagetotal
	        \advance\dimen0 by \IMSmarkvadjust
 \dimen2=\hsize \advance\dimen2 by .25 true in
	        \advance\dimen2 by \IMSmarkhadjust

  \openin2=publishd.tex
  \ifeof2\setbox0=\hbox to 0pt{}
  \else 
     \setbox0=\hbox to 3.1 true in{
                \vbox to \ht6{\hsize=3 true in \parskip=0pt  \noindent  
                {\SBI Published in modified form:}\hfil\break
                \input publishd.tex 
                \vfill}}
  \fi
  \closein2
  \ht0=0pt \dp0=0pt
 \ht6=0pt \dp6=0pt
 \setbox8=\vbox to \dimen0{\vfill \hbox to \dimen2{\copy0 \hss \copy6}}
 \ht8=0pt \dp8=0pt \wd8=0pt
 \copy8
 \message{*** Stony Brook IMS Preprint #1, #2. #3 ***}
}

\newcommand\Or{{\mathcal O}}
\newcommand\po{{\mathcal P}}
\newcommand\ra{{\mathcal R}}
\newcommand\ssm{\smallsetminus}
\newcommand\QED{$\quad\square$}
\newcommand\Ref{\hangindent=1pc \hangafter=1 \noindent}
\newcommand\q{{\mathcal Q}}
\newcommand\I{{~\mathcal I}}
\newcommand\Z{{\mathbb Z}}
\newcommand\D{{\mathbb D}}
\newcommand\R{{\mathbb R}}
\newcommand\Q{{\mathbb Q}}
\newcommand\C{{\mathbb C}}
\newcommand\rot{{\rm rot}}
\newcommand\Rot{{\rm Trans}}
\newcommand\bs{\bigskip}
\newcommand\nin{\noindent}
\newcommand\hts{{\rm~(higher~terms)}}
\newcommand\cc{{\bigcirc\!\!\!\!\!{\bf c}\,}}
\newcommand\+{\,+\,}
\newcommand\iso{{\buildrel \approx\over\to}}
\newcommand\Per{{\rm Per}}

\newcommand\cl{\centerline}
\newcommand\ms{\medskip}
\newcommand\ssk{\smallskip}
\newcommand\QP{\smallskip\leftskip=.4in\rightskip=.4in\noindent}
\DeclareOldFontCommand{\bit}{\normalfont\sffamily\slshape}{\mathsf}
\newcommand\tf{\bfseries\itshape}

\begin{document}
\markright{Periodic Orbits, Externals Rays and the Mandelbrot Set}

\title[Periodic Orbits, Externals Rays and the Mandelbrot Set]
{Periodic Orbits, Externals Rays and the Mandelbrot Set:
An Expository Account}\ms

\author[J. Milnor]{John Milnor}\ssk 

\maketitle
\thispagestyle{empty}
\SBIMSMark{1999/3}{May 1999}{}

{\QP\small{\tf Abstract} -
A presentation of some fundamental results
from the Douady-Hubbard theory of the Mandelbrot set, based on the
idea of ``orbit portrait'': the pattern of external rays landing on a
periodic orbit for a quadratic polynomial map.\par}
 
{\QP\small{\tf R\'esum\'e (Orbites p\'eriodiques, rayons externes et
l'ensemble de Mandelbrot: un compte-rendu)} - Nous expliquons
quelques r\'esultats fondamentaux de Douady-Hubbard sur l'ensemble de
Mandelbrot en utilisant l'id\'ee de ``portrait orbital''
c'est-\`a-dire le mod\`ele des rayons externes qui aboutissent sur une
orbite p\'eriodique d'une application polynomiale quadratique.\par}\bs

\cl{\tf
Dedicated to Adrien Douady on the occasion of his sixtieth birthday.}
\bs

\begin{quote}{ Contents:

1. Introduction

2. Orbit portraits

3. Parameter rays

4. Near parabolic maps

5. The period $n$ curve in $($parameter$\times$dynamic$)$ space

6. Hyperbolic components

7. Orbit forcing

8. Renormalization

9. Limbs and the satellite orbit

Appendix A. Totally disconnected Julia sets and the exterior of $M$.

Appendix B. Computing rotation numbers.

References\ssk}\end{quote}

\section{Introduction}

\begin{figure}
\cl{\psfig{figure=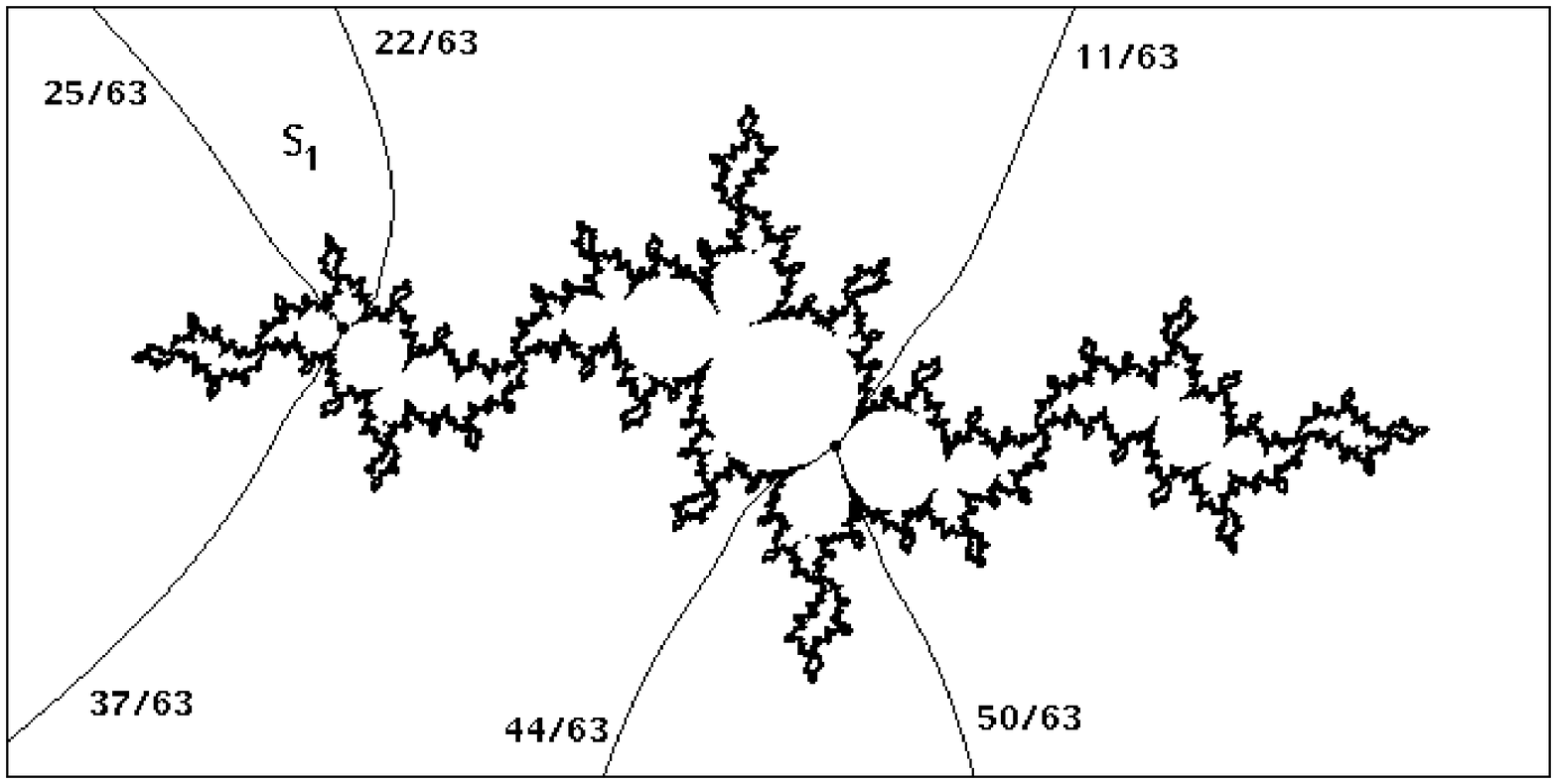,height=2in}}
\ssk
\begin{quotation}
{\bit Figure 1. Julia set for $z\mapsto z^2+ 
({1\over 4}\,e^{2\pi i/3}-1)$
showing the six rays landing on a period two parabolic orbit. The associated
orbit portrait has characteristic arc $\I=(22/63,25/63)$
and valence $v=3$ rays per orbit point.\ssk}
\end{quotation}\end{figure}

A key point in Douady and Hubbard's study of the Mandelbrot set $M$ is the
theorem that every parabolic point $c\ne 1/4$ in $M$
is the landing point for exactly two external rays with angles which are
periodic under doubling. (See [DH2]. By
definition, a parameter point is {\bit parabolic\/} if and only if the
corresponding quadratic map has a periodic orbit with some root of unity as
multiplier.)
This note will try to provide a proof of this result and some of its
consequences which relies as much as possible on
elementary combinatorics, rather than on more difficult analysis.
It was inspired by \S2 of the recent thesis
of Schleicher [S1], which contains very substantial simplifications of the
Douady-Hubbard proofs with a much more compact argument, and is highly
recommended. (See also [S2], [LS].)
The proofs given here are rather different from those of Schleicher, and are
based on a combinatorial study of the angles of
external rays for the Julia set which land on periodic orbits.
(Compare [A], [GM].)
As in [DH1], the basic idea is to find properties of $M$
by a careful study of the dynamics for parameter values outside of $M$.
The results in this paper are mostly well known; there is a particularly
strong overlap with [DH2]. The only claim to
originality is in emphasis, and the organization of the proofs. (Similar
methods can be used for higher degree polynomials with only one critical point.
Compare [S3], [E], and see [PR] for a different approach. For
a theory of polynomial maps which may have many
critical points, see [K].)

We will assume some familiarity with the classical Fatou-Julia
theory, as described for example in [Be], [CG], [St], or [M2].\ssk

{\tf Standard Definitions.} (Compare Appendix A.)
Let $K=K(f_c)$ be the {\bit filled Julia set\/},
that is the union of all bounded orbits, for the quadratic map
$$	f(z)~=~f_c(z)~=~z^2+c~. $$
Here both the parameter $c$ and the dynamic variable $z$
range over the complex numbers. The {\bit Mandelbrot set\/} $M$ can
be defined as the compact subset of the {\bit parameter plane\/}
(or $c$-plane) consisting of
all complex numbers $c$ for which $K(f_c)$ is connected.
We can also identify the complex number $c$ with one particular point
in the {\bit dynamic plane\/} (or $z$-plane),
namely the {\bit critical value\/}
$f_c(0)=c$ for the map $f_c$. The parameter $c$ belongs to $M$
if and only if the orbit $f_c:0\mapsto c\mapsto c^2+c\mapsto\cdots$ is
bounded, or in other words if and only if $0, c\in K(f_c)$.
Associated with each of the compact sets $K=K(f_c)$ in the dynamic plane
there is a {\bit potential function\/} or {\bit Green's function\/} $G^K:
\C\to[0,\infty)$ which vanishes precisely on $K$, is harmonic off
$K$, and is asymptotic to $\log|z|$
near infinity. The family of {\bit external rays\/} of $K$
can be described as the orthogonal trajectories of the level curves $G^K=
{\rm constant}$. Each such ray which extends to infinity can be specified
by its angle at infinity $t\in\R/\Z$, and will be denoted by  $\ra_t^K$.
Here $c$ may be either in or outside of the Mandelbrot set.
Similarly, we can consider the potential function $G^M$ and the
external rays $\ra_t^M$ associated with the Mandelbrot set.
We will use the term {\bit dynamic ray\/} (or briefly
$K$-{\bit ray\/}) for an external
ray of the filled Julia set, and {\bit parameter ray\/} (or briefly $M$-{\bit
ray\/}) for an external ray of the Mandelbrot set. (Compare [S1], [S2].)\ssk

\begin{figure}[htb]
\cl{\psfig{figure=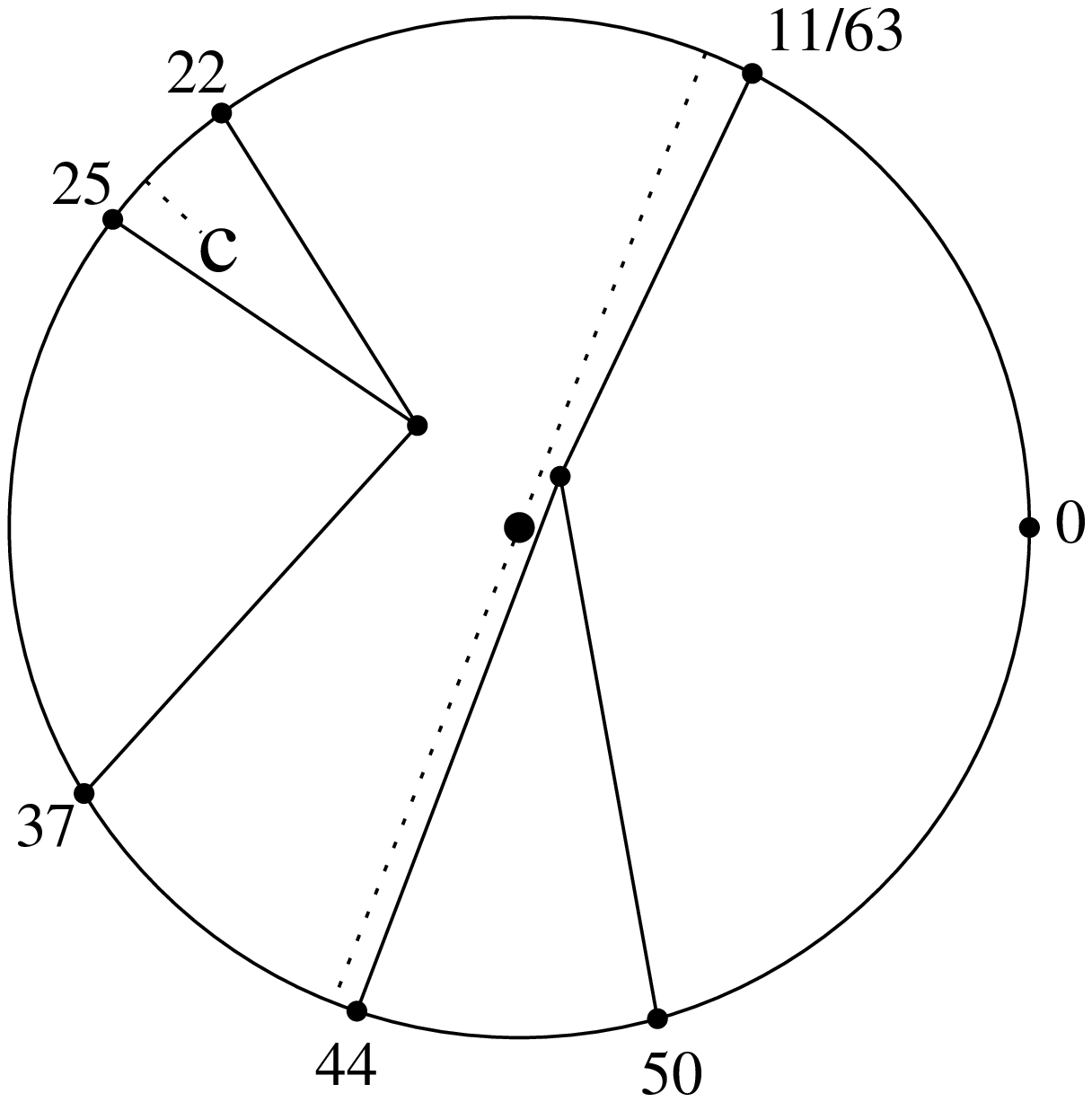,height=2.4in}}
\cl{\bit Figure 2. Schematic diagram illustrating the orbit portrait $(1)$.}
\end{figure}

{\tf Definition.} Let $\Or=\{z_1\,,\,\ldots\,,\,z_p\}$
be a periodic orbit for $f$.
Suppose that there is some rational angle $t\in\Q/\Z$
so that the dynamic ray $\ra_t^{K(f)}$ lands at
a point of $\Or$.
Then for each $z_i\in\Or$ the collection
$A_i$ consisting of all angles of dynamic
rays which land at the point $z_i$ is a finite and non-vacuous
subset of $\Q/\Z$. 
The collection $\{A_1\,,\,\ldots\,,\,A_p\}$
will be called
the {\bit orbit portrait\/} $\po=\po(\Or)$. As an example, Figure 1 shows
a quadratic Julia set having a parabolic orbit with portrait
$$	\po~=~\big\{\;\{22/63\,,\,25/63\,,\,37/63\}\;,~
	\;\{11/63\,,\,44/63\,,\,50/63\}\;\big\}~. \eqno (1)$$
It is often convenient to represent such a portrait by a schematic diagram,
as shown in Figure 2. (For details,
and an abstract characterization of orbit portraits,
see \S2.) 

The number of elements in each $A_i$ (or in other words the number
of $K$-rays which land on each orbit point)
will be called the {\bit valence\/} $v$.
Let us assume that $v\ge 2$. Then the $v$ rays landing at $z$
cut the dynamic plane up into $v$ open regions which will be called the
{\bit sectors\/} based at the orbit
point $z\in\Or$. 
The {\bit angular width\/} of a sector $S$
will mean the length of the open arc $I_S$ consisting of
all angles $t\in \R/\Z$ with $\ra_t^K\subset S$. (We use the word `arc'
to emphasize that we will identify $\R/\Z$ with the `circle at infinity'
surrounding the plane of complex numbers.)
Thus the sum of the angular widths
of the $v$ distinct sectors based at an orbit point $z$
is always equal to $+1$. The following result will be proved in 2.11.\ssk


\nin{\tf Theorem 1.1. The Critical Value Sector $S_1$.} {\it Let $\Or$
be an orbit of period $p\ge 1$ for $f=f_c$. If there are $v\ge 2$
dynamic rays landing at each point of $\Or$, then there is one and only one
sector $S_1$ based at some point $z_1\in\Or$ which contains the
critical value $c=f(0)$, and whose closure contains no point other than
$z_1$ of the orbit $\Or$. This critical value sector $S_1$
can be characterized, among all of the $pv$ sectors based at the
various points of $\Or$, as the unique sector of smallest angular width.}\ssk

It should be emphasized that this description is correct whether the filled
Julia set $K$ is connected or not.\ssk

Our main theorem can be stated as follows.
Suppose that there exists some polynomial $f_{c_0}$ which admits
an orbit $\Or$ with portrait $\po$, again having valence $v\ge 2$.
Let\break
$0<t_-<t_+<1$ be the angles of the two dynamic rays $\ra^K_{t_\pm}$
which bound the critical value sector $S_1$ for $f_{c_0}$.\ssk

{\nin{\tf Theorem 1.2. The Wake $W_\po$.} \it The two corresponding
parameter rays $\ra_{t_\pm}^M$ land at a single point ${\bf r}_\po$ of the
parameter plane. These rays, together with their landing point, cut the plane
into two open subsets $W_\po$ and
$\C\ssm \overline W_\po$ with the following
property: A quadratic map $f_c$ has a repelling orbit with portrait $\po$
if and only if $c\in W_\po$, and has a parabolic orbit with portrait $\po$ if
and only if $c={\bf r}_\po$.\ssk}

\begin{figure}[htb]
\cl{\psfig{figure=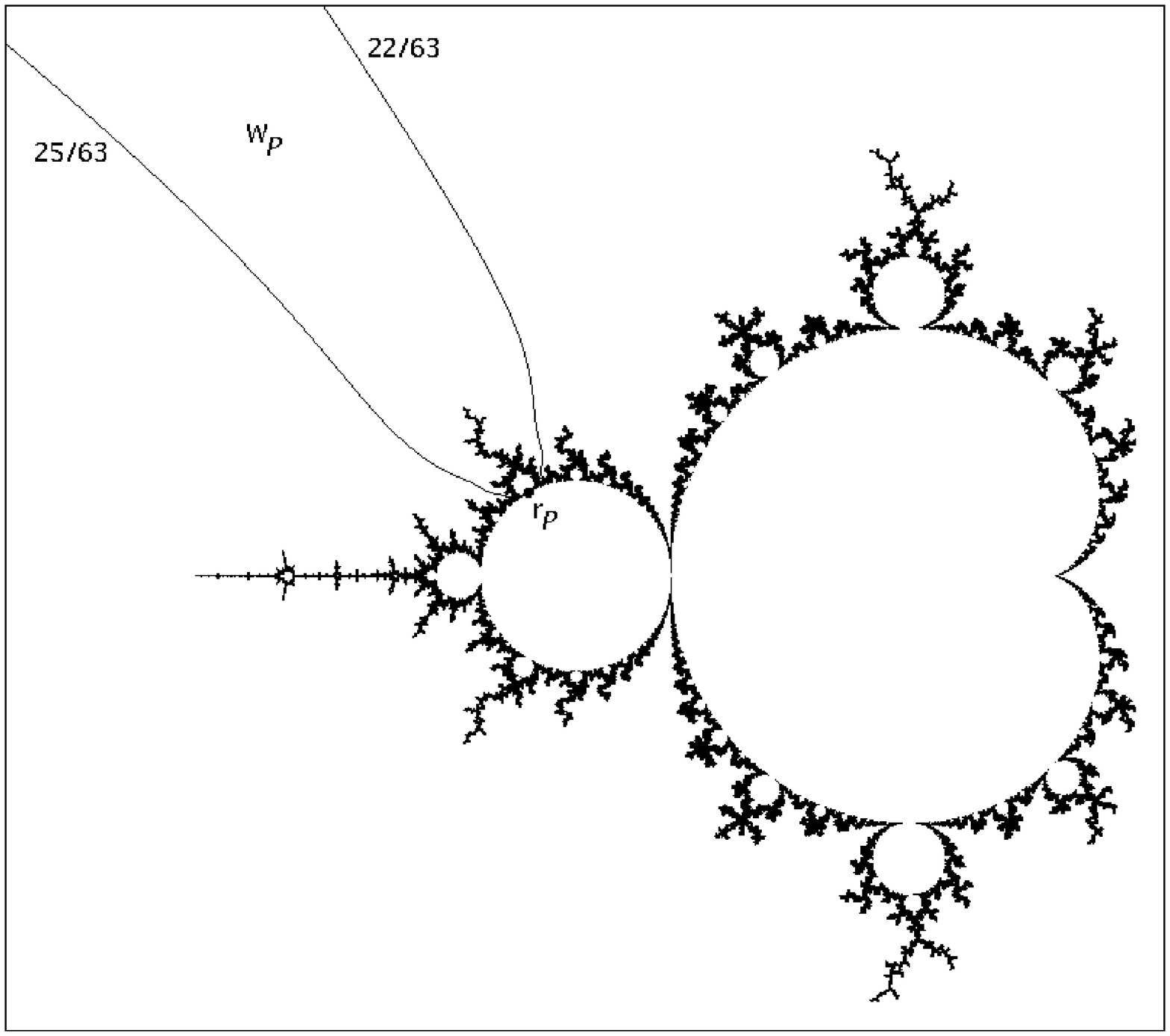,height=2.5in}}
\ssk
\begin{quote}{\bit Figure 3. The boundary of the Mandelbrot set, showing the
wake $W_\po$ and the root point ${\bf r}_\po=
{1\over 4}\,e^{2\pi i/3}-1$ associated with the orbit portrait
of Figure 1, with characteristic arc $\I_\po=(22/63,25/63)$.\par}
\vskip -.1in
\end{quote}\end{figure}

In fact this will follow by combining the assertions 3.1, 4.4, 4.8, and 5.4
below.\ms

{\tf Definitions.} This open set $W_\po$ will be called the
{\bit$\po$-wake\/} in parameter space (compare Atela [A]),
and ${\bf r}_\po$ will be called the {\bit root point\/} of this wake.
The intersection
$M_\po=M\cap\overline W_\po$ will be called the {\bit $\po$-limb\/} of
the Mandelbrot set. The open arc $I_{S_1}=(t_-\,,\,t_+)$ consisting
of all angles of dynamic rays $\ra^K_t$ which are contained in the interior
of $S_1$,
or all angles of parameter rays $\ra_t^M$ which are contained in $W_\po$,
will be called the {\bit characteristic arc\/} $\I=\I_\po$
for the orbit portrait $\po$. (Compare 2.6.)

In general, the orbit portraits with valence $v=1$ are of little interest to
us. These portraits certainly exist. For example, for the base map
$f_0(z)=z^2$ which lies outside of every wake,
{\it every\/} orbit portrait has valence $v=1$. As we follow a path
in parameter space which crosses into the wake $W_\po$ through its root point,
either one orbit
with a portrait of valence one degenerates to form an orbit of lower
period with portrait $\po$, or else two different orbits with portraits
of valence one fuse together to form an orbit with portrait $\po$. (If we
cross into $W_\po$ through a parameter ray
$\ra_{t_\pm}^M$, the picture is similar except
that the landing point of the dynamic
ray $\ra_{t_\pm}^K$ jumps discontinuously. If $t_+$ and $t_-$ belong
to the same cycle under angle doubling, then the landing points of both of
these dynamics rays jump discontinuously.)

However, there is one exceptional portrait of valence one:
The {\bit zero portrait\/} $\po=\{\{0\}\}$ will play an important role.
It is not difficult to check that the dynamic ray $\ra_0^{K}$ of angle $0$
for $f_c$ lands at a well defined fixed point
if and only if the parameter value $c$ lies in the complement of the
parameter ray $\ra_0^M=\ra_1^M=(1/4,\infty)$. Furthermore, this fixed
point necessarily has portrait $\{\{0\}\}$. Thus the
wake, consisting of all $c\in\C$ for which $f_c$ has a repelling fixed
point with portrait $\{\{0\}\}$, is just the complementary region $\C\ssm
[1/4\,,\,\infty)$. The characteristic
arc $\I_{\{\{0\}\}}$ for this portrait, consisting of all angles $t$
such that $\ra_t^K\subset W_{\{\{0\}\}}$, is the open interval $(0,1)$,
and the root point ${\bf r}_{\{\{0\}\}}$, the unique
parameter value $c$ such that $f_c$ has a parabolic fixed point with
portrait $\{\{0\}\}$, is the landing point $c=1/4$ for the zero
parameter ray.\ssk

\begin{figure}[htb]
\cl{\psfig{figure=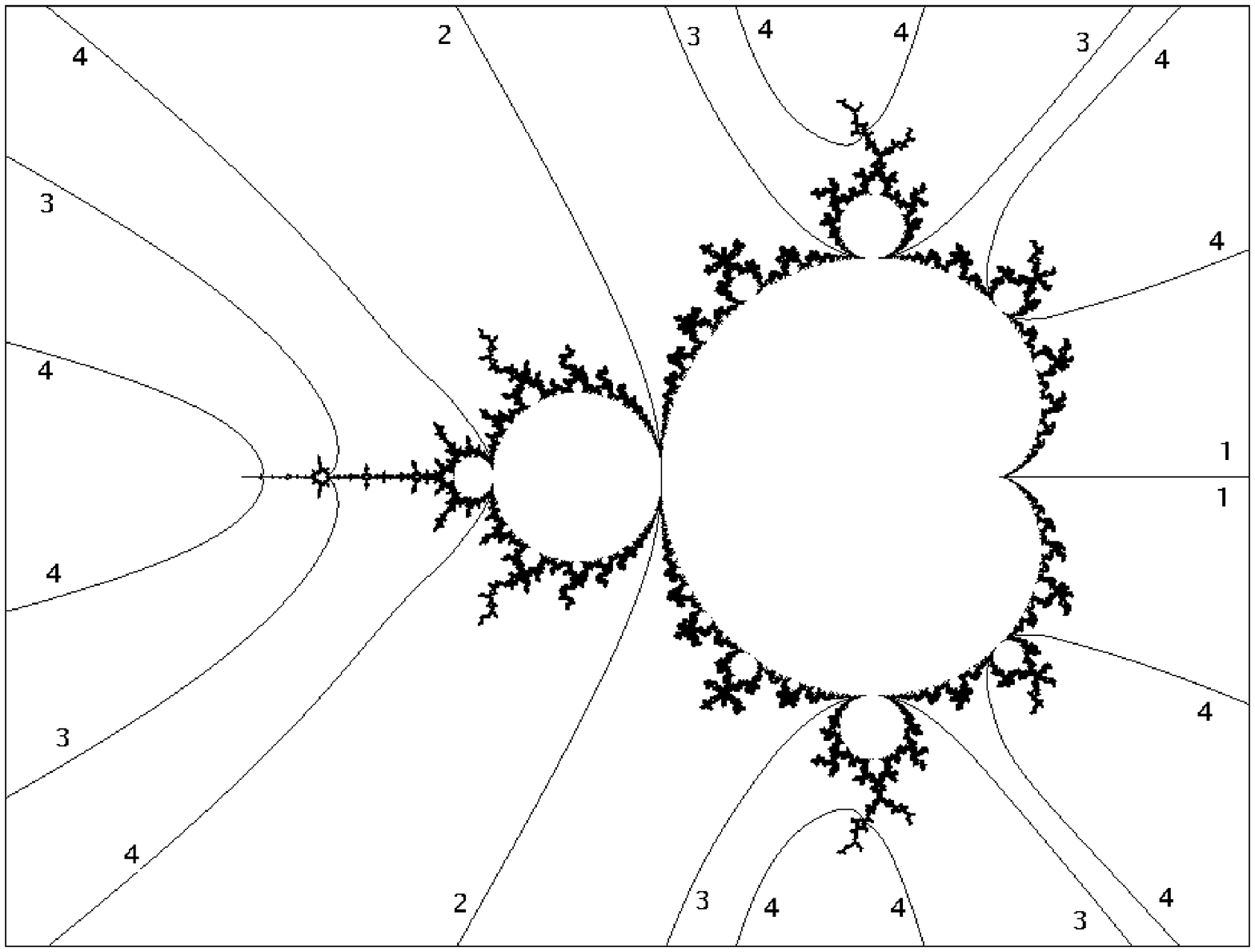,height=2.8in}}
\ssk
\cl{\bit Figure 4. Boundaries of the wakes of ray period four or less.}
\end{figure}

{\tf Definition.}
It will be convenient to say that a portrait $\po$ is {\bit non-trivial\/}
if it either has valence $v\ge 2$ or is equal to this zero portrait.\ssk

{\tf Remark.} An alternative characterization would be the following. An
orbit portrait $\{A_1\,,\,\ldots\,,\,\discretionary{}{}{} A_p\}$ is
non-trivial if and only if it 
is {\bit maximal\/}, in the sense that there is no orbit portrait
$\{A'_1\,,\,\ldots\,,\,A'_q\}$ with $A'_1\mathop{\supset}\limits_{\not=}
A_1$. This statement follows easily from 1.5 and 2.7 below. Still another
characterization would be that $\po$ is non-trivial if and only if it is the
portrait of some parabolic orbit. (See 5.4.)\ssk

{\nin{\tf Corollary 1.3. Orbit Forcing.} \it If $\po$ and
$\q$ are two distinct non-trivial orbit portraits,
then the boundaries $\partial W_\po$ and $\partial W_\q$ of the
corresponding wakes are disjoint subsets of $\C$.
Hence the closures
$\overline W_\po$ and $\overline W_\q$ are either disjoint
or strictly nested. In particular, if $\I_\po\subset \I_\q$ with $\po\ne\q$,
then it follows that $\overline W_\po\subset W_\q$.\ssk}

Thus whenever $\overline \I_\po\subset \I_\q$,
the existence of a repelling or parabolic orbit with portrait $\po$
{\bit forces\/} the existence of a {\it repelling\/}
orbit with portrait $\q$. We will write briefly $\po\Rightarrow\q$.
On the other hand, if $\I_\po\cap \I_\q=\emptyset$
then no $f_c$ can have both an orbit with portrait $\po$ and an orbit with
portrait $\q$.

See Figure 5 for a schematic description of orbit forcing relations for
orbits with ray period 4 or less,
corresponding to the collection of wakes illustrated in Figure 4. (Evidently
this diagram, as well as analogous diagrams in which higher periods are
included, has a tree structure, with no loops.)
\ssk

\begin{figure}[htb]
\cl{\psfig{figure=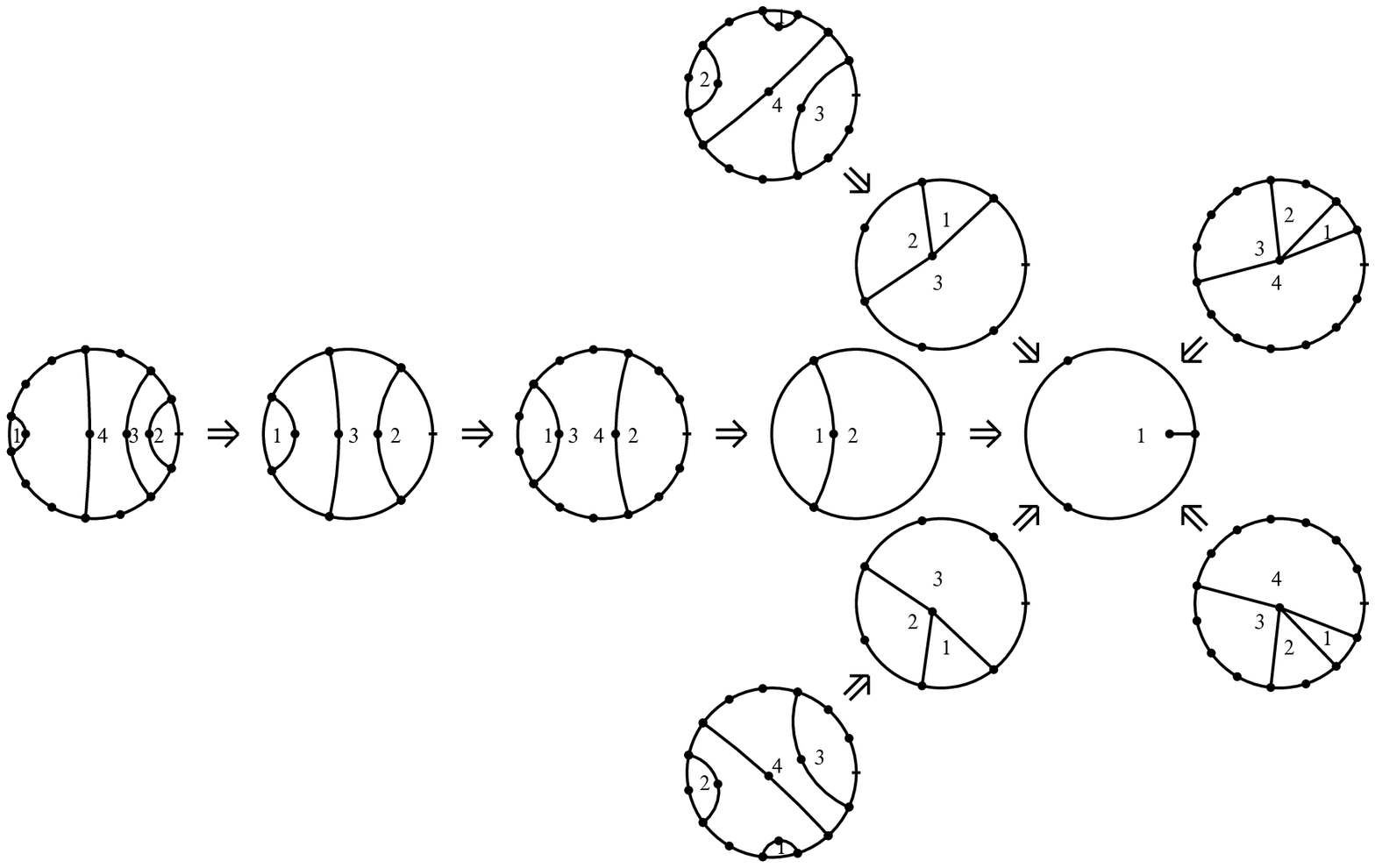,height=3.3in}}\ssk
\begin{quote}
{\bit Figure 5. Forcing tree for the non-trivial orbit portraits of ray period
$n\le 4$.  Each disk in this figure contains a schematic diagram
of the corresponding orbit portrait, with the first $n$
forward images of the critical value sector labeled. (Compare Figure 4;
and compare the ``disked-tree model''
for the Mandelbrot set in Douady [D5].)\ssk}
\end{quote}\end{figure}

{\tf Proof of 1.3,
assuming 1.2.} First note that $W_\po$ and $W_\q$ cannot have a
boundary ray in common. For the landing point of such a common ray would
have to have one parabolic orbit with portrait $\po$ and one parabolic orbit
with portrait $\q$. But a quadratic map, having only one critical point,
cannot have two distinct parabolic orbits. In fact this argument shows that
$\partial W_\po\cap\partial W_\q=\emptyset$. Note that the parameter point
$c=0$ (corresponding to the map $f_0(z)=z^2$) does not belong to any wake
$W_\po$ with $\po\ne\{\{0\}\}$.
Since rays cannot cross each other, it follows easily that either
$$\overline W_\po\subset W_\q\,,\quad{\rm or}
\quad\overline W_\q\subset W_\po\,,\quad{\rm or} \quad\overline W_\po\cap
 \overline W_\q=\emptyset~,$$
as required.\QED\ssk

For further discussion and a more direct proof, see \S7.\ssk

To fill out the picture, we also need the following two statements. To any
orbit portrait $\po=\{A_1\,,\,\ldots\,,\,A_p\}$ we associate not only its
{\bit orbit period\/} $p$ but also its {\bit ray period\/} $rp$, that is
the period of the angles $t\in A_i$ under doubling modulo one. In many cases,
$rp$ is a proper multiple of $p$. (Compare Figure 1.) Suppose in particular
that $c\in M$ is a parabolic parameter value, that is suppose that $f_c$
has a periodic orbit where the multiplier is an $r$-th root of unity,
$r\ge 1$. Then one can show that the ray period for the associated
portrait is equal to the product $rp$. (See for example [GM].) 
This is also the period of the Fatou component containing the critical point.
This ray period $rp$ is the most important parameter associated with a
parabolic point $c$ or with a wake $W_\po$.

It follows from 1.2 that every non-trivial portrait which occurs at all must
occur as the portrait of some uniquely determined
parabolic orbit. The converse statement will be proved in 4.8:

{\ssk\nin{\tf Theorem 1.4. Parabolic Portraits are Non-Trivial.} \it
If $c$ is any parabolic point in $M$,
then the portrait $\po=\po(\Or)$ of its parabolic orbit is a non-trivial
portrait. That is,
if we exclude the special case $c=1/4$, then at least two
$K$-rays must land on each parabolic orbit point.\ssk}

It then follows immediately from 1.2 that the parabolic parameter point
$c$ must be equal to the root point ${\bf r}_\po$ of an associated wake.
It also follows from 1.2 that the angles of the $M$-rays which bound a wake
$W_\po$ are always periodic under doubling. In \S5 
we use a simple counting argument to prove the converse statement.
(This imitates Schleicher, who uses a similar counting argument in
a different way.)

{\ssk\nin{\tf Theorem 1.5. Every Periodic Angle Occurs.} \it
If $t\ne 0$ in $\R/\Z$ is periodic under doubling, then $\ra_t^M$
is one of the two boundary rays of some (necessarily unique) wake.\ssk}

Further consequences of these ideas will be developed in \S6
which shows that each wake contains a uniquely associated hyperbolic component,
\S8 which describes how each wake contains an associated small copy of the
Mandelbrot set, and \S9 which shows that each limb is connected even
if its root point is removed. There are two appendices giving further
supporting details.\ssk

{\bf Acknowledgement.} I want to thank M. Lyubich and D. Schleicher for their
ideas, which play a basic role in this presentation. I am particularly
grateful to Schleicher, to S. Zakeri, and to Tan Lei
for their extremely helpful criticism of the manuscript.
Also, I want to thank both the
Gabriella and Paul Rosenbaum Foundation and the National Science
Foundation (Grant DMS-9505833)
for their support of mathematical activities at Stony Brook.


\section{Orbit Portraits.} This section will begin the proofs by
describing the basic properties of orbit portraits. We will need the following.
Let $f(z)=z^2+c$ with filled Julia set $K$.

{\ssk\nin{\tf Lemma 2.1. Mapping of Rays.} \it If a dynamic ray $\ra_t^K$
lands at a point $z\in \partial K$, then the image ray $f(\ra_t^K)=
\ra_{2t}^K$ lands at the image point $f(z)$.  Furthermore,
if three or more rays $\ra_{t_1}^K\,,\,\ra_{t_2}^K\,,\,
\ldots\,,\,\ra_{t_k}^K$ land at $z\ne 0$, 
then the cyclic order of the angles $t_i$ around the circle $\R/\Z$ is
the same as the cyclic order of the doubled angles $2t_i~~({\rm mod}~\Z)$
around $\R/\Z$.\ssk}


{\tf Proof.} Since each $\ra_{t_j}^K$ is assumed to be a
smooth ray, it cannot pass through any precritical point. Hence $\ra_{2t_j}^K$
also cannot pass through a precritical point, and must be a smooth ray landing
at $f(z)$. Now suppose that we are given three or more rays with angles
$0\le t_1<t_2<\cdots<t_k<1$, all landing at $z$. These rays, together
with their landing point, cut the plane up into sectors $S_1\,,\,\ldots\,,\,
S_k$, where each $S_i$ is bounded by $\ra_{t_i}^K$ and $\ra_{t_{i+1}}^K$
(with subscripts modulo $k$). The cyclic ordering of these various rays
can be measured within an arbitrarily small neighborhood of the landing point
$z$, since any transverse arc which crosses $\ra_{t_i}^K$ in the positive
direction must pass from $S_{i-1}$ to $S_i$. Since $f$ maps a
neighborhood of $z$ to a neighborhood of $f(z)$ by an orientation
preserving diffeomorphism, it follows that the image rays must have the
same cyclic order.\QED

Now let us impose the following.\ssk

{\tf Standing Hypothesis 2.2.} $\Or=\{z_1\,,\,\ldots\,,\,z_p\}$ is a periodic
orbit for a quadratic map $f_c(z)=z^2+c$, with orbit points
numbered so that $f(z_j)=z_{j+1}$, taking subscripts modulo $p$.
Furthermore there is at least
one rational angle $t\in\Q/\Z$ so that the dynamic ray $\ra_t^K$
associated with $f$ lands at some point of this orbit $\Or$.\ssk

If $c$ belongs to the Mandelbrot set $M$, or in other words if the
filled Julia set $K$ is
connected, then this condition will be satisfied if and only if the orbit
$\Or$ is either repelling or parabolic. (Compare [Hu], [M3].) On the other
hand, for $c\not\in M$, all periodic orbit are repelling, but the
condition may fail to be satisfied either because the rotation number is
irrational (compare [GM, Figure 16]),
or because the $K$-rays which `should' land on
$\Or$ bounce off precritical points en route ([GM, Figure 14]).\ssk

As in \S1, let $A_j\subset\R/\Z$ be the set of all angles of
$K$-rays which land on the point $z_j\in\Or$.\bs

\nin{\tf Lemma 2.3. Properties of Orbit Portraits.} {\it If this
Standing Hypothesis 2.2 is satisfied, then:\ssk

\nin$(1)$ Each $A_j$ is a finite subset of $\Q/\Z$.\ssk

\nin$(2)$  For each $j$ modulo $p$, the doubling map $t\mapsto
2\,t\;({\rm mod~}\Z)$ 
carries $A_j$ bijectively onto $A_{j+1}$ preserving cyclic order
around the circle,\ssk

\nin$(3)$  All of the angles in
$A_1\cup\cdots\cup A_p$ are periodic under doubling, with a common period
$rp$, and\ssk
 
\nin$(4)$  the sets $A_1\,,\,\ldots\,,\, A_p$ are {\bit pairwise unlinked};
that is,  for each $i\ne j$ the sets $A_i$ and $A_j$
are contained in disjoint sub-intervals of $\R/\Z$.}

\ms

\begin{figure}[htp]
\cl{\psfig{figure=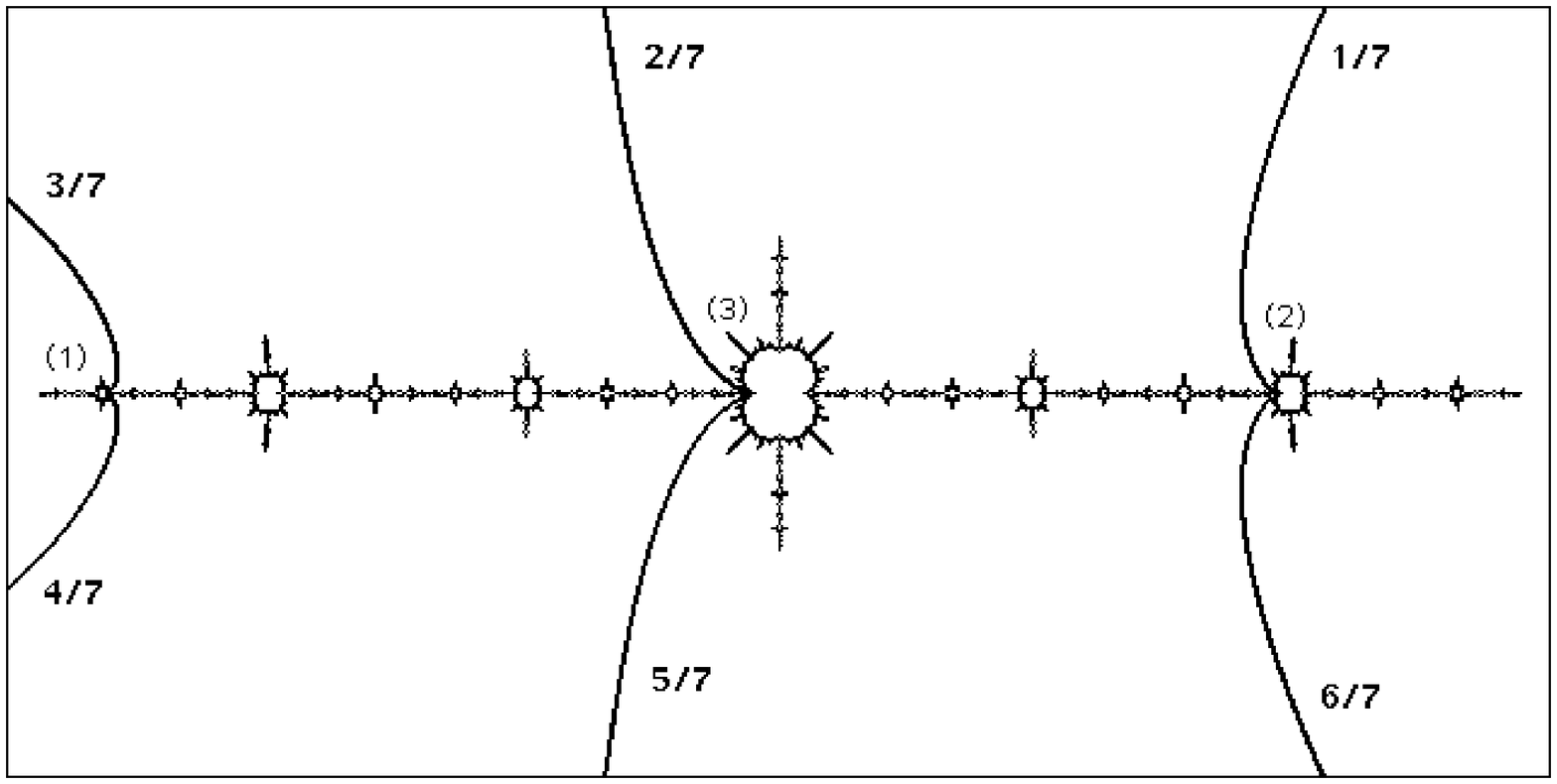,height=2in}}
\ssk
\begin{quote}{\bit Figure 6. Julia set for $z\mapsto z^2-7/4$, showing the six
$K$-rays landing on a period three parabolic orbit. Each number $(j)$
in parentheses is close to the orbit point $z_j$ (and also to $f^{\circ
 j}(0)$).\ssk}
\end{quote}

\ssk
\cl{\psfig{figure=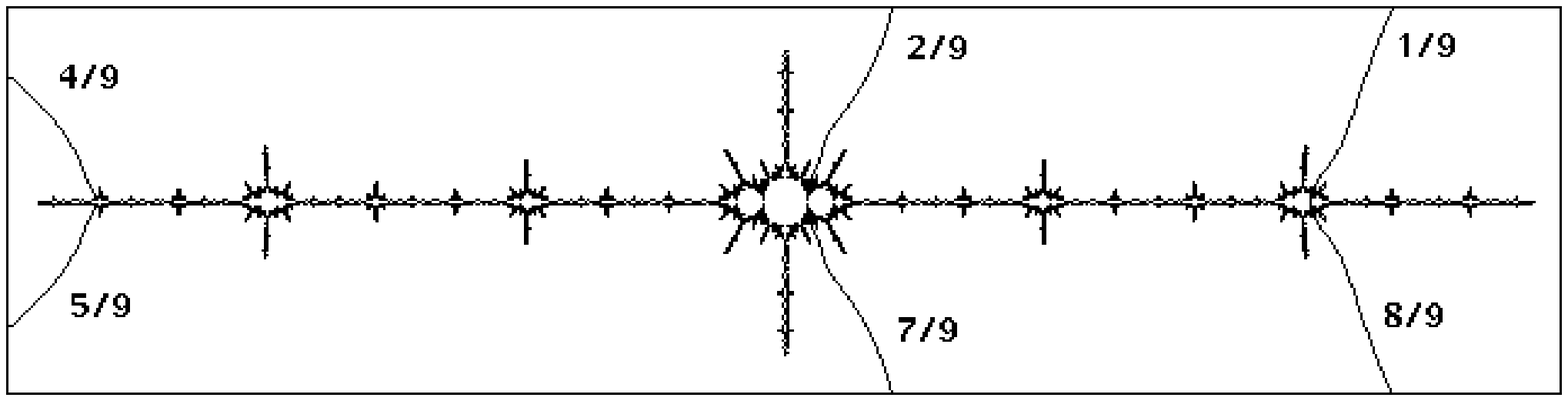,height=1.3in}}
\ssk\begin{quote}
{\bit Figure 7. Julia set for $z\mapsto z^2-1.77$, showing the six
$K$-rays landing on a period three orbit. In contrast to Figure 6, these six
rays are permuted cyclically by the map.\bs}\end{quote}

\ssk

\cl{\psfig{figure=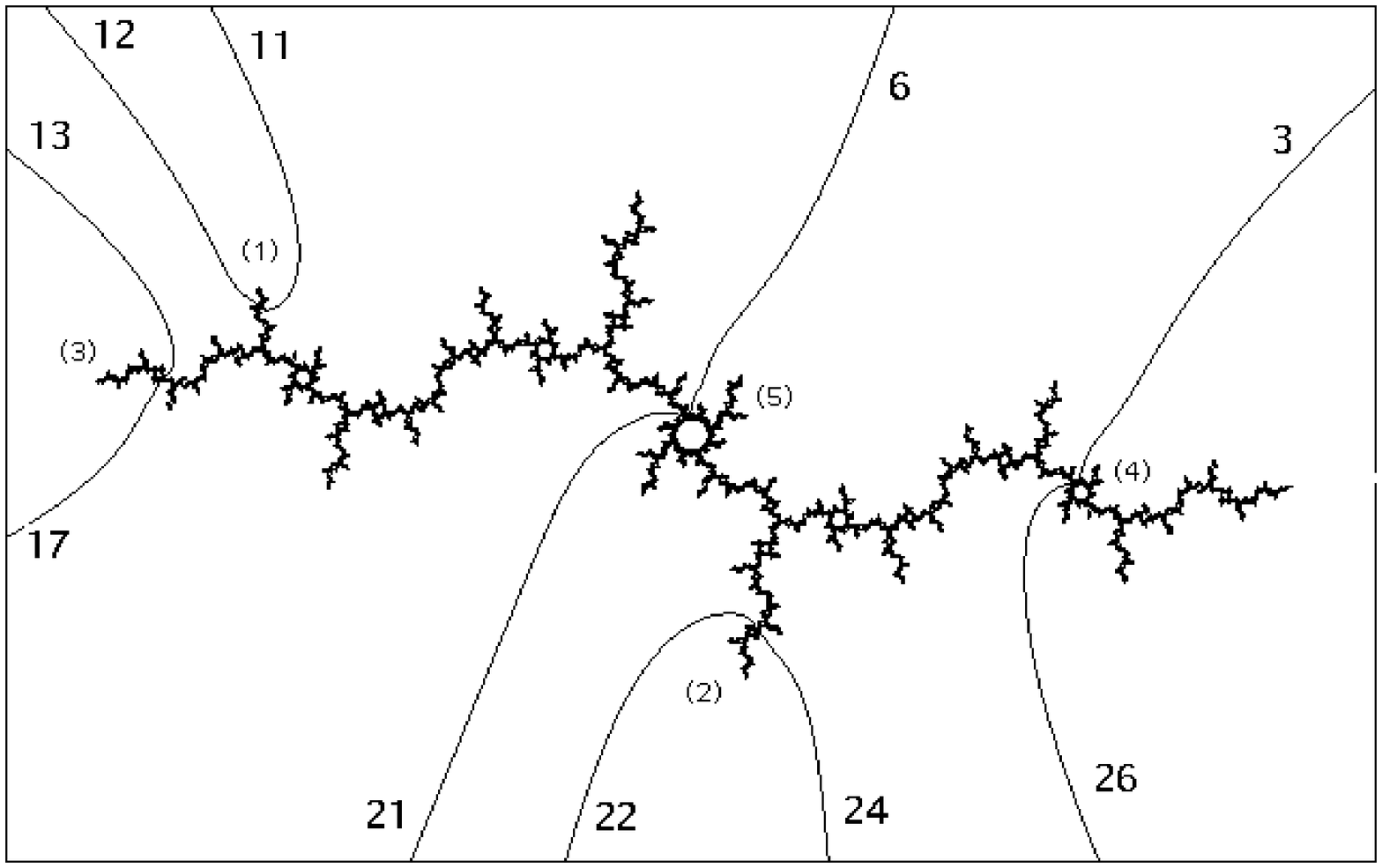,height=2.6in}}
\ssk
\begin{quote}{\bit Figure 8. Julia set $J(f_c)$ for $c=-1.2564+.3803\,i$, showing
the ten rays landing on a period 5 orbit. Here the angles are in units of
$1/31$.\ssk}\end{quote}
\end{figure}

As in \S1, the collection
$\po=\{A_1\,,\,\ldots\,,\,A_p\}$ is called the {\bit orbit portrait\/}
for the orbit $\Or$.
As examples, Figure 6 shows an orbit of period and ray period three,
with portrait
$$	\po~=~\big\{\;\{3/7\,,\,4/7\}\;,\;\{6/7\,,\,1/7\}\;,\;\{5/7\,,\,2/7\}\;
	\big\}~, $$
Figure 7 shows a period three orbit with ray period six, and with portrait
$$	\po~=~\big\{\;\{4/9\,,\,5/9\}\;,\;\{8/9\,,\,1/9\}\;,\;\{7/9\,,\,2/9\}\;
	\big\}~, $$
while Figure 8 shows an orbit of period and ray period five, with portrait
$$	\po~=~\left\{\,\left\{{11\over 31}\,,\,{12\over 31}\right\}\,,\,\left\{{22\over 31}\,,\,{24\over 31}\right\}\,,\,\left\{{13\over 31}\,,\,{17\over 31}\right\}\,,\,
\left\{{26\over 31}\,,\,{3\over 31}\right\}\,,\,\left\{{21\over 31}\,,\,{6\over 31}\right\}\,\right\}~. $$
\ms

\begin{figure}[htb]
\cl{\psfig{figure=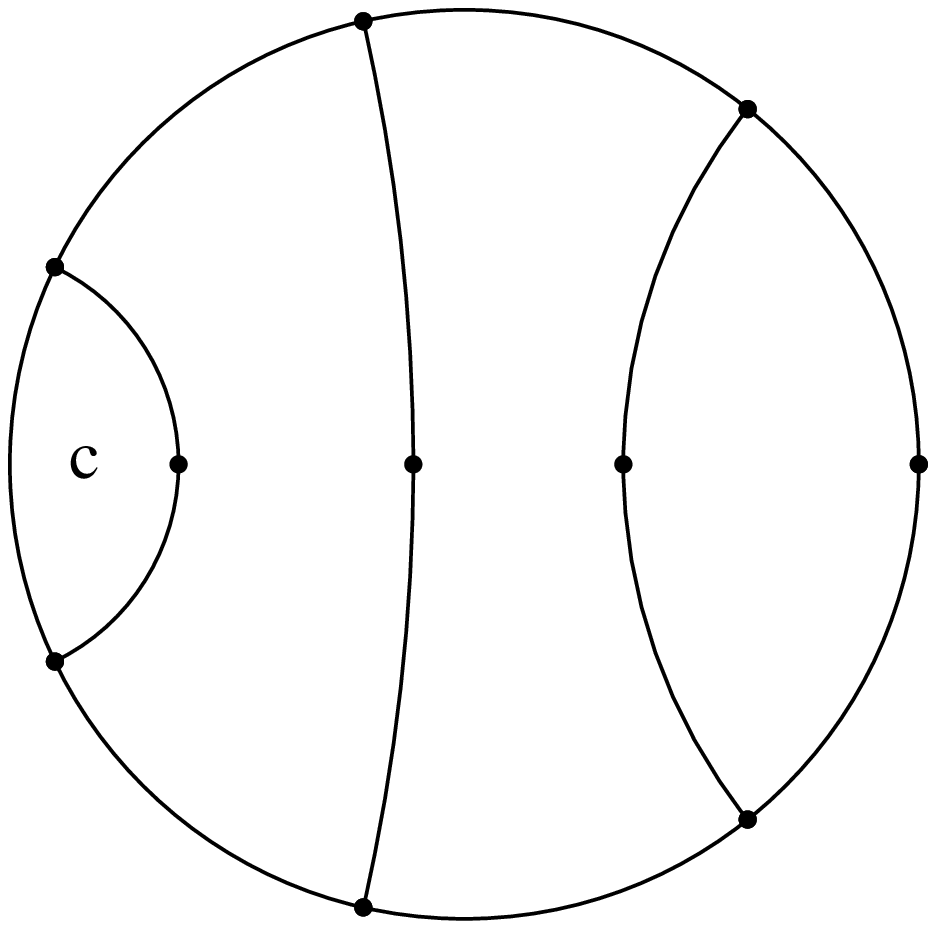,height=1.5in}\qquad
  \psfig{figure=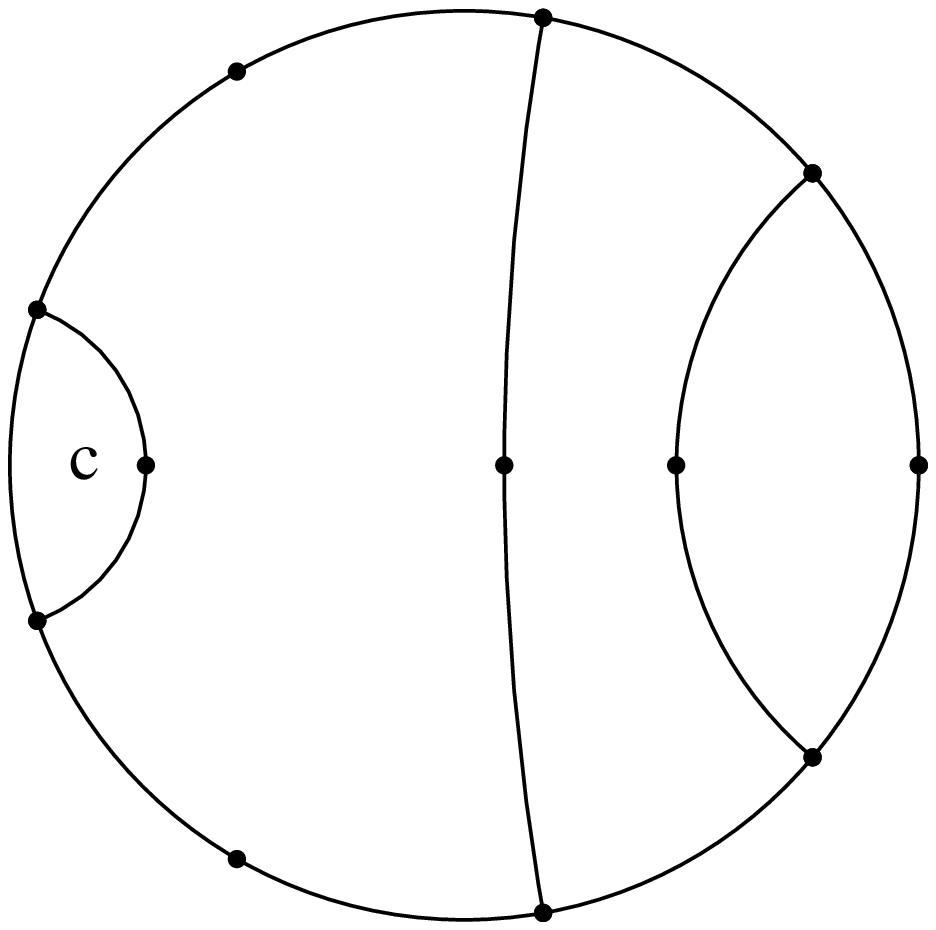,height=1.5in}\qquad
  \psfig{figure=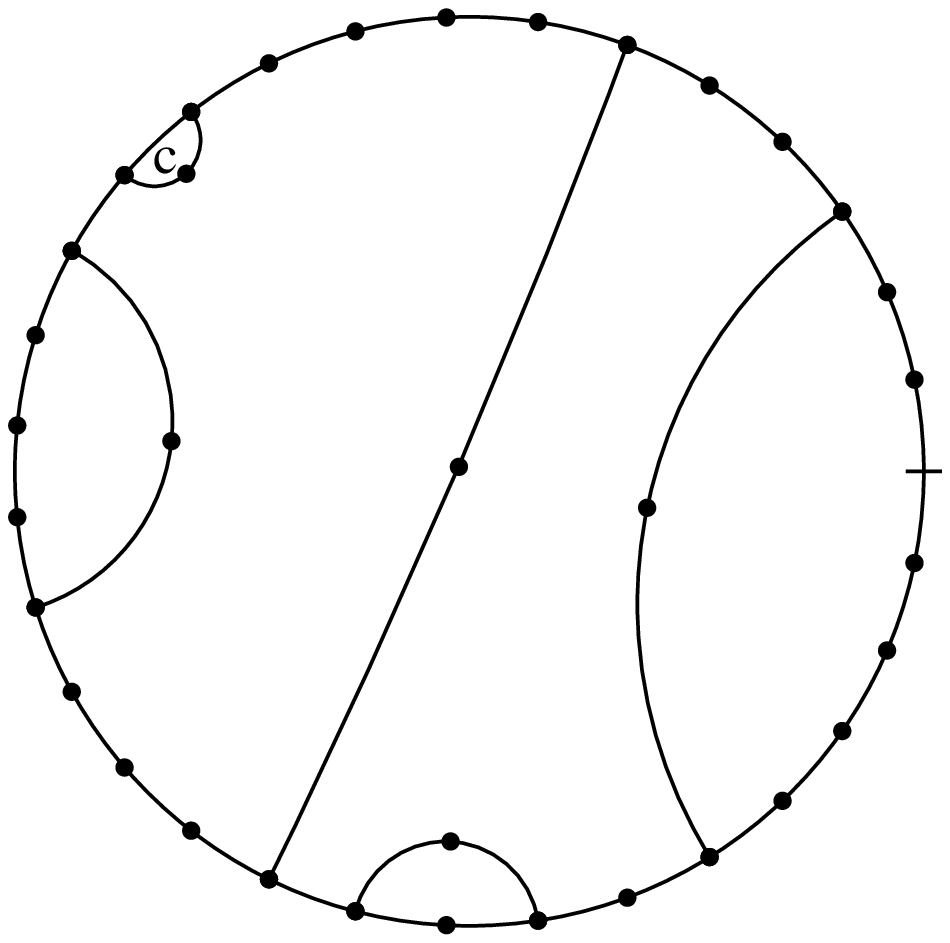,height=1.5in}\qquad }\ssk
\begin{quote}
{\bit Figure 9. Schematic diagrams associated with the orbit portraits of
Figures 6, 7, 8. The angles are in units of $1/7$, $1/9$ and
$1/31$ respectively.\par}\end{quote}
\end{figure}

{\tf Proof of 2.3.} Since some $A_i$ contains a rational number modulo
$\Z$, it follows from 2.1 that some $A_j$ contains an angle $t_0$
which is periodic under doubling. Let the period be $n\ge 1$,
so that $2^nt_0\equiv t_0~~({\rm mod~}\Z)$. 
Applying 2.1 $n$ times, we see that the mapping $\eta(t)\equiv
2^nt~~({\rm mod~}\Z)$
maps the set $A_j\subset\R/\Z$ injectively into itself, preserving cyclic
order and fixing $t_0$.
In fact we will show that every element of $A_j$ is fixed by $\eta$.
For otherwise, if $t\in A_j$ were not fixed, then
choosing suitable representatives modulo $\Z$ we would have
for example $t_0=\eta(t_0)<t<\eta(t)<t_0+1$. Since $\eta$ preserves cyclic order,
it would then follow inductively that
$$	t_0\,<\,t\,<\,\eta(t)\,<\,\eta^{\circ 2}(t)\,<\,\eta^{\circ 3}(t)
\,<\,\cdots\,<\,t_0+1~. $$
Hence the successive images of $t$ would converge to a fixed point of $\eta$.
But this
is impossible since every fixed point of $\eta$ is repelling. Thus $\eta$
fixes every point of $A_j$. But the fixed points of $\eta$ are precisely
the rational numbers of the form $i/(2^n-1)$, so it follows that $A_j$ is a
finite set of rational numbers. It follows easily that all of the $A_k$
are pointwise fixed by $\eta$. This proves (1), (2) and (3) of 2.3; and (4)
is clearly true since rays cannot cross each other.\QED\ssk

It is often convenient to compactify the complex numbers by adding a circle of
points $e^{2\pi it}\infty$ at infinity, canonically parametrized by
$t\in \R/\Z$.
Within the resulting closed topological disk $\cc$, we can form a
{\bit diagram\/} $\mathcal D$ illustrating the orbit portrait $\po$ by
drawing all of the $K$-rays joining the circle at infinity to $\Or$.
These various rays
are disjoint, except that each $z\in\Or$ is a common endpoint for exactly
$v$ of these rays.

Note that this diagram $\mathcal D$ deforms continuously, preserving its
topology, as we move the parameter point $c$, provided that the periodic
orbit $\Or$ remains repelling, and provided that the associated
$K$-rays do not run into precritical points. (Compare [GM, Appendix B].)

In fact, given $\po$, we can construct a diagram homeomorphic to
${\mathcal D}$ as follows. Start with the unit circle, and
mark all of the points $e(t)=e^{2\pi it}$ corresponding to angles $t$
in the union ${\bf A}_\po=A_1\cup\cdots\cup A_p$. Now for each $A_i$,
let $\hat z_i$ be the center of gravity of the corresponding points $e(t)$,
and join each of these points to $\hat z_i$ by a straight line segment.
It follows easily from Condition (4) that these line segments will
not cross each other. (In practice, in drawing such diagrams, we will
not usually
use straight lines and centers of gravity, but rather use some topologically
equivalent picture, fixing the boundary circle, which is easier to see.
Compare Figures 2, 5, 9.)\ssk

 It will be convenient to temporarily introduce the term {\bit formal
orbit portrait\/} for a collection $\po=\{A_1\,,\,\ldots\,,\,A_p\}$ of subsets
of $\R/\Z$ which satisfies the four conditions of 2.3, whether or not
it is actually associated with some periodic orbit. In fact we will prove the
following.

{\ssk\nin{\tf Theorem 2.4. Characterization of Orbit Portraits.} \it If $\po$
is any formal orbit portrait, then there exists a quadratic polynomial $f$
and an orbit $\Or$ for $f$ which realizes this portrait $\po$.\ms}

This will follow from Lemma 2.9 below.
To begin the proof, let us study the way in which the angle doubling
map acts on a formal orbit portrait.
As in \S1, the number of angles in each $A_j$ will be called the {\bit
valence\/} $v$ for the formal portrait $\po$. It is easy to see that any
formal portrait of valence $v=1$ can be realized by an appropriate
orbit for the map $f(z)=z^2$. Hence it suffices to study the case
$v\ge 2$. For each $A_j\in\po$
the $v$ connected components of the complement $\R/\Z\ssm A_j$ are
connected open arcs with total length $+1$. These will be called the
{\bit complementary arcs\/} for $A_j$.

{\ssk\nin{\tf Lemma 2.5. The Critical Arcs.} \it For each $A_j$
in the formal orbit portrait $\po$, all but
one of the complementary arcs is carried diffeomorphically by the angle
doubling map onto a complementary arc for $A_{j+1}$. However, the
remaining complementary arc for $A_j$ has length greater than $1/2$. Its image
under the doubling map covers one particular complementary arc for
$A_{j+1}$ twice, and every other complementary arc for $A_{j+1}$
just once.\ssk}

{\tf Definition.} This longest complementary arc will be called the {\bit
critical arc\/} for $A_j$. The arc which it covers twice under
doubling will be called the {\bit critical value arc\/} for $A_{j+1}$.
(This language will be justified in 2.9 below.)\ssk

{\tf Proof of 2.5.}
If $I\subset\R/\Z$ is a complementary arc for $A_j$
of length less
than $1/2$, then clearly the doubling map carries $I$ bijectively onto
an arc $2I$ of twice the length, bounded by two points of $A_{j+1}$.
This image arc cannot contain any other point of $A_{j+1}$, since the
doubling map from $A_j$ to $A_{j+1}$ preserves cyclic order. It follows
easily that these image arcs cannot overlap. Since we cannot fit $v$
arcs of total length $+2$ into the circle without overlap,
and since there cannot be any complementary arc of length exactly $1/2$,
it follows that there must be exactly one ``critical'' complementary
arc for $A_j$ which has length greater than $1/2$.
Suppose that it has length
$(1+\epsilon_j)/2$. Then the $v-1$ non-critical arcs for $A_j$
have total length $(1-\epsilon_j)/2$, and their images under doubling
form $v-1$ complementary arcs for $A_{j+1}$ with total length
$1-\epsilon_j$. Since the doubling map is exactly two-to-one, it follows
easily that it maps the critical arc for $A_j$ onto the entire circle,
doubly covering one ``critical value arc'' for $A_{j+1}$ which has
length $\epsilon_j$, and covering
every other complementary arc for $A_{j+1}$ just once.\QED

{\ssk\nin{\tf Lemma 2.6. The Characteristic Arc for $\po$.} \it Among the
complementary arcs for the various $A_j\in\po$, there exists a unique
arc $\I_\po$
of shortest length. This shortest arc is a critical value
arc for its $A_j$, and is contained in all of the other critical
value arcs.\ssk}

{\tf Definition.} This shortest complementary arc $\I_\po$
will be called the {\bit characteristic arc\/} for $\po$. (Compare 2.11.)
\ssk

{\tf Proof of 2.6.} There certainly exists at least one complementary arc
$\I_\po$ of minimal length $\ell$
among all of the complementary arcs for all
of the $A_j\in\po$. This $\I_\po$ must be a critical value arc,
since otherwise it would have the form $2J$ where $J$ is some complementary
arc of length $\ell/2$. Suppose then that $\I_\po$ is the critical value
arc for $A_{j+1}$, doubly covered by the critical arc $I_c$
for $A_j$. Since $\I_\po$ is minimal, it follows from 2.3(4) that this
open arc $\I_\po$ cannot contain any point of the union ${\bf A}_\po=
A_1\cup\cdots\cup A_p$. Hence its preimage under doubling also cannot contain
any point of ${\bf A}_\po$. This preimage consists of two arcs $I'$
and $I''=I'+1/2$, each of length $\ell/2$. Note that both
of these arcs
are contained in $I_c$. In fact the arc $I_c$ of length
$(1+\ell)/2$ is covered by these two open
arcs of length $\ell/2$ lying at either end, together with the
closed arc $I_c\ssm(I'\cup I'')$ of length $(1-\ell)/2$ in the middle.

Now consider any $A_k\in\po$ with $k\not\equiv j$. It follows from the
unlinking property 2.3(4) that the entire set
$A_k$ must be contained either in the arc
$(\R/\Z)\ssm I_c$ of length $(1-\ell)/2$, or in $I_c$ and hence in the
arc $I_c\ssm(I'\cup I'')$ which also has
length $(1-\ell)/2$. In either case,
it follows that the union of all non-critical arcs for $A_k$ is
contained in this same arc of length $(1-\ell)/2$, and hence that
the image of this union under doubling is contained in the arc
$$	2((\R/\Z)\ssm I_c)~=~ 2(I_c\ssm(I'\cup I''))~=~
 (\R/\Z)\ssm \I_\po $$
of length $1-\ell$. Therefore, the critical value arc for $A_{k+1}$
contains the complementary arc $\I_\po$, as required. It follows that
this minimal arc $\I_\po$ is unique. For if there were an $\I'_\po$
of the same length, then this argument would show that each of these two
must contain the other, which is impossible.\QED

{\tf Remark.} This characteristic arc never contains the angle zero.
In fact let $I_c$ be the critical arc whose image under doubling
covers $I_\po$ twice. If $0\in I_\po$, then it is not hard to see that
one endpoint of $I_c$ must lie in $\I_\po$ and the other endpoint must
lie outside, in $1/2+\I_\po$. But this is impossible by 2.3(4) and the
minimality of $I_\po$.\ssk

Recall that the union ${\bf A}_\po=A_1\cup\cdots\cup A_p$ contains $pv$
elements, each of which has period $rp$ under doubling. Hence this union
splits up into
$$	{pv\over rp}~=~{v\over r} $$
distinct cycles under doubling. If $\po$ is the portrait of a periodic orbit
$\Or$, then the ratio
$v/r$ can be described as the {\bit number of cycles of
$K$-rays\/} which land on the orbit $\Or$. As examples,
we have $v=r=3$ for Figure 1 and $v=r=2$ for Figure 7 so that there is only
one cycle under doubling, but $v=2$ and $r=1$
for Figures 6 and 8 so that there are two distinct cycles. In fact we next
show that there are at most two cycles in all cases.
\ssk

{\ssk\nin{\tf Lemma 2.7. Primitive versus Satellite.} \it Any formal orbit
portrait of valence $v>r$ must have $v=2$ and $r=1$. It follows that there
are just two posibilities:

\nin{\bf Primitive Case.} If $r=1$, so that every ray which lands on the
period $p$ orbit is mapped
to itself by $f^{\circ p}$, then at most two rays land on each orbit point.

\nin{\bf Satellite Case.}
If $r>1$, then $v=r$ so that exactly $r$ rays land on each orbit point,
and all of these rays belong to a single cyclic orbit under angle doubling.
\ssk}

This terminology will be justified in \S6. (Compare Figure 12.)\ssk


{\tf Proof of 2.7.} Suppose that $v>r$ and $v\ge 3$.
Let $\I_\po$ be the characteristic arc. We suppose that
$\I_\po$ is the critical value arc in the complement of $A_1$.
Let $I_-$ the complementary arc for $A_1$ which is
just to the left of $I_\po$ and let $I_+$
be the complementary arc just to the right of $I_\po$.
To fix our ideas, suppose that $I_-$ has length $\ell(I_-)\ge\ell (I_+)$.
Since $I_+$ is not the critical value arc for $A_1$, we see, arguing
as in 2.6, that it must
be the image under iterated doubling of the critical value arc $I'$
for some $A_j$. That is, we have $I_+=2^m I'$ for some $m\ge 1$.
Hence $\ell(I')<\ell(I_+)$.

The hypothesis that $v>r$ implies that the two endpoints of $\I_\po$ belong
to different cycles under doubling. Thus the left endpoints of $I'$ and
$\I_\po$ belong to distinct cycles, hence $I'\ne \I_\po$. Therefore, by 2.6,
$I'$ strictly contains $\I_\po$. This arc $I'$ cannot strictly
contain the neighboring arc $I_+$, since it is shorter than $I_+$.
Hence it must have an endpoint in $I_+$, and therefore, by 2.3(4), it must
have both endpoints in $I_+$. But this implies that $I'$ contains
$I_-$, which is impossible since $\ell(I')<\ell(I_+)\le\ell(I_-)$.
Thus, if $v>r$ it follows that $v\le 2$, hence $r=1$ and $v=2$,
as asserted.\QED\ssk

{\ssk\nin{\tf Lemma 2.8. Two Rays determine $\po$.}
\it Let $\po=\{A_1\,,\,\ldots\,,\, A_p\}$ be a formal orbit portrait of
valence $v\ge 2$, and let $\I_\po=
(t_-\,,\,t_+)$ be its characteristic arc, as described above. Then
a quadratic polynomial $f_c$ has
an orbit with portrait $\po$ if and only if the two $K$-rays
with angles $t_-$ and $t_+$ for the
filled Julia set of $f_c$ land at a common point.\ssk}

{\tf Proof.} If $f_c$ has an orbit with portrait $\po$, this is true by
definition. Conversely, if these rays land at a common point $z_1$, then
the orbit of $z_1$ is certainly periodic. Let $\po'$ be the portrait for
this actual orbit. We will denote its period by $p'$, its valence by $v'$,
and so on. Note that the ray period $rp$ is equal to $r'p'$,
the common period of the angles $t_-$ and $t_+$ under doubling.

{\bf Primitive Case.} Suppose that $r=1$ so that $v/r=2$, and so
that each of these angles $t_\pm$
has period exactly $p$ under doubling. If $p'<p$ hence $r'>1$, then
it would follow from 2.7 applied to the portrait $\po'$ that $t_-$ and
$t_+$ must belong to the same cycle under doubling, contradicting the
hypothesis that $v/r=2$.

{\bf Satellite Case.} If $r>1$ hence $v=r$, then $t_-$ and $t_+$
do belong to the same cycle under doubling, say $2^kt_-\equiv t_+~({\rm mod~}\Z)$.
Clearly it follows that $r'>1$ hence $v'=r'$. Furthermore, it
follows easily that multiplication by $2^k$ acts transitively on
$A_1$, and hence that all of the rays $\ra_t^K$ with $t\in A_1$ land at
the same point $z_1$. In other words $A_1\subset A'_1$.
This implies that $r\le r'$ hence
$p\ge p'$. If $p$ were strictly greater than $p'$, then it would
follow that $A_{1+p'}$ is also
contained in $A'_1$. But the two sets $A_1$ and $A_{1+p'}$ are unlinked
in $\R/\Z$. Hence there is no way that multiplication by $2^p$ can
act non-trivially on $A_1\cup A_{1+p'}$ carrying each of these two sets
into itself and preserving cyclic order on their union. This contradiction
implies that $A_1=A'_1$ and $p=p'$, and hence that $\po=\po'$, as required.
\QED

Now let $c$ be some parameter value outside the Mandelbrot set. Then,
following Douady and Hubbard, the point $c$, either in the dynamic plane or
in the parameter plane,
lies on a unique external ray, with the same well defined angle $t(c)\in\R/\Z$
in either case. (Compare Appendix A.)

{\ssk\nin{\tf Lemma 2.9. Outside the Mandelbrot Set.} \it Let
$\po=\{A_1\,,\,\ldots\,,\,A_p\}$ be a formal orbit portrait with
characteristic arc $\I_\po$, and let
$c$ be a parameter value outside of the Mandelbrot set.
Then the map $f_c(z)=z^2+c$ admits a periodic
orbit with portrait $\po$ if and only if the {\bit external angle}
$t(c)$ belongs to this open arc $\I_\po$.\ssk}

{\tf Proof.} The two dynamic rays $\ra_{t(c)/2}^K$
and $\ra_{(1+t(c))/2}^K$ meet at the critical point $0$, and together
cut the dynamic plane into two halves. Furthermore,
 every point of the Julia set $\partial K=K$
is uniquely determined by its symbol sequence with respect to this partition.
Correspondingly, the two diametrically opposite points $t(c)/2$
and $(1+t(c))/2$ on the circle $\R/\Z$ cut the circle into two semicircles,
and almost every point $t\in\R/\Z$ has a well defined symbol sequence with
respect to this partition under the doubling map. Two rays $\ra_t^K$ and
$\ra_u^K$ land at a common point of $K$ if and only if the external angles
$t$ and $u$ have the same symbol sequence.

First suppose that the angle $t(c)$ lies in the characteristic arc
$\I_\po$. Then, with notation as in the proof of 2.6, the two points
$t(c)/2$ and $(1+t(c))/2$ lie in the two components $I'$ and $I''$
of the preimage of $\I_\po$. For every $A_j\in\po$, all of the
points of $A_j$ lie in a single component of $\R/\Z\ssm(I'\cup I'')$.
Hence the rays $\ra_t^K$ with $t\in A_j$ land at a common point $z_j
\in K$. It follows from 2.8 that these points lie in an orbit with portrait
$\po$, as required.

On the other hand, if $t(c)$ lies outside of $\overline \I_\po$, then it
is easy to check that the two endpoints of $\I_\po$ are separated by the points
$t(c)/2$ and $(1+t(c))/2$. Hence these two endpoints, both belonging to
$A_1\in\po$, land at different points of $K$. Hence $f_c$ has no orbit
with portrait $\po$.

Finally, in the limiting case where $t(c)$ is precisely equal to one of the
two endpoints $t_\pm$ of $\I_\po$, since these angles are periodic under
doubling, it follows that the ray $\ra_{t_\pm}^K$ passes through a precritical
point, and hence does not have any well defined landing point in $K$.
This completes the proof of 2.9.\QED\ssk

Evidently the Realization Theorem 2.4 is an immediate corollary. Since we have
proved 2.4, we can now forget about the distinction between ``formal''
orbit portraits and portraits which are actually realized. We can describe
further properties of portraits and their associated diagrams as follows.

{\tf Definition 2.10.} Suppose that we start
with any periodic orbit $\Or$ with valence $v\ge 2$ and period $p\ge 1$,
and fix some point $z_i\in\Or$. As in \S1, the $v$ rays landing at $z_i$
cut the dynamic plane $\C$ up into $v$ open subsets
which we call the {\bit sectors\/} based at $z_i$. Evidently there is a
one-to-one correspondence between sectors based at $z_i$ and complementary
arcs for the corresponding set of angles $A_i\subset\R/\Z$,
characterized by the property that $\ra_t^K$ is contained in the open
sector $S$ if and only if $t$ is contained in the corresponding
complementary arc. By definition,
the {\bit angular size\/} $\alpha(S)>0$ of a sector is the length
of the corresponding complementary arc, which we can think of as its
``boundary at infinity''. It follows that
$\sum_S\,\alpha(S)=1$, where the sum extends over the $v$ sectors
based at some fixed $z_i\in\Or$.

{\tf Remark.} The angular size of a sector has nothing to do with the
angle between
the rays at their common landing point, which is often not even defined.

Altogether there are $pv$ rays landing at the various
points of the orbit $\Or$. Together these rays cut the plane up into
$pv-p+1$ connected components. The closures of these components will
be called the pieces of the {\bit preliminary puzzle\/}
associated with the diagram $\mathcal D$ or the associated portrait $\po$.
Note that every closed sector $\overline S$ can be expressed as
a union of preliminary puzzle pieces, and that every preliminary puzzle
piece is equal to the intersection of the closed sectors containing it. This
construction will be modified and developed further in Sections 7 and 8.

For every point $z_i$ of the orbit, note that just one of the $v$
sectors based at $z_i$ contains the critical point $0$. We will call
this the {\bit critical sector\/} at $z_i$, while the others will be called
the {\bit non-critical sectors\/} at $z_i$. Another noteworthy sector at
$z_i$ (not necessarily distinct from the critical sector)
is the {\bit critical value sector\/}, which contains $f(0)=c$.\ms

{\nin{\tf Lemma 2.11. Properties of Sectors.} \it The diagram ${\mathcal D}\subset
\cc$ associated with any orbit $\Or$ of
valence $v\ge 2$ has the following properties:\ssk

\nin $(a)$ For each $z_i\in\Or$, the critical sector at $z_i$
has angular size strictly greater than $1/2$.
It follows that the $v-1$ non-critical sectors at
$z_i$ have total angular size less than $1/2$.\ssk

\nin $(b)$ The map $f$ carries
a small neighborhood of $z_i$ diffeomorphically onto
a small neighborhood of $z_{i+1}=f(z_i)$, carrying each sector
based at $z_i$
locally onto a sector based at $z_{i+1}$, and preserving the cyclic
order of these sectors around their base point. The critical sector
at $z_i$ always maps locally, near $z_i$,
onto the critical value sector based at $z_{i+1}$.\ssk

\nin $(c)$ Globally, each non-critical sector $S$ at $z_i$
is mapped homeomorphically by $f$ onto a sector
$f(S)$ based at $z_{i+1}$, with angular size given by
$\alpha(f(S))~=~2\,\alpha(S)$.
However, the critical sector at $z_i$ maps so as to cover the entire plane,
covering the critical value sector at $z_{i+1}$ twice with a ramification
point at $0\mapsto c$, and covering every other sector just once.\ssk

\nin $(d)$ Among all of the $pv$ sectors based at the various points
of $\Or$, there is a unique sector of smallest angular size, corresponding
to the characteristic arc $\I_\po$. This smallest sector contains the
critical value, and does not contain any other sector.\ssk}

(As usual, the index $i$ is to be construed as an integer modulo $p$.)
The proof, based on 2.6 and the fact that $f$ is exactly two-to-one except
at its critical point, is straightforward and will be left to the reader.
Evidently Theorem 1.1 follows.\QED


Now let us take a closer look at the dynamics of the diagram $\mathcal D$ or of
the associated portrait $\po$. The iterated map $f^{\circ p}$ fixes each
point $z_i\in\Or$, permuting the various rays which land on $z_i$ but
preserving their cyclic order.
Equivalently, the $p$-fold iterate of the doubling map carries each finite
set $A_i\subset\Q/\Z$ onto itself by a bijection which preserves the cyclic
order. For any fixed $i$ mod $p$, we can number the angles in $A_i$
as $0\le t^{(1)}<t^{(2)}<\cdots<t^{(v)}<1$. It then follows that
$$	2^p\,t^{(j)}~\equiv~t^{(j+k)}\quad({\rm mod~}\Z)~, $$
taking superscripts modulo $v$, where $k$ is some fixed residue
class modulo $v$.

\nin{\tf Definition 2.12.} The ratio $k/v~({\rm mod~}\Z)$ is called the
combinatorial {\bit rotation number\/} of our orbit portrait. It is easy to
check that this rotation number does not depend on the choice of orbit point
$z_i$.
Let $d$ be the greatest common divisor of $v$ and $k$. The we can
express the rotation number as a fraction $q/r$ in lowest terms,
where $k=qd$ and $v=rd$. (In the special case of rotation number zero,
we take $q=0$ and $r=1$.)

In all cases, note that the denominator
$r\ge 1$ is equal to the period of the angles $t^{(j)}\in A_i$
under the mapping $t\mapsto 2^p\,t~({\rm mod~}\Z)$ from $A_i$ to itself.
It follows easily that the period of $t^{(j)}$ under angle doubling is equal
to the product $rp$. Thus this definition of $r$ as the denominator
of the rotation number is compatible with our earlier notation $rp$
for the ray period.\ms

{\tf Notation Summary.} Since we have been accumulating quite a bit
of notation, here is a brief summary:

{\bit Orbit period\/} $p$: the number
of distinct element in our orbit $\Or$,

{\bit Ray period\/} $rp$:
the period of each angle $t\in A_1\cup\cdots\cup A_p$ under doubling.

{\bit Rotation number\/} $q/r$: describes the action of
multiplication by $2^p$ on each set $A_i$.

{\bit Valence\/} $v$: number of angles in each $A_i$, for a
total of $pv$ angles altogether.

{\bit Cycle number\/} $v/r$: the number of disjoint cycles of size
$rp$ in the union $A_1\cup\cdots\cup A_p$.\ms

According to 2.7, this cycle number is always equal to $1$ for a
satellite portrait, and is at most 2 in all cases. Thus,
in the case $v\ge 2$
there are just two possibilities as follows:\ssk

{\ssk\nin{\bf Primitive Case.} The rotation number is zero. There are $v=2$ rays
landing at each orbit point, for a total of $2p$ rays. These split up into
two cycles of $p$ rays each under doubling.\ssk}

{\ssk\nin{\bf Satellite Case.} The rotation number is $q/r\not\equiv 0$. There are
$v=r$ rays landing at each orbit point, for a total of $pv=rp$ rays
altogether. These $rp$ rays are permuted cyclically under angle doubling,
so that the number of cycles is $v/r=1$.\ssk}

As examples, Figures 6, 8 illustrate primitive portraits with rotation
number zero, while Figures 1, 7 show satellite portraits with rotation
number $1/3$ and $1/2$. We will see in
\S6 that primitive portraits correspond to primitive hyperbolic
components in the Mandelbrot set, that is, to those with a cusp point.



\section{Parameter Rays.} This section will prove
the following preliminary version of Theorem 1.2.

Let $\po$ be any orbit portrait of valence $v\ge 2$, and let
$\I_\po=(t_-\,,\,t_+)$ be its characteristic arc, where
$0<t_-<t_+<1$. If the quadratic
polynomial $f_c=z^2+c$ has an orbit $\Or$ with portrait $\po$,
recall that the two dynamic rays $\ra_{t_-}^K$ and $\ra_{t_+}^K$ for $f_c$
land at a common orbit point, and together bound a sector $S_1$ which has
minimal angular size among all of the sectors based at points of the
orbit $\Or$. This $S_1$ can also be characterized as the smallest of these
sectors which contains the critical value $c$. (Compare Lemmas
2.6, 2.9, 2.11.)\ssk

{\nin{\tf Theorem 3.1. Parameter Rays and the Wake.}
\it The two parameter rays $\ra_{t_-}^M$ and $\ra_{t_+}^M$
with these same angles
land at a common parabolic point in the Mandelbrot set. Furthermore,
these two rays, together
with their common landing point, cut the parameter plane into two open subsets
$W_\po$ and $\C\ssm\overline W_\po$ with the following property:
The quadratic map $f_c$ has a repelling orbit with portrait $\po$
if and only if $c\in W_\po$.\ssk}

{\tf Proof.} Let ${\bf A}_\po=A_1\cup\cdots\cup A_p$ be the set of all angles
for the orbit portrait $\po$, and let
$n=rp$ be the common period of these angles under doubling.
The set $F_n\subset M$ of possibly exceptional parameter values
will consist of those $c$ for which $f_c^{\circ n}$ has a fixed point of
multiplier $+1$. Since $F_n\subset\C$ is an algebraic variety and is not
the entire complex plane, it is necessarily a finite set.
As noted in [GM], if $c$
belongs to the Mandelbrot set but $c\not\in F_n$, then the various
dynamic rays $\ra_t^{K(f_c)}$ with $t\in {\bf A}_\po$ all
land on repelling periodic points, and the pattern of which of these
rays land at a common point remains stable under
perturbation of $c$ throughout some open neighborhood within parameter space.

Now suppose that $c$ lies outside of the Mandelbrot set. Then $c$,
considered as a point in parameter space,
belongs to some uniquely defined parameter ray $\ra_{t(c)}^M$, and considered
as a point in the dynamic plane for $f_c$, belongs to the dynamic ray
$\ra_{t(c)}^K$ with this same angle. In this case,
a dynamic ray $\ra_t^K$ for $f_c$ has a well defined landing point
in $K=K(f_c)$ if and only if the forward orbit $\{2t\,,\,4t\,,\,8t\,,\,
\ldots\}$ under doubling does not contain this critical value angle $t(c)$.
Since the angles in ${\bf A}_\po$ are periodic, it follows
that all of the dynamic rays $\ra_t^K$ with $t\in{\bf A}_\po$ have well
defined landing points in $K$ if and only if $t(c)\not\in{\bf A}_\po$.

Let $t\in{\bf A}_\po$ and let $c_0\in M$ be any accumulation point
for the parameter ray $\ra_t^M$. Since every neighborhood of $c_0$
contains parameter values $c\in\ra_t^M$ for which the dynamic
ray $\ra_t^{K(f_c)}$
does not land, it follows that $c_0$ must belong to $F_n$.
Thus every accumulation point for $\ra_t^M$ belongs to the finite set $F_n$.
Since the set of all accumulation points of a ray is connected, this
proves that $\ra_t^M$ must actually land at a single point of $F_n$.

These parameter rays $\ra_t^M$ with $t\in{\bf A}_\po$,
together with the points of $F_n$, cut the complex
parameter plane up into finitely many open sets $U_i$, and the pattern
of which of the corresponding dynamic rays
$\ra_t^K$ with $t\in{\bf A}_\po$
land at a common periodic point remains fixed as $c$ varies
through any $U_i$. Since every $U_i$ is unbounded, it follows from
Lemma 2.9 that for $c\in U_i$ the map $f_c$
has an orbit with portrait $\po$ if and only if $U_i$
is that open set which contains the points in $\C\ssm M$ with external
angle $t(c)$ in $(t_-\,,\,t_+)$. It follows that the two rays $\ra_{t_-}^M$
and $\ra_{t_+}^M$ must land at a common point of $F_n$, so as to separate
the parameter plane. For otherwise, if they had different landing points,
the connected set $U_i$ containing points with external angles in
$(t_-\,,\,t_+)$ would also contain points with other external angles,
which is impossible.

Define the {\bit root point\/} ${\bf r}_\po\in M$ to be this common
landing point, and define the {\bit
wake\/} $W_\po$ to be that connected component of
$\C\ssm(\ra^M_{t_-}\cup\ra^M_{t_+}\cup {\bf r}_\po)$
which does not contain $0$. For $c\in W_\po\ssm F_n$, it follows from
the discussion above that $f_c$ does have a repelling orbit with portrait
$\po$, while for $c\in\C\ssm (W_\po\cup F_n)$ it follows that $f_c$ does
not have any repelling orbit with portrait $\po$. Thus, to complete the
proof of 3.1, we need only consider those $f_c$ with $c$ in the finite
set $F_n$.

First suppose that some point $c_0\in F_n\ssm W_\po$ had a repelling orbit
with portrait $\po$. Then any nearby parameter value would have a nearby
repelling orbit with the same landing pattern for rays with angles in
${\bf A}_\po$. A priori it might seem possible that some extra ray, perhaps
one landing on a parabolic orbit for $f_{c_0}$, might land on this same
repelling orbit after perturbation. (Compare [GM, Fig.~12].) However,
this is ruled out by 2.8. Hence all nearby parameter
values must belong to $W_\po$, which is impossible.

Now consider parameter points $c\in W_\po$. I am indebted to
Tan Lei for pointing out the following very elegant argument due to
Peter Ha{\"\i}ssinsky which replaces
my own more complicated reasoning. As noted above, for every $c\in W_\po$
the rays $\ra_{t_-}^{K(f_c)}$ and $\ra_{t_+}^{K(f_c)}$ land
at well defined periodic points of the Julia set $J(f_c)$.
Let $z=z(m,t_\pm,c)$ be the unique point on the ray
$\ra_{t_\pm}^{K(f_c)}$ which has potential $G(z)=1/m$. The functions
$c\mapsto z(m,t_\pm,c)$ are evidently holomorphic;\break in fact
$z(m,t,c)=\phi_c^{-1}\Big(\exp(2\pi it+1/m)\Big)$, where the function
$$\phi_c:\C\ssm K\buildrel\cong\over\longrightarrow\C
 \ssm\overline\D\qquad\text{can be defined locally as}
 \qquad\phi_c(z)~=~\lim_{k\to\infty}\root{2^k}\of {f_c^{\circ k}(z)} $$
(choosing appropriate branches of the iterated square root) and hence
is holomorphic as a function of both variables.
These maps $c\mapsto z(m,t_-,c)$ or $c\mapsto z(m,t_+,c)$
form a normal family throughout $W_\po$, since they miss the three points
$0\,,\,c$ and $\infty$. Choosing
a convergent subsequence and passing to the limit as
$m\to\infty$, we see that the landing points also depend holomorphically on
$c$. Since the landing points for $\ra_{t_-}^{K(f_c)}$ and $\ra_{t_+}^{K(f_c)}$
coincide for $c\in W_\po\ssm F_n$, it follows by continuity that they
coincide for all $c\in W_\po$. It follows also that the multiplier $\lambda(c)$
of this common landing point depends holomorphically on the parameter
$c\in W_\po$, with $|\lambda(c)|\ge 1$ since the landing point must be
repelling or parabolic. But the absolute value of a non-zero holomorphic
function cannot have a local minimum unless it is constant,
so all of these landing points must actually be
repelling. Together with 2.8, this completes the proof of 3.1.\QED\smallskip

We will deal with parabolic orbits with portrait $\po$ in the next two
sections.\bigskip


\section{Near Parabolic Maps.} Let $\hat c$ be a parabolic point
in parameter space. This section will study the dynamic
behavior of the quadratic map
$f_c$ for $c$ in a neighborhood of $\hat c$. (Compare [DH2, \S14(CH)], [Sh2].)

Let $\Or$ be the parabolic orbit for $f_{\hat c}$ with period $p\ge 1$
and with representative point $\hat z$. Then the multiplier $\hat\lambda =
(f_{\hat c}^{\circ p})'(\hat z)$ is
a primitive $r$-th root of unity for some $r\ge 1$. Let $\po$
be the associated orbit portrait, with ray period $rp\ge p$.
We will first prove the following.

{\ssk\nin{\tf Theorem 4.1. Deformation Preserving the Orbit Portrait.} \it
There exists a smooth path in parameter space ending at the parabolic point
$\hat c$ and consisting of parameter values $c$ with the following
property: The associated
map $f_c$ has both a repelling orbit of period $p$ and
an attracting orbit of period $rp$. Furthermore, this repelling orbit
has portrait $\po$, and lies on the boundary of the immediate basin for
the attracting orbit. As $c$ tends to $\hat c$, these two orbits both
converge towards the original parabolic orbit $\Or$.\ms}

\begin{figure}[htb]
\centerline{\psfig{figure=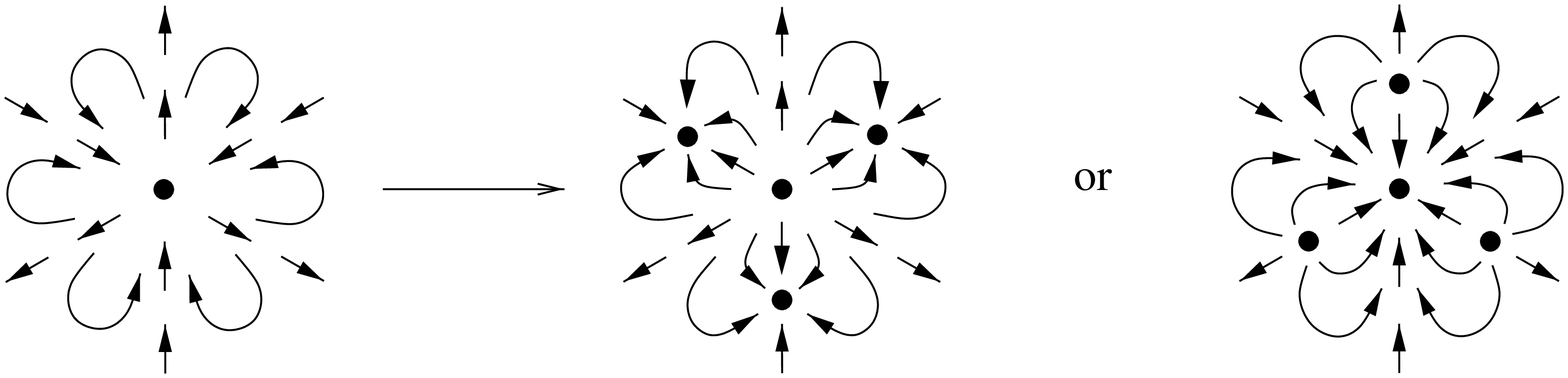,height=1.1in}}\smallskip
{\QP\bit Figure 10 (courtesy of S. Zakeri). The left sketch shows a parabolic fixed point with $r=3$,
the middle shows the modified version with an attracting orbit of period 3,
and the right shows a modified version with an attracting fixed
point. Here the arrows indicate the action of $f^{\circ 3}$.\par}
\end{figure}

(Compare Figure 10, middle.) The proof will depend on the following.\ssk

{\nin{\tf Lemma 4.2. Convenient Coordinates.} \it For any complex number
$\lambda$ close to $\hat\lambda$ there exists at least one
parameter value $c$ close to $\hat c$ and point $z_\lambda$ close to
$\hat z$ so that $z_\lambda$ is a periodic point for the map $f_c$
with period $p$ and with multiplier $\lambda$. Furthermore
there is a local holomorphic change of coordinate $z=\phi_\lambda(w)$ with
$z_\lambda=\phi_\lambda(0)$ so that the map $F=F_\lambda=\phi_\lambda^{-1}
\circ f_c^{\circ p}\circ\phi_\lambda$ takes the form
$$	F(w)~=~\lambda\,w\,  +\,R(\lambda, w) $$
for $w$ near zero, and so that its $r$-th iterate takes the form
$$	F^{\circ r}(w)~=~
\phi_\lambda^{-1}\circ f_c^{\circ rp}\circ\phi_\lambda~=~
\lambda^{\,r}w\,\big(1+w^r+R'(\lambda,w)\big)
 ~, \eqno (2)$$
where the remainder terms $R$ and $R'$
satisfy $~|R|\,,\,|R'|\le {\rm constant}\,|w|^{r+1}$
uniformly for $\lambda$ in some neighborhood of $\hat \lambda$ and for
$w$ in some neighborhood of zero.\ssk}

(In 4.5, we will sharpen this statement by showing that
the phrase ``at least one'' in 4.2 can be replaced by ``exactly
one''.)\ssk

{\tf Proof of 4.2 in the Primitive Case.} First suppose that $\po$ is a
primitive portrait, so that the multiplier
$\big(f_{\hat c}^{\circ p}\big)'(\hat z)$
is equal to $+1$ for $\hat z\in\Or$, with $r=1$.
In this case, $\hat z$ is a fixed point of multiplicity two
for the iterate $f_{\hat c}^{\circ p}$,
and splits into two nearby fixed points under
perturbation. (It cannot have a higher multiplicity, since a fixed point
of multiplicity $\mu>2$ would have $\mu-1\ge 2$ attracting Leau-Fatou petals,
each with at least one critical point in its basin, which is impossible for
a quadratic map.) As $c$
traverses a small loop around $\hat c$, these two fixed
points a priori may be (and in practice always will be) interchanged.
However, if we loop twice around $\hat c$, then each of these fixed points
must return to its original position. Thus, if we introduce a new parameter
$u$ by the equation $c=\hat c+u^2$, then we can choose these fixed
points as holomorphic functions, $z_\iota=z_\iota(u)$ for $\iota=1,2$,
with $z_1(0)=z_2(0)=\hat z$. Evidently the $u$-plane is a
two-fold branched cover of the $c$-parameter plane.
Let $\lambda_\iota(u)=\big(f_c^{\circ p}\big)'\big(z_\iota(u)\big)$ be the
multiplier for the orbit of $z_\iota$, and note that
$\lambda_1(0)=\lambda_2(0)=1$. Since the holomorphic function $u\mapsto
\lambda_1(u)$ cannot be constant, it takes on all values close to $+1$
as $u$ varies through a neighborhood of $0$.

Expanding the function $f_c^{\circ p}$ as a power series about its
fixed point $z_1$, we obtain
$$f_c^{\circ p}(z_1(u)+h)\,-\,z_1(u)\,~=~\,
 \lambda_1(u)\,h\,+\,a(u)\,h^2 \,+\,({\rm higher~terms~in~} h) \eqno (3)$$
for $h$ and $u$ close to zero, where $c=\hat c+u^2$. Here the
coefficient $a(u)$ is also a holomorphic function of $u$, with
$a(0)\ne 0$ since the fixed
point multiplicity is two. It follows that $a(u)\ne 0$ for $u$
sufficiently small. Denoting the expression (3) by $g_u(h)$, and
replacing the variable $h=z-z_1$ by $w=\alpha_u\,h$
where $\alpha_u=a(u)/\lambda_1(u)$, we see easily that the function
$$	F_u(w)~=~{\alpha_u}\;g_u\big( w/\alpha_u\big) $$
has the required form (2). \QED\ssk

{\tf Proof of 4.2 in the Satellite Case.} We now suppose that
$\hat\lambda$ is a primitive $r$-th root of unity, with $r>1$.
Then we can solve for the period $p$ point
$z=z(c)$ as a holomorphic function of $c$ for $c$ in some neighborhood of
$\hat c$, with $z(\hat c)=\hat z$.
Hence the multiplier $\lambda(c)=(f_c^{\circ p})'(z(c))$
will also be a holomorphic function of $c$, taking the value $\hat\lambda
\in\root r\of 1$ when $c=\hat c$. Similarly $\lambda(c)^r$
is a holomorphic function, taking the value $\lambda(\hat c)^r=1$ when
$c=\hat c$. This function $\lambda(c)^r$ clearly
cannot be constant, so it takes all values
close to $+1$ as $c$ varies through a neighborhood of $\hat c$.

We will construct a sequence of holomorphic changes of variable which conjugate
the map $z\mapsto f_c^{\circ p}(z)$ in a neighborhood of $z=z(c)$ to maps
$h\mapsto g_{c,k}(h)$ in a neighborhood of $h=0$, where $1\le k\le r$,
so that
$$	g_{c,k}(h)~=~\lambda(c) h\big(1+a_k(c)h^k+({\rm higher~terms~in~}h)
	\big) $$
for some constant $a_k(c)$. Here $c$ can be any point in some neighborhood
of $\hat c$. To begin the construction, let
$$	g_{c,1}(h)~=~f^{\circ p}(z(c)+h)\,-\,z(c)~.$$
This certainly has the required properties. Now inductively set
$$	g_{c,k+1}(h)~=~\phi^{-1}\circ g_{c,k}\circ\phi(h)\qquad{\rm where}
	\qquad\phi(h)~=~h+bh^{k+1} $$
for $1\le k<r$.
We claim that the constant $b=b(c)$ can be uniquely chosen so that $g_{c,k+1}$
will have the required form. In fact a brief computation shows that
$$	g_{c,k+1}(h)~=~\lambda h\big(1+(a+b-\lambda^kb)h^k+({\rm higher~
terms})~.$$
But $\lambda^k\ne 1$ since $\lambda$ is close to $\hat\lambda$, which is a
primitive $r$-th root of unity with $1\le k<r$. Hence there is a unique
choice of $b$ so that $a+b-\lambda^kb~=~0$, as required.

In particular, pushing this argument as far as possible, we can take
$k=r$ and replace $f_c^{\circ p}$ near $z=z(c)$ by $g_{c,r}(h)=\lambda
 h\big(1+ah^r+\cdots\big)$ near $h=0$. Hence we can replace $f_c^{\circ rp}$
near $z(c)$ by
$$	g_{c,r}^{\circ r}(h)~=~\lambda^r\,h\big(1\+a'\,h^{r}+\hts\big)~, $$
where computation shows that $a'=\big(1+\lambda^r+\lambda^{2r}+\cdots+
\lambda^{(r-1)r}\big)a$. Here the coefficient $a'$ of $h^{r}$
must be non-zero when $\lambda=\hat\lambda$, and hence for $\lambda$ close
to $\hat\lambda$. For otherwise, the Leau-Fatou flowers around the
points of the parabolic orbit would give rise to
more than one periodic cycle of attracting
petals for $f_c$. This is impossible, since each such cycle must contain a
critical point, and a quadratic polynomial has only one critical point.
Finally, after a scale change, replacing $g_{c,r+1}(h)$ by $F_c(w)=
\alpha_c\, g_{c,r+1}^{\circ r}
(w/\alpha_c)$ for suitably chosen $\alpha_c$, we obtain simply
$$	F_c^{\circ r}(w)~=~ \lambda^rw\big(1\+w^{r}\+({\rm higher~terms~in~}w)
\big)~, $$
as required.\QED\ssk

{\tf Proof of 4.1.} First note that we can choose a smooth path in parameter
space so that the multiplier $\lambda^r$ of Lemma 4.2
is real and belongs
to some interval $(1\,,\,1+\eta)$. This follows easily from the fact that
$\lambda$ is a non-constant holomorphic function of $c$ in the case
$r>1$, or of $u=\sqrt{c-\hat c}$ in the case $r=1$.
Note that the map $F^{\circ r}$ of 4.2 satisfies
$$ |F^{\circ r}(w)|~=~\lambda^r\cdot|w|\cdot\big(1+{\rm Re}(w^r)+\hts\big)
 \eqno (4)$$
and
$$	\arg(F^{\circ r}(w))~=~\arg(w)\+{\rm Im}(w^r)+\hts \eqno (5)$$
whenever $\lambda^r$ is real and positive; and note
also that $F^{\circ r}$ has a locally defined holomorphic inverse of the form
$$	F^{-r}(w)~~=~~w\,\big(1\,-\,w^r/\lambda^{2r}
\+\hts\big)/\lambda^r~, $$
which satisfies
$$ |F^{-r}(w)|~=~|w|\,\big(1-{\rm Re}(w^r)/\lambda^{2r}+
\hts\big)/\lambda^r \eqno (4')$$
and
$$	\arg(F^{-r}(w))~=~\arg(w)\,-\,{\rm Im}(w^r)/\lambda^{2r}+\hts~.
 \eqno (5')$$

\begin{figure}[htb]
\cl{\psfig{figure=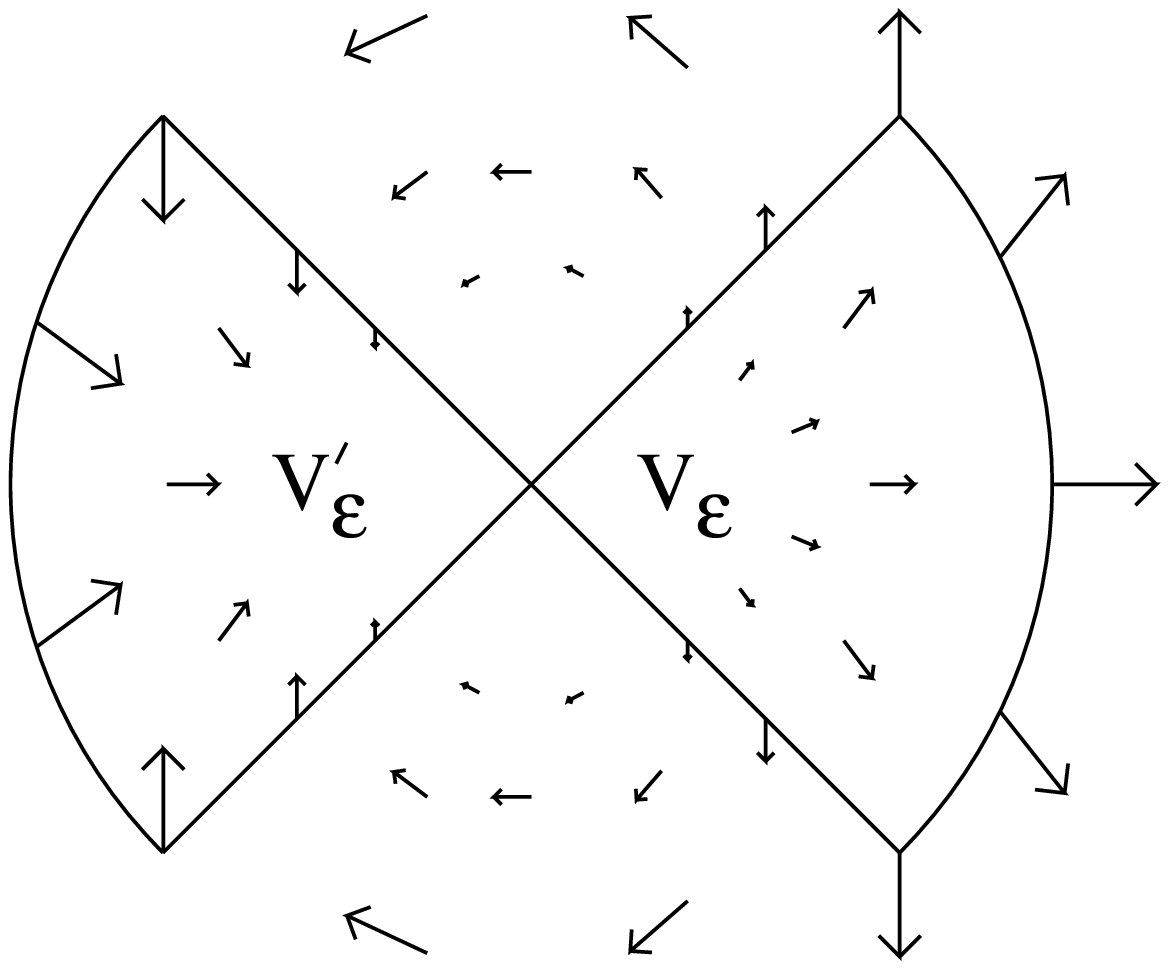,height=2in}}
\begin{quote}
{\bit Figure 11. A repelling petal $V_\epsilon$ and attracting petal
$V'_\epsilon$ for the map\break $F(w)\approx w+w^2$
(illustrating the primitive case, before perturbation).\par}
\vskip -.1in
\end{quote}\end{figure}

As a representative repelling petal for $F^{\circ r}$ let us choose
a small wedge shaped region $V_\epsilon$ described in polar
coordinates by setting $w=\rho\,e^{2\pi it}$ with $0\le \rho\le\epsilon$
and $|t|<1/(8r)$. (Compare Figure 11 for the case $r=1$.)
If $\lambda^r\ge 1$ with $\lambda^r$ sufficiently
close to 1, it follows easily from $(4')$ and $(5')$
that $V_\epsilon$ maps into itself under
$F^{-r}$, with all orbits converging towards the boundary
fixed point at $w=0$.
If a dynamic ray for $f_{\hat c}$ lands at $\hat z$, then
it must land through one of the $r$ repelling petals, for example
through the image of $V_\epsilon$ in the $z$-plane.
For $c$ sufficiently close to $\hat c$, this image must still contain
a full segment,
from some point $z$ to $f_c^{\circ rp}(z)$, of the perturbed ray, hence
this perturbed ray must still
land at the repelling point which corresponds to $w=0$.

Note that no new rays land at this point, after perturbation. There are only
finitely many rays which have period $p$. But every dynamic ray of
period $p$ for $f_{\hat c}$ with angle not in the set ${\bf A}_\po$
of angles for $\po$ must land on some disjoint repelling point, and this
condition will be preserved under perturbation. Thus the perturbed orbit,
for $\lambda^r>1$, still has portrait $\po$.

As an attracting petal for $F^{\circ r}$ we can choose the set
$V'_\epsilon=e^{\pi i/r}V_\epsilon$ consisting of all
$w=\rho e^{2\pi it}$ with $0\le\rho\le\epsilon$ and ${3\over 8r}\le t\le
{5\over 8r}$. If $\lambda^r> 1$ with $\lambda^r$ close to 1,
then using (4) and (5) we can check that $F^{\circ r}$
maps $V'_\epsilon$ into
itself. However, the origin is a repelling point, so orbits cannot converge to
it. In fact, if $K$ is the compact set obtained from $V'_\epsilon$
by removing a very small neighborhood of the origin, then $F^{\circ r}$
maps $K$ into its own interior. It follows easily that all orbits in
$V'_\epsilon\ssm\{0\}$ converge to an interior fixed point. This must be
a strictly attracting point, and must correspond to an attracting orbit
of period $rp$ for the map $f_c$.\QED

{\ssk\nin{\tf Corollary 4.3. Parabolic Points as Root Points.} \it If
$f_{\hat c}$ has a parabolic orbit whose portrait $\po$ is non-trivial,
then $\hat c$ must be equal to the root
point ${\bf r}_\po$ of the $\po$-wake.\ssk}

{\tf Note:} The hypothesis that $\po$ is non-trivial is actually redundant.
(See 4.8.) It will be
shown in 5.4 that every parabolic point is the root point of only one wake,
so that the root point of the $\po$-wake always has portrait equal to
$\po$.\ssk

{\tf Proof of 4.3.}
Since $f_{\hat c}$ has a parabolic orbit with portrait $\po$,
it certainly cannot have a {\it repelling\/} orbit with portrait $\po$.
Hence it cannot be inside the $\po$-wake by 3.1. On the other hand, by
4.1 it must belong to the boundary of the
$\po$-wake. By construction, the root point ${\bf r}_\po$ is the
only boundary point of $W_\po$ which belongs to the Mandelbrot set.\QED


Here is a complementary statement to 4.1, in the case $r>1$.

{\ssk\nin{\tf Lemma 4.4. A Deformation Breaking the Portrait.}
\it Under the hypothesis of 4.1, there also exists
a smooth path of parameter values $c$, converging to $\hat c$, so that
each $f_c$ has an attracting orbit of period $p$, and a repelling orbit
of period $rp$ which lies on the boundary of its immediate basin.
Furthermore, the dynamic rays with angles in ${\bf A}_\po=A_1\cup\cdots\cup
A_p$ all land on this repelling orbit.\ssk}

(Compare Figure 10, right.)
For such values of $c$ (still assuming that $r>1$), it follows
that there is {\it no\/} periodic orbit with portrait $\po$. Together
with 4.1, this gives an alternative proof that $\hat c$ is on the boundary
of the $\po$-wake.

The proof of 4.4 is completely analogous to the proof of 4.1, and will be
left to the reader: One simply
deforms so that $\lambda^r<1$, instead of $\lambda^r>1$.\QED\ssk

The following assertion helps to make the statement
of 4.2 more precise.

{\ssk\nin{\tf Lemma 4.5. Local Uniqueness.} \it Under the hypothesis of 4.2,
there exist unique single valued
functions $~c=c(\lambda)~$ and $~z=z(\lambda)\,, $
defined and holomorphic for
$\lambda$ in a neighborhood of $\hat\lambda$, so that $z(\lambda)$ is a
periodic point of period $p$ and multiplier $\lambda$ for the map
$f_{c(\lambda)}$, with $\hat c=c(\hat\lambda)$ and $\hat z=z(\hat\lambda)$.
This function $c(\lambda)$ is univalent in the
satellite case, but has a simple
critical point at $\hat\lambda$ in the primitive
case.\ssk}

The implications of this lemma for the geometry of the Mandelbrot set
will be described in 6.1 and 6.2.\ssk

{\tf Proof of 4.5.} First consider the satellite case, with $\hat
\lambda\ne 1$. Then clearly
the period $p$ orbit and its multiplier $\lambda(c)$ depend
smoothly on $c$ throughout some neighborhood of $\hat c$.
We will show that the derivative $d\lambda/dc$
is non-zero at $\hat c$. For otherwise, we could write
$$	\lambda^r(c) ~=~ 1 + a(c-\hat c)^k+\hts $$
with $k\ge 2$. Hence
we could vary $c$ from $\hat c$ in two or more different directions so
that $\lambda^r>1$ and in two or more intermediate directions so that
$\lambda^r<1$. The former points would be within the $\po$-wake and the
later points would be outside it; but this configuration is impossible by 3.1.
Thus $d\lambda/dc\ne 0$, and it follows by the Inverse Function Theorem
that the inverse mapping $\lambda\mapsto c(\lambda)$ is well
defined and holomorphic throughout a neighborhood of $\hat\lambda$,
as required.

In the primitive case, the situation is different, but the proof is similar.
In this case, setting $c=\hat c+u^2$, we must
express the multiplier $\lambda_1$ for one of the two nearby period
$p$ points as a holomorphic function of $u$, and show that
the derivative $d\lambda_1/du$ is non-zero at $u=0$. Otherwise,
if the derivative $d\lambda_1(u)/du$ were equal to zero
for $u=0$, then we could write
$$      \lambda_1(u) ~=~ 1 + a\,u^k+\hts $$
for some $k\ge 2$. It would follow that
we could vary $u$ from $0$ in two or more different directions so
that $\lambda_1>1$ and in two or more separating directions so that
$\lambda_2>1$. All of these points would be within the $\po$-wake,
but the rays landing on the periodic point
$z_1$ would have to jump discontinuously so
as to land on $z_2$ as we pass from $\lambda_1>1$ to $\lambda_2>1$,
and such points of discontinuity must be outside the $\po$-wake.
Even allowing for the fact that the $u$-plane is a two-fold covering of
the $c$-plane, such a configuration is incompatible with 3.1.
Therefore, $\lambda_1$ and $u$ must determine each other
holomorphically in a neighborhood of $\hat\lambda\leftrightarrow 0$.
In particular, it follows that the parameter value $c=\hat c+u^2$
can be expressed as a holomorphic function of $\lambda_1$,
with a simple critical point at $\lambda_1=\hat\lambda$.\QED
\ssk

To conclude this section, we will prove that the portrait of a parabolic
periodic point is always non-trivial. We will use a somewhat simplified form
of the Hubbard tree construction to show that every parabolic orbit
with ray period $rp\ge 2$ must have portrait with valence $v\ge 2$.
First some general remarks about locally connected subsets of the plane.\ssk

{\ssk\nin{\tf Lemma 4.6. A Canonical Retraction.}
\it Let $K\subset\C$ be compact, connected, locally
connected, and full, and let $U$
be a connected component of the interior of $K$. Then
the closure $\overline U$ is homeomorphic to the closed unit disk,
and there is a unique retraction $\rho_U$ from
$\C$ onto $\overline U$ which carries each external ray, and also
each connected component of the complement $K\ssm\overline U$, to a
single point of the circle $\partial U$. There are at least two distinct
external rays landing at a point $z_0\in\partial U$ if and only if
$K\ssm\{z_0\}$ is disconnected, or if and only if there is some connected
component $X$ of $K\ssm\overline U$ with $\rho_U(X)=\{z_0\}$.\ssk}

{\tf Proof.} (Compare [D5].) The statement that $\overline U$ is a disk
follows easily from well known results of
Carath\'eodory. Furthermore, according to Carath\'eodory,
there is a unique retraction from $\C$ onto $K$ which maps each external
ray to its landing point. Composing this with the retraction $K\to\overline
U$ which maps each component $X$ of $K\ssm\overline U$ to the unique
intersection point $z_0\in\overline X\cap\overline U$, we obtain the required
retraction $\rho_U$.

For any such $X$, note that there must be at least one
maximal open interval of
angles $t$ such that the ray $\ra_t^K$ lands in $X$. The endpoints of
such a maximal interval are the angles for the required pair of rays landing on
$z_0$. Conversely, if there were two rays landing on $z_0$ but no
component $X$ attached in between, then there would be an entire open interval
of angles $t$ so that $\ra_t^K$ lands at $z_0$. But this is impossible
by a classical theorem of F. and M. Riesz. (See for example [M2, App. A].)\QED

In particular, let $K=K(f)$ be the filled Julia set for a hyperbolic
quadratic polynomial. (We are actually interested in the parabolic case,
but will work first with the hyperbolic case,
since that will suffice for our purposes, and since it is much easier to
prove local connectivity in the hyperbolic case.)\ssk

{\ssk\nin{\tf Lemma 4.7. The Dynamic Root Point.} \it
Suppose that $f=f_c$ has an attracting orbit of period
$n\ge 2$. Let $K$ be its filled Julia set, and let $U_0$ and
$U_1\subset K$ be
the Fatou components containing the critical point $0$ and the critical
value $c$ respectively. Then the canonical
retraction $\rho_{U_1}:\C\to\overline U_1$ carries the component $U_0$
to the unique point ${\bf r}_c\in\partial U_1$ which is fixed by
$f^{\circ n}$. Hence at least two dynamic rays land at this point.\ssk}

(See for example Figures 1, 6.)
Following Schleicher, I will call ${\bf r}_c$ the {\bit dynamic root point\/}
for the Fatou component $U_1$.\ssk

{\tf Proof.} Let $U_0\to U_1\iso U_2\iso\cdots\iso U_n=U_0$ be the Fatou
components containing the critical orbit. Then $f^{\circ n}$ maps each
circle $\partial U_j$ onto itself by an expanding map of degree two.
Hence there is a canonical homeomorphism $a_j:\partial U_j\to\R/\Z$ which
conjugates $f^{\circ n}$ to the angle doubling map on the standard circle.
For each $z\in\C\ssm U_j$, the image $a_j(\rho_{U_j}(z))$
will be called the {\bit internal angle\/}
of the point $z$ with respect to $U_j$. The map $f$ from $\partial U_j$
to $\partial U_{j+1}$ preserves the internal angles of boundary
points for $0< j<n$, but
doubles them for the case $j=0$ of the critical component.

Define the {\bit $t$-wake\/} $L_t(U_j)$ to be the set of all
$z\in\C\ssm U_j$ with $a_j(\rho_{U_j}(z))=t\in\R/\Z$. These wakes are
pairwise disjoint sets with union equal to $\C\ssm U_j$. In general
$f$ maps to $t$-wake of $U_j$ homeomorphically onto the $t$-wake
of $U_{j+1}$ for $0<j<n$, and onto the $2t$-wake of $U_{j+1}$ when
$j=0$. However, there is one exceptional value of $t$ for each $U_j$
with $0<j<n$. Namely, if the wake $L_t(U_j)$ contains the critical component
$U_0$ then it certainly cannot map homeomorphically, and its image may
be much larger than $L_{t}(U_{j+1})$.

Let ${\mathcal A}_j\subset\R/\Z$ be the finite set consisting of all
angles $t\in\R/\Z$ such that the wake $L_t(U_j)$ contains one of the
components $U_k$ (where necessarily $j\ne k$). Then
it follows that ${\mathcal A}_1\subset{\mathcal A}_2\subset\cdots
\subset{\mathcal A}_n$ and $2{\mathcal A}_n={\mathcal A}_0$. On the other hand, since
$K$ is full, the various $U_j$ must be connected together in a tree-like
arrangement (the {\bit Hubbard tree\/}). There cannot be any cycles.
Hence at least one of the ${\mathcal A}_i$ must consist of a single angle.
It follows easily that ${\mathcal A}_1=\{0\}$, and the conclusion follows.\QED

{\ssk\nin{\tf Corollary 4.8. Parabolic Orbit Portraits are Non-Trivial.} \it
If $c$ is any parabolic point of the Mandelbrot set
other than $c=1/4$, and if $\Or$ is the parabolic orbit for
$f_c$, then at least two dynamic rays land on each point of $\Or$.\ssk}

(This is just a restatement of Theorem 1.4 of \S1.)\ssk

{\tf Proof.} In the satellite case this is trivially true, while in the
primitive case it follows from 4.7, using 4.1 to pass from the
parabolic to the hyperbolic case.
This completes the proof of Theorem 1.4.\QED

\section{The Period $n$ Curve in (Parameter$\times$Dynamic) Space.}
It is convenient to define a sequence of numbers $\nu_2(n)$
inductively by the formula
$$	2^k~=~\sum_{n|k} \nu_2(n)~,\qquad{\rm or}\quad \nu_2(k)~=~\sum_{n|k}\mu
(k/n)2^n~, $$
to be summed over all divisors $n\ge 1$ of $k$, where $\mu(k/n)\in\{\pm 1, 0\}$
is the M\"obius function. In fact we will be mainly
interested in the quotients $\nu_2(n)/2$ and $\nu_2(n)/n$.
The first few values are\ssk

\cl{\vbox{\halign{~#~&~#~&~#~&~#~&~#~&~#~&~#~&~#~&~#~&~#~&~#~&~#~&~#~\cr
	$n$&~&1&2&3&4&5&6&7&8&9&10\cr
	$\nu_2(n)/2$&~&1&1&3&6&15&27&63&120&252&495\cr
	$\nu_2(n)/n$&~&2&1&2&3&6&9&18&30&56&99&.\cr}}}\ssk

Define the {\bit period $n$ curve\/} $\Per_n\subset
\C^2$ to be the locus of zeros of the polynomial
$Q_n(c,z)$ which is defined by the formula
$$	f_c^{\circ k}(z)-z~=~\prod_{n|k} Q_n(c,z)~,\quad{\rm or}\quad
Q_k(c,z)~=~\prod_{n|k} \big(f_c^{\circ k}(z)-z\big)^{\mu(k/n)}~, $$
taking the product over all divisors $n$ of $k$. For example,
$$	Q_1(c,z)~=~z^2+c-z\,,\qquad
 Q_2(c,z)~=~{(z^2+c)^2+c-z\over z^2+c-z}~=~ z^2+z+c+1~. $$
Note that each point $(c,z)\in\Per_n$ determines a periodic orbit
$$	z~=~z_0\mapsto z_1\mapsto\cdots\mapsto z_n~=~z_0 $$
for the map $f_c$. Let $\lambda_n=\lambda_n(c,z)=
\partial f_c^{\circ n}(z)/\partial z=2^nz_1\cdots z_n$.
For a generic choice of $c$, this orbit has period exactly $n$, and $\lambda_n$
is the multiplier.
However, if $z$ is a parabolic periodic point for $f_c$ with ray period
$n=rp>p$, then $(c,z)$ belongs both to $\Per_n$ with $\lambda_n=1$, and to
$\Per_p$ with $\lambda_p\in\root r\of 1$. (In fact,
the two curves $\Per_n$ and $\Per_p$ intersect transversally at $(c,z)$.)

{\tf Remarks.} Compare [M4] for a somewhat analogous discussion
for cubic polynomials. The fact that $Q_n$ is really a polynomial can be
verified by expressing $f_c^{\circ k}(z)-z$ as a product of irreducible
polynomials, and checking that each of these irreducible
factors has a well defined period $n$ dividing $k$. The factors are all
distinct since $\partial(f_c^{\circ j}(z)-z)/\partial z\ne 0$ at every zero
of this polynomial when $|c|$
is large. It is shown in [Bou], and also in [S1], [LS], that the
algebraic curve $\Per_n$ (or the polynomial $Q_n$)
is actually irreducible; however, we will not make any use of that fact.

{\ssk\nin{\tf Lemma 5.1. Properties of the Period $n$ Curve.} \it
This algebraic curve $\Per_n\subset\C^2$ is
non-singular. The projection $(c,z)\mapsto c$ is a proper map of degree
$\nu_2(n)$ from $\Per_n$ to the parameter plane, while the projection
$(c,z)\mapsto z$ is a proper map of degree $\nu_2(n)/2$ to the dynamic
plane. Finally, the function $(c,z)\mapsto\lambda_n(c,z)$ is a
proper map of degree $\,n \nu_2(n)/2$ to the $\lambda_n$-plane.\ssk}

Note that the cyclic group of order $n$, which we will denote by $\Z_n$,
acts on $\Per_n$, a generator carrying $(c,z)$ to $(c\,,\,f_c(z))$.

{\ssk\nin{\tf Lemma 5.2. Properties of $\Per_n/\Z_n$.} \it The quotient
$\Per_n/\Z_n$ is a smooth algebraic curve consisting of all
pairs $(c,\Or)$ where $\Or$ is a periodic orbit for $f_c$ which is either
non-parabolic of period $n$, or parabolic with attracting petals of
period $n$. At any point where $\lambda_n\ne 1$, the coordinate
$c$ can be used as local uniformizing parameter, while in a neighborhood
of a point with $\lambda_n=1$, the multiplier $\lambda_n=\lambda_n(c,z)$ serves
as a local uniformizing parameter for this curve.
The projection maps $(c,\Or)\mapsto c$ and $(c,\Or)\mapsto \lambda_n$
are proper, with degrees $\nu_2(n)/n$ and $\nu_2(n)/2$ respectively.\ssk}

The proof that $\Per_n$ and $\Per_n/\Z_n$ are
non-singular will be divided into three cases, as follows.\ssk

{\bf Generic Case.} First consider a point $(\hat c,\hat z)\in \Per_n$ with
$\lambda_n(\hat c,\hat z)\ne 1$. Then, by the
Implicit Function Theorem, we can solve the equation $f_c^{\circ n}(z)=z$
locally for $z$ as a smooth function of $c$. It follows that both of the curves
$\Per_n$ and $\Per_n/\Z_n$ are locally smooth, with $c$ as local
uniformizing parameter.\ssk

{\bf Primitive Parabolic Case.} Now consider a point $(\hat c,\hat z)\in
\Per_n$ with $\lambda_n(\hat c,\hat z)=1$, where $\hat z$ has
period exactly $n$ under $f_{\hat c}$. According to the proof of 4.5,
if we set $c=\hat c+u^2$, then both $z$ and $\lambda_n=\lambda_n(c,z)$
can be expressed locally as smooth functions of $u$ with
$d\lambda_n/du\ne 0$. It follows that both $\Per_n$ and $\Per_n/\Z_n$
are locally smooth at this point, and that we can use
either $u$ or $\lambda_n$ as local
uniformizing parameter. (Similarly $dz/du\ne 0$, so we could use $z$ as
local uniformizing parameter for $\Per_n$. However $dc/du$ is zero when $u=0$,
so $c$ cannot be used as local parameter.)\ssk

{\bf Satellite Parabolic Case.} Again suppose that $\lambda_n(\hat c,\hat z)=1$,
but now assume that the period $p$ of $\hat z$ is strictly less than
the ray period $n=rp$. For $c$ near
$\hat c$, let $z=z(c)$ be the equation of the unique period $p$
point near $\hat z$. Using the change
of variable $w=\alpha (z-z(c))+\hts$ of 4.2, the map $f_c^{\circ n}$
corresponds to
$$  w\mapsto F^{\circ r}(w)~=~\lambda^rw\big(1\,+\,w^r\,+\hts\big)~,
	\eqno (6)$$
where $\lambda=\lambda(w)$ is the multiplier of this period $p$ orbit.
The equation for a fixed point is\break $w=\lambda^rw\,(1+w^r+\hts)$.
Dividing by $w$
(since we want the fixed point with $w\ne 0$ or with $z\ne z(c)$), this
becomes
$$	1=\lambda^r(1+w^r+\hts)\qquad{\rm or}\qquad\lambda^r~=~1-w^r+\hts~.$$
Thus we can express $\lambda$ as a
holomorphic function of $w$, with a critical point at $w=0$. Therefore,
by 4.5, we can also express $c$ as a holomorphic function of $w$. Since
$w$ is defined as a holomorphic function of $z$ and $c$ with
$\partial w/\partial z\ne 0$, it follows that
$\Per_n$ is locally smooth with local uniformizing parameter $z$ or $w$.

Now note that there is a unique local change of coordinate $w\mapsto\phi(w)$
with $\phi'(0)=1$ so that $\lambda^r=1-\phi(w)^r$. Since the
expression $\phi(w)^r$ is invariant under the $\Z_n$ action of $\Per_n$,
it follows easily that this action can be described by the formula
$\phi(w)\mapsto\hat\lambda\phi(w)$. It follows that $\phi(w)^r
=1-\lambda^r$ is a local uniformizing parameter for the quotient curve
$\Per_n/\Z_n$. Therefore, either $\lambda$ or $c$ can also be taken
as local uniformizing parameter. In particular, it follows that
the multiplier $\lambda_n$
of the period $n=rp$ orbit can be expressed as a smooth function of
the multiplier $\lambda=\lambda_p$ of the period $p$ orbit. Note that
$$	d\lambda_n/d(\lambda^r)~=~ -r \eqno (7)$$
at the parabolic point. (Compare [CM (4.3)].)
This can be verified by direct computation from (6), or
by using the holomorphic fixed point formula [M2]
for the function $f_c^{\circ n}$ to show that the expression
$$	{r\over 1-\lambda_n}\+{1\over 1-\lambda^r} $$
depends smoothly on the parameter $c$ throughout some neighborhood of the
parabolic point. Therefore $\lambda_n$ can also be used as local uniformizing
parameter for $\Per_n/\Z_n$.\ssk

The degrees of the various
projection maps can easily be computed algebraically,
by counting solutions to the appropriate polynomial equations. Here is
a more geometric argument, which also
provides a quite explicit description of the
ends of the curve $\Per_n$, and hence proves that these mappings are
proper. Let us consider the limiting
case as $|c|\to\infty$. Setting $c=-v^2$ with $|v|>2$, let $\pm\Delta$
be the open
disk of radius $1$ centered at $\pm v$. It is not difficult to check
that both $\Delta$ and $-\Delta$ map holomorphically onto a disk
$f(\Delta)$ which contains $\overline\Delta\cup(-\overline\Delta)$.
The (filled) Julia set $K$
can then be described explicitly as follows. Given an arbitrary sequence
of signs $\epsilon_0\,,\,\epsilon_1\,,\,\ldots$, there is one an only one
orbit $z_0\mapsto z_1\mapsto\cdots$ in $K$ with $z_j\in\epsilon_j\Delta$
for every $j\ge 0$. This is proved using the Poincar\'e metric for the
inverse maps $f(\Delta)\to\pm\Delta\subset f(\Delta)$. In particular,
the number of solutions of
period $n$ is equal to the number of sign sequences of period $n$,
which is easily seen to be $\nu_2(n)$. Thus the degree
of the projection to the $c$-plane is $\nu_2(n)$.
It follows also that the product $z_1\cdots z_n=\lambda/2^n$ is
given asymptotically by
$$ \lambda/2^n~\sim~\pm  z^n~\sim~\pm  v^n~=~\pm (-c)^{n/2}
	\qquad{\rm as}\qquad |v|\to\infty~. $$
Thus the degree of the projection to the $\lambda$-plane is
$n$ times the degree
of the projection to the $z$-plane, and is $n/2$ times the degree
of the projection to the $c$-plane.\QED\ssk

Thus we have a diagram of smooth algebraic curves and proper holomorphic
maps with degrees as indicated:
$$\begin{matrix}
 \Per_n&\buildrel n\over\longrightarrow&\Per_n/\Z_n&  \buildrel
    \nu_2(n)/2\over{-\!\!-\!\!\!\longrightarrow}&\lambda_n\hbox{\rm -plane}\cr
 \qquad\downarrow{\scriptstyle \nu_2(n)/2} & &
     \qquad\downarrow{\scriptstyle \nu_2(n)/n}\cr
 z\hbox{\rm -plane}& & c{\rm \hbox{-}plane}\cr\end{matrix}$$
For a generic choice of $c$, it follows that the map $f_c$ has exactly
$\nu_2(n)/n$ periodic orbits of period $n$, while for generic
choice of $\lambda_n$
there are exactly $\nu_2(n)/2$ pairs $(c,\Or)$
consisting of a parameter value $c$
and a period $n$ {\it orbit\/} of multiplier $\lambda_n$ for the map $f_c$.
The discussion shows that the correspondence $(c,\Or)\mapsto (c,\lambda_n)$
yields a smooth immersion of $\Per_n/\Z_n$ into $\C^2$.
(Caution: Presumably some $f_c$ may have two different period $n$
orbits with the same multiplier, so this immersion
may have self-intersections.)\ssk

{\ssk\nin{\tf Corollary 5.3. Counting Parabolic Points.} \it The number of
parabolic points in the Mandelbrot set with ray period $rp=n$ is equal to
$\nu_2(n)/2$.\ssk}

{\tf Proof.} This is the same as the number of points in the pre-image of
$+1$ under the projection $(c,\Or)\mapsto \lambda_n(c,\Or)$ from
$\Per_n/\Z_n$ to the $\lambda_n$-plane. According to 5.2, the degree of this
projection is $\nu_2(n)/2$, and $+1$ is a regular value. The conclusion
follows.\QED\ssk

We are now ready to prove the main results, as stated in \S1.

{\ssk\nin{\tf Corollary 5.4.} \it There are exactly two parameter rays
which angles which are periodic under doubling landing
at each para\-bolic point $\hat c\ne 1/4$. Hence
distinct wakes have distinct root points; and for
each non-trivial portrait $\po$, the root point of the $\po$-wake has a
parabolic orbit with portrait $\po$.\ssk}

(For angles which are not periodic, compare 9.4.)

{\ssk\nin{\tf Corollary 5.5.} \it Every parameter ray $\ra_t^M$ whose
angle has period $n\ge 2$ under doubling forms one of the two
boundary rays for one and only one wake $W_\po$, where $\po$ is some portrait
with ray period $n$.\ssk}

{\tf Proof of 5.4 and 5.5.} According to 5.3, the number of parabolic points
$\hat c$ with ray period $n\ge 2$ is equal to $\nu_2(n)/2$, and according to
Theorem 1.4 each such point is the landing point of at least two rays,
which necessarily have ray period $n$. Thus altogether
there are at least $\nu_2(n)$ distinct rays of period $n$.
On the other hand, since the map $t\mapsto 2^n t~~({\rm mod~}\Z)$ has $2^n-1$
fixed points, it follows inductively that the number of angles with period
exactly $n\ge 2$ is precisely  equal to $\nu_2(n)$. Thus
there cannot be more than two rays landing at any such point $\hat c$. It
follows that $\hat c$ is the root point of at most one wake. For if
$\hat c$ were the root point of two different wakes, then (even if they
shared a boundary ray) it would be the landing point for at least three
different parameter rays. Using 4.3, it now follows that each such
$\hat c$ is the root point ${\bf r}_\po$ for exactly one wake $W_\po$,
and furthermore that each $f_{{\bf r}_\po}$ has a parabolic orbit
with portrait $\po$.

Here we have assumed that $n\ge 2$. However, for $n=1$ there is clearly
just one parameter ray $\ra_0^M=(1/4,\infty)$ which is fixed under doubling,
and its landing point $\hat c=1/4$ is the unique parabolic point with
ray period $n=1$. This completes the proof of
5.4 and 5.5. Clearly Theorems 1.2 and 1.5,
as stated in \S1, follow immediately.\QED\ms

To conclude this section, here is a more explicit
description of the first few period $n$ curves:\ssk

{\bf Period 1.} The curve $\Per_1=\Per_1/\Z_1\cong\C$ can be identified
with the $\lambda_1$-plane. It is a $2$-fold branched cover of the
$c$-plane, ramified at the root point $r_{\{\{0\}\}}=1/4$,
and can be described by the equations $z=\lambda_1/2\,,\;c=z-z^2$.
Note that the unit disk $|\lambda_1|<1$
in the $\lambda_1$-plane maps homeomorphically onto the region bounded
by the cardioid in the $c$-plane.\ssk

{\bf Period 2.} The quotient
$\Per_2/\Z_2\cong\C$ can be identified either with the $\lambda_2$-plane
or with the
$c$-plane, where $\lambda_2=4\,(1+c)$. The curve $\Per_2\cong\C$
is a $2$-fold branched cover with coordinate $z$,
branched at the point $\lambda_2=1$
which corresponds to the period $2$ root point
$c=r_\po=-3/4$ with portrait $\po={\{\{1/3,2/3\}\}}$.
It is described by the equation $z^2+z+(c+1)=0$, with
$\Z_2$-action 
$z\leftrightarrow f_c(z)=-z-1$.\ssk

{\bf Period 3.} (See [Giarrusso and Fisher].)	
The quotient $\Per_3/\Z_3\cong\C$ can be identified with a
$2$-fold branched cover of the $c$-plane, branched at the root point
$r_\po=-7/4$ of the real period $3$ component, where
$\po=\{\{3/7,4/7\}, \{6/7,1/7\},\{5/7,2/7\}\}$. If we choose a parameter
$u$ on this quotient by setting $c=-(u^2+7)/4$, then computation shows that
the multiplier is given by the cubic expression $\lambda_3=u^3-u^2+7\,u+1$.
The curve $\Per_3$ itself is conformally
isomorphic to a thrice punctured Riemann sphere. It can be described
as a $3$-fold cyclic branched cover of this $u$-plane, branched
with ramification index $3$ at the two points $u=(1\pm\sqrt{-27})/2$
where $\lambda_3=1$.

\section{Hyperbolic Components.} By definition, a {\bit hyperbolic component\/}
$H$ of {\bit period\/} $n$ in the Mandelbrot set is a connected
component of the open set consisting of all parameter values
$c$ such that $f_c$ has a (necessarily unique)
attracting orbit of period $n$. We will first study the geometry of
a hyperbolic component near a parabolic boundary point. 

{\ssk\nin{\tf Lemma 6.1. Geometry near a Satellite Boundary Point.} \it Let $\hat
 c$ be a parabolic point with orbit portrait $\po$ having ray period
$rp>p$. Then $\hat c$ lies on the boundary of exactly two hyperbolic
components. One of these has period $rp$ and lies inside the $\po$-wake,
while the other has period $p$ and lies outside the $\po$-wake. Locally
the boundaries of these components are smooth curves which
meet tangentially at $\hat c$.\ssk}

\begin{figure}[htb]
\cl{\psfig{figure=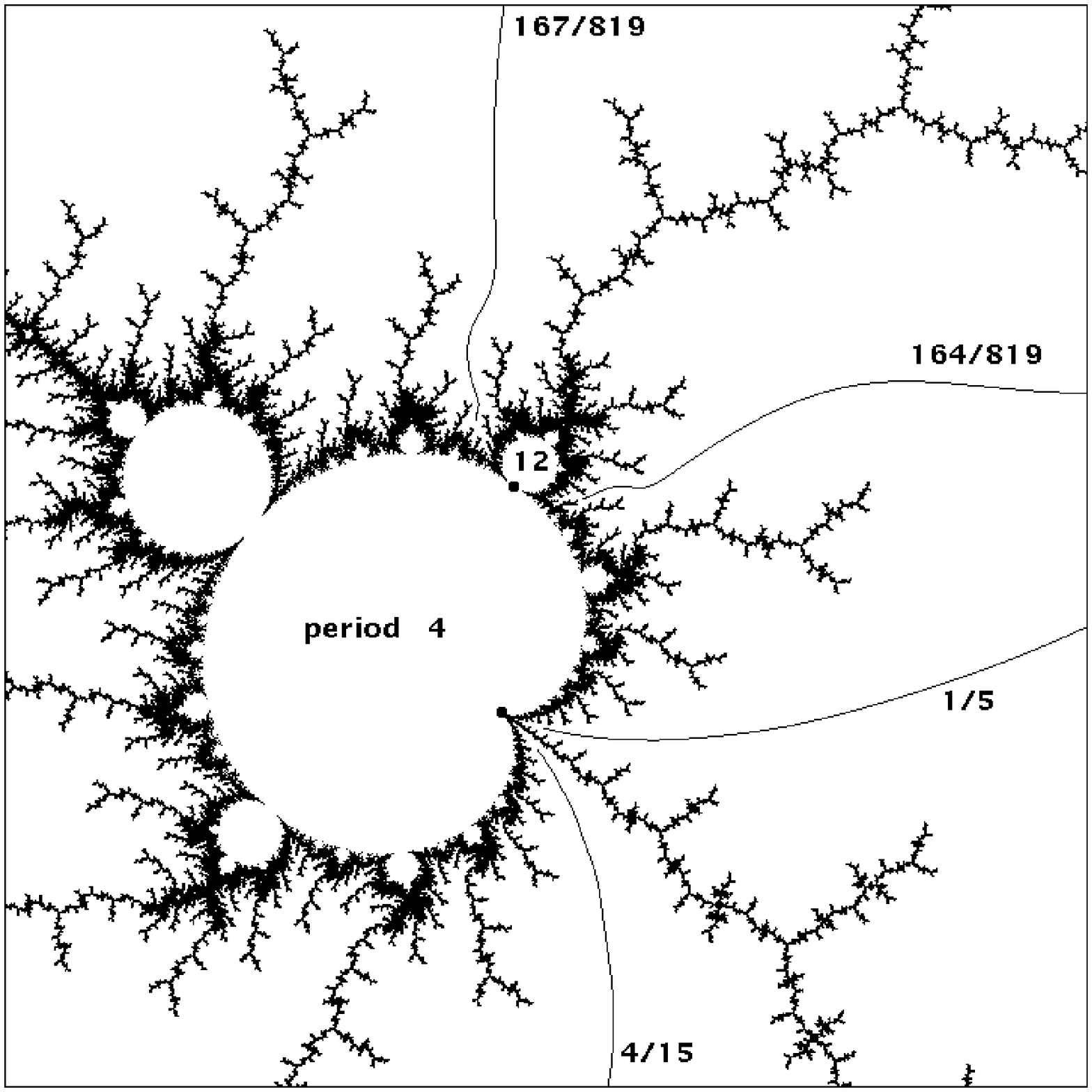,height=3in}}\ssk
\begin{quote}{\bit Figure 12. Detail of the Mandelbrot boundary, showing the
rays landing at the root points of a primitive period 4 component and a
satellite period 12 component.}\vskip -.5in
\end{quote}\end{figure}

{\tf Proof.} According to 4.1, $\hat c$ lies on the boundary of a
hyperbolic component $H_{rp}$ of period $rp$ which lies inside the
$\po$-wake,
while according to 4.4 it lies on the boundary of a component $H_p$
of period $p$ which lies outside the $\po$-wake. Let $\Or_{rp}$ and $\Or_p$
be the associated periodic orbits, with multipliers $\lambda_{rp}$ and
$\lambda_p$. According to 4.5, the multiplier $\lambda_p$ can be used as
a local uniformizing parameter for the $c$-plane near $\hat c$. Therefore
the boundary $\partial H_p$, with
equation $|\lambda_p|=1$, is locally smooth. Similarly, it follows from
equation (7) of \S5, that we can take $\lambda_{rp}$ as local uniformizing
parameter, so the locus $|\lambda_{rp}|=1$ is also locally smooth. These
two boundary
curves are necessarily tangent to each other since the two hyperbolic
components cannot overlap, or by direct computation from (7).

To see that there are no other components with $\hat c$ as boundary point,
first note that all periodic orbits for the map $f_{\hat c}$, other than its
designated parabolic orbit, must be strictly repelling. For any orbit with
multiplier $|\lambda|\le 1$ must either attract the critical orbit (in the
attracting or parabolic case) or at least be in the $\omega$-limit set of the
critical orbit (in the Cremer case), or have Fatou component boundary in
this $\omega$-limit set (in the Siegel disk case). Since the unique
critical orbit converges to the parabolic orbit, all other periodic orbits
must be repelling.

Now choose some large integer $N$. If we choose $c$ sufficiently
close to $\hat c$, then all 
repelling periodic orbits of period $\le N$
for $f_{\hat c}$ will deform to repelling periodic orbits of the same
period for $f_c$. Thus any non-repelling orbit of period $\le N$ for $f_c$
must be one of the two orbits $\Or_p$ and $\Or_{rp}$ which arise from
perturbation of the parabolic orbit.
In other words, any hyperbolic component $H'$ of period $\le N$
which intersects
some small neighborhood of $\hat c$ must be either $H_p$ or $H_{rp}$.
In particular, any hyperbolic component which has $\hat c$ as boundary point
must coincide with either $H_p$ or $H_{rp}$.\QED

By definition, the component $H_{rp}$ is a {\bit satellite\/} of $H_p$,
attached at the parabolic point $\hat c$.
(It follows from (7) that $|d\lambda_{rp}/d\lambda_p|=r^2$ at $\hat c$,
so to a first approximation the component $H_{p}$ is $r^2$ times as big as
its satellite $H_{rp}$. Compare [CM].)

{\ssk\nin{\tf Lemma 6.2. Geometry near a Primitive Boundary Point.} \it If
the portrait $\po$ of the parabolic point $\hat c$ has ray period
$rp=p$, then $\hat c$ lies on the boundary of just one hyperbolic
component $H$, which has period $p$ and lies inside the $\po$-wake.
The boundary of $H$ near $\hat c$ is a smooth curve,
except for a cusp at the point $\hat c$ itself.\ssk}

{\tf Proof.} As in the proof of 4.2, we set $c=\hat c+u^2$ and find
a period $p$ point $z(u)$ with multiplier $\lambda(u)$
which depends smoothly on $u$, with $d\lambda/du\ne 0$. Hence the
locus $|\lambda(u)|=1$ is a smooth curve in the $u$-plane, while its image
in the $c$-plane has a cusp at $c=\hat c$. The rest of the argument is
completely analogous to the proof of 6.1.\QED

{\ssk\nin{\tf Lemma 6.3. The Root Point of a Hyperbolic Component.} \it Every
parabolic point of ray period $n=rp$ is on the
boundary of one and only one hyperbolic component of period $n$.
Conversely, every hyperbolic component of period $n$ has one and only
one parabolic point of ray period $n$ on its boundary. In this way, we obtain
a canonical one-to-one correspondence between parabolic points and
hyperbolic components in parameter space.\ssk}

{\tf Proof.} The first statement follows immediately from 6.1 and 6.2.
Conversely, if $H$ is a
hyperbolic component of period $n$, then we can map $H$ holomorphically
into the open unit disk ${\bf D}$ by sending each $c\in H$ to the multiplier of
the unique attracting orbit for $f_c$.
In order to extend to the closure $\overline H$, it is convenient to
lift to the curve $\Per_n/\Z_n$, using the proper holomorphic map
$(c\,,\,\Or)\mapsto c$
of \S5. Evidently $H$ lifts biholomorphically to an open
set $H^\natural\subset\Per_n/\Z_n$, which then maps holomorphically to the
$\lambda_n$-plane under the projection $(c,\Or)\mapsto\lambda_n(c,\Or)$.
(Here $H^\natural$ is a connected component of the set of
$(c,\Or)$ such that $\Or$ is an attracting period $n$ orbit for $f_c$.)
Since the projection to the $\lambda_n$-plane is open and proper, it follows
easily that the closure $\overline H^\natural$ maps {\it onto\/} the
closed disk $\overline{\bf D}$. In particular, there exists a point
$(\hat c\,,\,\hat\Or)$ of $\overline H^\natural$ with
$\lambda_n(\hat c,\hat O)=+1$. Evidently this $\hat c$ is a parabolic
boundary point of $H$ with ray period dividing $n$, and it follows from
6.1 and 6.2 that it must have ray period precisely $n$.

According to 4.7, for each $c\in H$ there is a unique repelling
orbit of lowest period on the boundary of the immediate basin for
the attracting orbit of $f_c$. Furthermore, according to 4.1, the portrait
$\po=\po_H$ for this orbit is the same as the portrait for the parabolic orbit
of $f_{\hat c}$. Since there is only one parabolic point with specified
portrait by Theorem 1.2, this proves that there can only one such point
$\hat c \in\partial H$.\QED\ssk

{\tf Definition.} This distinguished parabolic point on the boundary
$\partial H$ of a hyperbolic component is called the {\bit root
point\/} of the hyperbolic component
$H$. We know from 1.2 and 1.4 that the parabolic points
of ray period $n$ can be indexed by the non-trivial orbit portraits of ray
period
$n$. {\it Hence the hyperbolic components of period $n$ can also be indexed
by non-trivial portraits of ray period $n$.\/}
We will write $H=H_\po$ (or $\po=\po_H$)
if $H$ is the hyperbolic component with root point $r_\po$. We will say that
$H$ is a {\bit primitive component\/} or a {\bit satellite component\/}
according as the associated portrait is primitive or satellite.\ssk

\nin
{\tf Remark 6.4.} Of course there are many other parabolic points in $\partial
H$. For each root of unity $\mu=e^{2\pi iq/s}\ne 1$ a similar
argument shows that there is at least one point $(\hat c_\mu,\Or_\mu)\in
\partial H^\natural$ with $\lambda_n(\hat c_\mu\,,\,\Or_\mu)=\mu$. In fact
the following theorem implies
that $\hat c_\mu$ is unique. This $\hat c_\mu$ is the root
point for a hyperbolic component $H'$ of period $sn>n$, with associated orbit
portrait $\po'$ of period $n$ and rotation number $q/s$. By definition,
$\po'$ is the $(q/s)$-{\bit satellite\/} of $\po$, and
$H'$ is the $(q/s)$-{\bit satellite\/} of $H$.\ssk

We next prove the following basic result of Douady and Hubbard.
Again let $H$ be a hyperbolic component of period $n$ and let
$H^\natural\subset\Per_n/\Z_n$ be the set of pairs
$(c,\Or)$ with $c\in H$, where $\Or$ is the attracting orbit for $f_c$.

{\ssk\nin{\tf Theorem 6.5. Uniformization of Hyperbolic Components.} \it
The closure $\overline H$ is
homeomorphic to the closed unit disk $\overline{\bf D}$. In fact there is
a canonical homeomorphism
$$	\overline\D~\cong~\overline H^\natural~\to~\overline H $$
which carries each point $\lambda\ne 1$ in $\overline\D$ to the unique point
$c\in\overline H$ such that $f_c$ has a period $n$ orbit of multiplier
$\lambda$. This homeomorphism extends holomorphically over a neighborhood
of $\overline\D$, with just one critical point
$1\in\overline\D$ mapping to the root point $\hat c\in H$ in the primitive
case, and with no critical points in the satellite case.
The closures of the various hyperbolic components
are pairwise disjoint, except for the tangential contact between a
component and its satellite as described in $6.1$.\ms}

{\tf Proof.} Recall that $\lambda_n:\Per_n/\Z_n\to\C$ is a proper holomorphic
map of degree $\nu_2(n)/2$. We will first show that there are no critical
values of $\lambda_n$ within the closed unit disk $\overline\D$. This will
imply that the inverse image $\lambda_n^{-1}(\overline\D)$ is the disjoint
union of $\nu_2(n)/2$ disjoint sets $\overline H^\natural$, each of which
maps diffeomorphically onto $\overline\D$. First note that there are no
critical values of $\lambda_n$ on the boundary circle $\partial\D$. In the case
of a root of unity $\mu\in\partial\D$, every $(c,\Or)$ with $\lambda_n(c,\Or)=
\mu$ must be parabolic, and it follows from 6.1 and 6.2 that the derivative of
$\lambda_n$ at $(c,\Or)$ is non-zero.
Consider then a point $(\hat c,\Or)\in\partial H^\natural$ such that
$\lambda_n(\hat c,\Or)$ is not a root of unity. According to 5.2, we can
use $c$ as local uniformizing parameter throughout a neighborhood of
$(\hat c,\Or)$. If this were a critical point of $\lambda_n$, then it would
follow that we could find two different line segments emerging from $\hat c$
which map into $\D$, separated by two line segments which map outside of $\D$.
In other words, one of the following two possibilities would have to occur.

{\bf Case 1.} There are two different hyperbolic components
with $\hat c$ as non-root boundary point. Each of these
components must have a root point, and be contained in its associated wake.
But these two components cannot be separated by any rational parameter ray,
hence each one must be contained in the wake of the other, which is impossible.

{\bf Case 2.} The single hyperbolic component $H$ must approach $\hat c$ from
two different directions, separated by two directions which lie outside of $H$.
In other words. There must be a simple closed loop $L\subset\overline H$ which
encloses points lying outside of $\overline H$. Now the collection of
iterates $f_c^{\circ k}(0)$ must be uniformly bounded for $c\in L$, and hence
also for all $c$ in the region bounded by $L$. Thus this entire region must
lie within the interior of the Mandelbrot set, which is impossible since
this region contains parabolic points.

Thus both cases are impossible, and $\lambda_n$ must be locally injective
near the boundary of $H^\natural$. It follows easily that $H^\natural$ maps onto
$\D$ by a proper map of some degree $d\ge 1$, and similarly that the boundary
$\partial H^\natural$ wraps around the boundary circle $\partial\D$ exactly
$d$ times. Now a counting argument shows that this degree is $+1$.
In fact the number of $H$ or $H^\natural$ of period $n$
is equal to $\nu_2(n)/2$ by 6.3 and 5.3.
Since the degree of the map $\lambda_n$ on $\Per_n/\Z_n$ is also
$\nu_2(n)/2$ by 5.2, it follows that each $H^\natural$ must map with degree
$d=1$. Therefore $\lambda_n$ maps each $\overline H^\natural$ biholomorphically
onto $\overline\D$.

Next consider the projection $(c,\Or)\mapsto c$ from the compact
set $\overline
H^\natural$ onto $\overline H$. This is one-to-one, and hence a homeomorphism,
by a theorem of Douady and Hubbard which asserts that a polynomial of degree
$d$ can have at most $d-1$ non-repelling cycles. (Compare [Sh1]. Alternatively,
it follows from the classical Fatou-Julia theory that a polynomial with one
critical point can have at most one attracting cycle. If two
distinct points of $\partial H^\natural$ mapped to a single point of
$\partial H$, then, as in Case 2 above, a path between these
points in $H^\natural$ would map to a loop in $\overline H$ which could enclose
no boundary points of $H$, leading to a contradiction.)

According to 5.2, the parameter $c$ can be used as local uniformizing
parameter for $\Per_n/\Z_n$ unless $\lambda_n=1$. Hence the
only possible critical value for the projection $\overline H^\natural
\to\overline H$ is the root point. In fact, by 6.1 and 6.2, the root
point is actually a critical value if and only if $H$ is a primitive
component.

Finally suppose that two different hyperbolic components have a common
boundary point. If this boundary point is parabolic, then one of these
components must be a
satellite of the other by 6.1 and 6.2. If the point
were non-parabolic, then the
argument of Case 1 above would yield a contradiction.
This completes the proof of 6.5.\QED

\section{Orbit Forcing.} Recall that an orbit portrait
is {\bit non-trivial\/} if either it has valence $v\ge 2$, or it is the
zero portrait $\{\{0\}\}$.
The following statement follows easily from 1.3.
However, it seems of interest to give a direct and more constructive proof;
and the methods used will be useful in the next section.

{\ssk\nin{\tf Lemma 7.1. Orbit Forcing.} \it Let $\po$ and $\q$ be distinct
non-trivial orbit portraits. If their characteristic arcs satisfy
$I(\po)\subset I(\q)$, then every $f_c$ with a (repelling or parabolic)
orbit of
portrait $\po$ must also have a repelling orbit of portrait $\q$.\ssk}

Compare Figure 5, and see 1.3 and for further discussion. The proof of 7.1
begins as follows.\ssk

{\bf Puzzle Pieces.}
Recall from 2.10 that the $pv$ rays landing on a periodic orbit for $f=f_c$
separate the dynamic plane into $pv-p+1$ connected components, the closures
of which are called the (unbounded) {\bit preliminary
puzzle pieces\/} associated with the given orbit portrait. (As in [K],
we work with puzzle pieces which are closed but not compact. The associated
{\bit bounded\/} pieces can be obtained by intersecting each unbounded puzzle
piece with
the compact region enclosed by some fixed equipotential curve.)

Most of these preliminary puzzle pieces $\Pi$ have the {\bit Markov property\/}
that $f$ maps $\Pi$ homeomorphically onto some union of preliminary
puzzle pieces. However, the puzzle piece containing the critical point
is exceptional: Its image under $f$ covers the critical value puzzle
piece twice, and also covers some further puzzle pieces once.
To obtain a modified puzzle with more convenient properties, we will subdivide
this exceptional piece into two connected sub-pieces.

Let $\Pi_1$ be the preliminary puzzle piece containing the
critical value. Then $\partial\Pi_1$
consists of the two rays whose angles bound the characteristic arc for
$\po$, together with their common landing point, say $z_1$.
The pre-image $\Pi_0=f^{-1}(\Pi_1)$ is bounded by two rays landing at
the point $z_0=f^{-1}(z_1)\cap\Or$, together with two rays landing at
the symmetric point $-z_0$. Note that $\Pi_0$ is a connected set containing
the critical point, and that the map $f$ from $\Pi_0$ onto $\Pi_1$
is exactly two-to-one, except at the critical point $0$, which maps to $c$.

The $pv$ rays landing on $\Or$, together with these two additional
rays landing on $-z_0$, cut the complex plane up into $pv-p+2$ closed
subsets which we will call the pieces of the corrected {\bit puzzle\/}
associated with $\po$. These will be numbered as
$\Pi_0\,,\,\Pi_1\,,\,\ldots\,,\,\Pi_{pv-p+1}$, with $\Pi_0$ and $\Pi_1$
as above. The central piece $\Pi_0$ will be called
the {\it critical puzzle piece\/}, and
$\Pi_1$ will be called the {\it critical value puzzle piece\/}.
This corrected puzzle satisfies the following.

{\ssk\nin{\bf Modified Markov Property.}
\it The puzzle piece $\Pi_0$ maps onto $\Pi_1$ by a 2-fold branched
covering, while
every other puzzle piece maps homeomorphically onto a finite union of
puzzle pieces.\ssk}

We can represent the allowed
transitions by a {\bit Markov matrix\/} $M_{ij}$, where
$$	M_{ij}~=~\begin{cases} 1 & \text{if}~~ \Pi_i~~
 \text{maps homeomorphically, with~~}
 f(\Pi_i)\supset\Pi_j\cr
	0 & \text{if}~~ f(\Pi_i)~~ \text{and}~~ \Pi_j~~
 \text{have no interior points in common,}
\end{cases}$$
and where $M_{01}=2~$ since $\Pi_0$ double covers $\Pi_1$.
Since $f$ is quadratic, note that the sum of entries  in any column is equal
to $2$. Equivalently, this same data can be represented by a {\bit Markov
graph\/}, with one vertex for each puzzle piece, and with $M_{ij}$
arrows from the $i$-th vertex to the $j$-th.

\begin{figure}[htb]
\cl{\psfig{figure=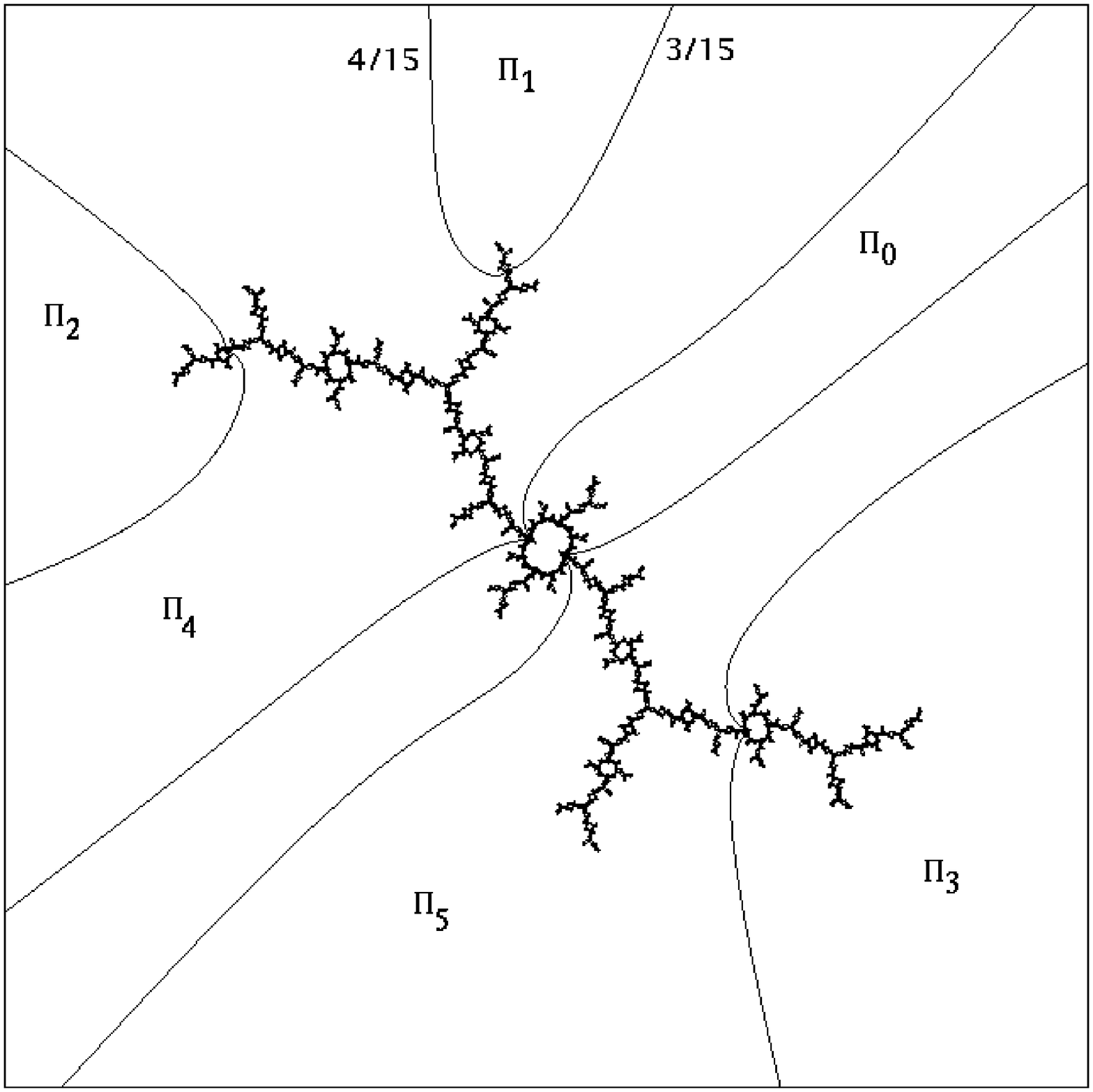,height=2.5in}\qquad
	\psfig{figure=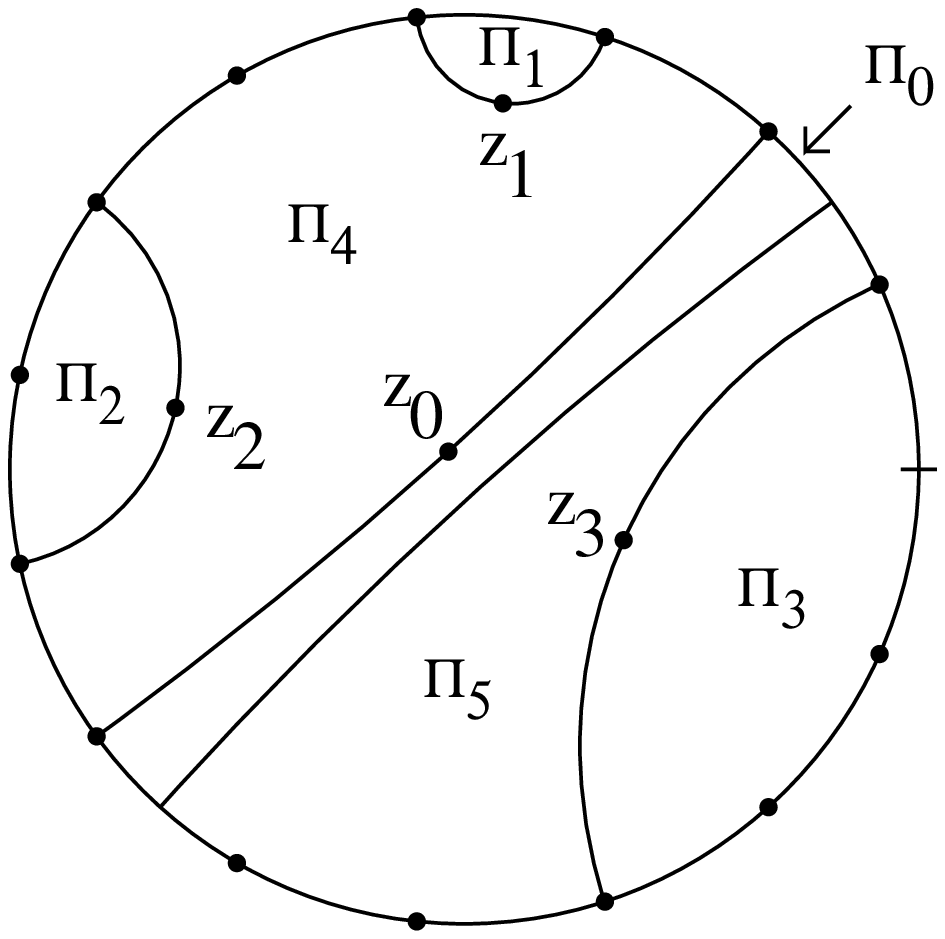,height=2.2in}}

\begin{quote}{\bit Figure 13. Julia set with a parabolic orbit of period four
with characteristic arc $I(\po)=(3/15\,,\,4/15)$, showing the six
corrected puzzle pieces; and a corresponding schematic diagram. (For the
corresponding preliminary puzzle, see the top of Figure 5.)}\end{quote}
\end{figure}
\begin{figure}[htb]
\cl{\psfig{figure=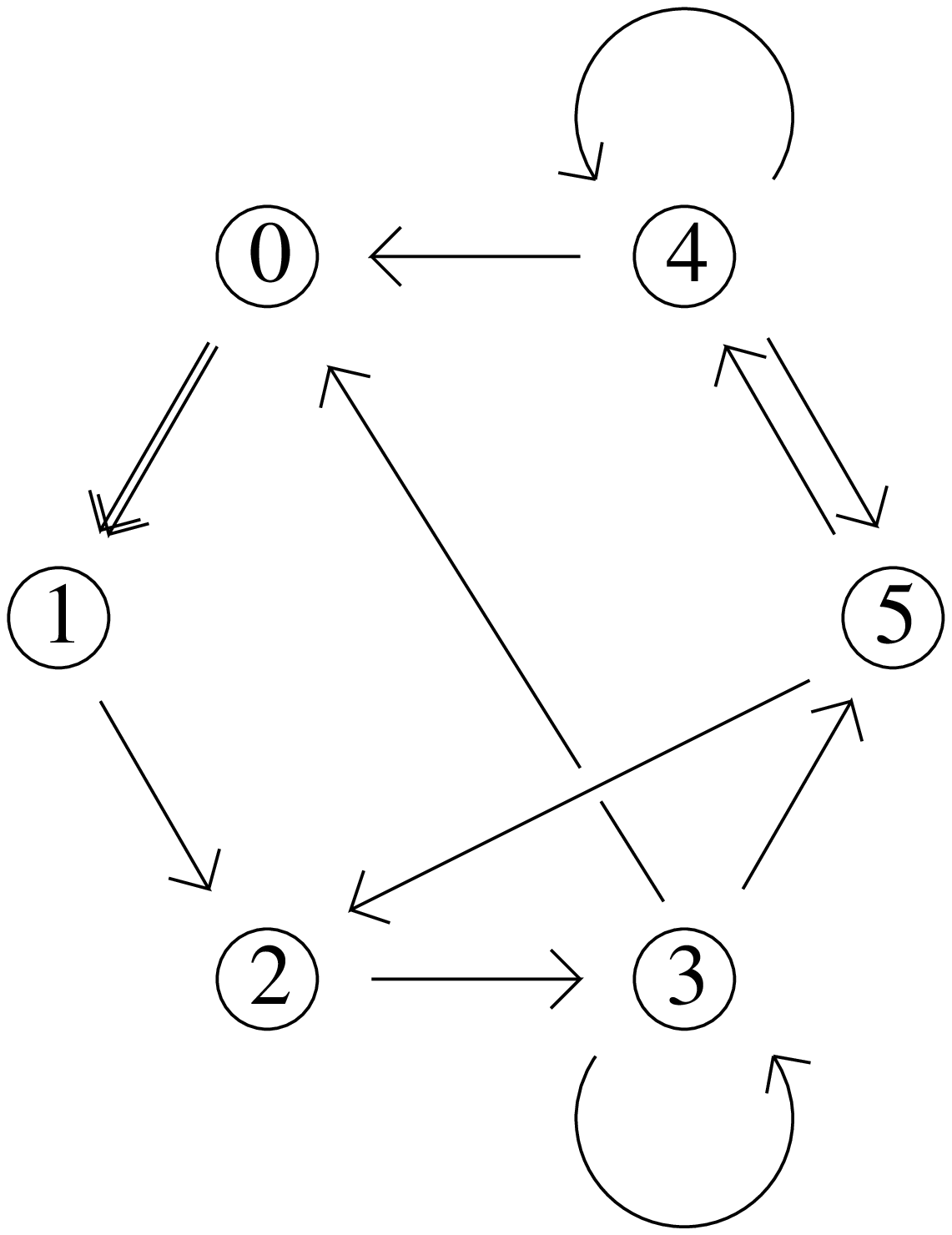,height=2.2in}}\ssk
\begin{quote}{\bit Figure 14.  Markov graph associated with the matrix
$(8)$, with one vertex for each
puzzle piece. Since $f$ is quadratic,
there are two arrows pointing to each vertex.\vskip -.4in}
\end{quote}\end{figure}

As an example, for the puzzle shown in Figure 13, we obtain the
Markov graph of Figure 14, or the following Markov matrix
$$ \big[M_{ij}\big]~=~\left[\begin{matrix}
	0&2&0&0&0&0\cr
	0&0&1&0&0&0\cr
	0&0&0&1&0&0\cr
	1&0&0&1&0&1\cr
	1&0&0&0&1&1\cr
	0&0&1&0&1&0\cr
\end{matrix}\right]~. \eqno(8)$$

\ssk
To illustrate the idea of the proof of 7.1,
let us show that any $f$ having an
orbit with this portrait $\po$ must also have a repelling orbit with portrait
$\q=\{\{1/7,2/7,4/7\}\}$. (Compare the top implication in Figure 5.)
Inspecting the next to last row of the matrix (8), we see that $f(\Pi_4)
=\Pi_0\cup\Pi_4\cup\Pi_5$.
Therefore, there is a branch $g$
of $f^{-1}$ which maps the interior of $\Pi_4$
holomorphically onto some proper subset of itself. This mapping $g$ must
strictly decrease the Poincar\'e metric for the interior of $\Pi_4$. On
the other
hand, it is easy to check that the $1/7\,,\,2/7$ and $4/7$ rays are all
contained in the interior of $\Pi_4$. Hence their landing points, call them
$w_1\,,\,w_2$ and $w_3$, are also contained in $\Pi_4$, necessarily in the
interior, since the points of $K\cap\partial\Pi_4$ have period four. Now
$$	g~:~w_1~\mapsto w_3~\mapsto~w_2~\mapsto~w_1~,$$
and all positive distances are strictly decreased. Thus if the distance
from $w_i$ to $w_j$ were greater than zero, then applying $g$
three times we would obtain a contradiction. This proves that $w_1=w_2=w_3$,
as required. This fixed point must be repelling, since $g$ clearly cannot
be an isometry.

A similar argument proves the following statement. Suppose that $f=f_c$
has an orbit $\Or$ with some given portrait $\po$. By a {\bit Markov
cycle\/} for $\po$ we will mean an infinite sequence
of non-critical puzzle pieces $\Pi_{i_1}\,,\,\Pi_{i_2}\ldots$
which is periodic, $i_j=i_{j+m}$ with period $m\ge 1$, and which satisfies
$f(\Pi_{i_\alpha})\supset\Pi_{i_{\alpha+1}}$, so that
$M_{i_\alpha\,i_{\alpha+1}}=1$, for every $\alpha$ modulo $m$.

{\ssk\nin{\tf Lemma 7.2. Realizing Markov Cycles.} \it Given such a Markov
cycle, there is one and only one periodic orbit
$z_1\mapsto\cdots\mapsto z_m$ for $f_c$ with period dividing $m$ so that
each $z_\alpha$ belongs to $\Pi_{i_\alpha}$,
and this orbit is necessarily repelling unless it coincides with the given
orbit $\Or$ (which may be parabolic). In particular, for any
angle $t$ which is periodic under doubling, if the dynamic ray
with angle $2^\alpha\,t$ lies in $\Pi_{i_\alpha}$ for all integers
$\alpha$, then this ray must land at the point $z_\alpha$.\ssk}

(Note that the period of $t$ may well be some multiple of $m$, as in the
example just discussed.)

{\tf Proof.} There is a unique branch of $f_c^{-1}$ which carries the interior
of $\Pi_{i_{\alpha+1}}$ holomorphically onto a subset of $\Pi_{i_\alpha}$.
Let $g_{i_\alpha}$ be the composition of these $m$ maps, in the appropriate
reversed order so as to carry the interior of $\Pi_{i_\alpha}$ into itself.

A similar construction applies to the associated external angles. Let $J_i
\subset\R/\Z$ be the set of all angles of dynamic rays which are contained
in $\Pi_i$. Thus each $J_i$ is a finite union of closed arcs,
and together the $J_i$ cover $\R/\Z$ without overlap. Now there is
a unique branch of the 2-valued map $t\mapsto t/2$ which carries
$J_{i_{\alpha+1}}$ into $J_{i_\alpha}$ with derivative $1/2$
everywhere. Taking an $m$-fold composition, we map each $J_{i_\alpha}$
into itself with derivative $1/2^m$. This composition may well permute
the various connected components of $J_{i_\alpha}$. However, some iterate
must carry some component of $J_{i_\alpha}$ into itself, and hence have a
unique fixed point $t$
in that component. The landing point of the corresponding
dynamic ray will be a periodic point $z_\alpha\in\Pi_{i_\alpha}$.

{\bf Case 1.}
If this landing point belongs to the interior of $\Pi_{i_\alpha}$,
then it is fixed by some iterate of our map $g_{i_\alpha}$. This
map $g_{i_\alpha}$ cannot be an isometry, hence it must contract the
Poincar\'e metric. Therefore
every orbit under $g_{i_\alpha}$ must converge towards
$z_\alpha$. Thus $z_\alpha$ is an attracting fixed point
for $g_{i_\alpha}$, and hence
is a repelling periodic point for $f$.

{\bf Case 2.} If the landing point belongs to the boundary of
$\Pi_{i_\alpha}$ then it must belong to $\Or\cup\{-z_0\}$, and hence
to the original orbit $\Or$ since $-z_0$ is not periodic. Evidently
this case will occur only when the angle $t$ belongs to the union
$A_1\cup\cdots\cup A_p$ of angles in the given portrait $\po$. \QED

{\tf Note.} It is essential for this argument that our given Markov cycle
$\{\Pi_{i_\alpha}\}$ does not involve the critical puzzle piece $\Pi_0$.
In fact, as an immediate corollary we get the following statement:

{\ssk\nin{\tf Corollary 7.3. Non-Repelling Cycles.} \it Any non-repelling
periodic orbit for $f$ must intersect the critical puzzle piece $\Pi_0$
as well as the critical value puzzle piece $\Pi_1$.}
\ssk

{\tf Proof of 7.1.} If $I(\po)\subset I(\q)$, then it follows from Lemma
2.9 that there exists a map $\hat f$ having  both an orbit with portrait $\po$
and an orbit with portrait $\q$. The latter orbit determines a Markov
cycle in the puzzle associated with $\po$. (The condition $I(\po)\subset I(\q)$
guarantees that this cycle avoids the critical puzzle piece.)
Now for any map $f$ with
an orbit of portrait $\po$, we can use this Markov cycle, together with
7.2, to construct the required periodic orbit and to guarantee that
the rays associated with the portrait $\q$ land on it, as required.\QED

In fact an argument similar to the proof of 7.2 proves a much sharper
statement. Let $\Or$ be a repelling periodic orbit with non-trivial portrait
$\po$.

{\ssk\nin{\tf Lemma 7.4. Orbits Bounded Away From Zero.}
\it Given an infinite sequence of non-critical puzzle
pieces $\{\Pi_{i_k}\}$ for $k\ge 0$ with $f(\Pi_{i_k})\supset
\Pi_{i_{k+1}}$, there is one and only one point $w_0\in K(f)$ so that the
orbit $w_0\mapsto w_1\mapsto\cdots$ satisfies $w_k\in\Pi_{i_k}$ for
every $k\ge 0$. It follows that the action of $f$
on the compact set $K_\po$ consisting of
all $w_0\in K(f)$ such that the forward orbit $\{w_k\}$ never hits the
interior of $\Pi_0$ is topologically
conjugate to the one-sided subshift  of finite type, associated
to the matrix $[M_{ij}]$ with $0$-th row and column deleted. In particular,
the topology of $K_\po$ depends only on $\po$, and not on the particular
choice of $f$ within the $\po$-wake.\ssk}

{\tf Proof Outline.} First replace each puzzle piece
$\Pi_i$ by a slightly thickened puzzle piece, as described in [M3].
(Compare \S8, Figure 18.)
The interior of this thickened piece is an open neighborhood
$N_i\supset\Pi_i$, with the property that
$f(N_i)\supset \overline N_j$ whenever $f(\Pi_i)\supset\Pi_j$.
It then follows that there is a branch of $f^{-1}$ which maps $N_j$ into
$N_i$, carrying $K\cap\Pi_j$ into $K\cap\Pi_i$, and
reducing distances by at least some fixed ratio $r<1$
throughout the compact set $K\cap\Pi_j$. Further details
are straightforward.\QED

Presumably this statement remains true for a parabolic orbit,
although the present proof does not work in the parabolic case.
(Compare [Ha].) 

\section{Renormalization.} One remarkable property of the Mandelbrot
boundary is that it is densely filled with small copies of itself. (See
Figures 11, 14 for a magnified picture of one such small copy.)
This section will provide a rough outline, without proofs,
of the Douady-Hubbard theory of renormalization, or the inverse operation of
tuning, which provides a dynamical
explanation for these small copies. It is based on
[D4] as well as [DH3], [D3]. (Compare [D1], [M1]. For the Yoccoz interpretation
of this construction, see [Hu], [M3], [Mc], [Ly].
For a more general form of renormalization, see [Mc], [RS].)\ssk

\begin{figure}[htb]
\cl{\psfig{figure=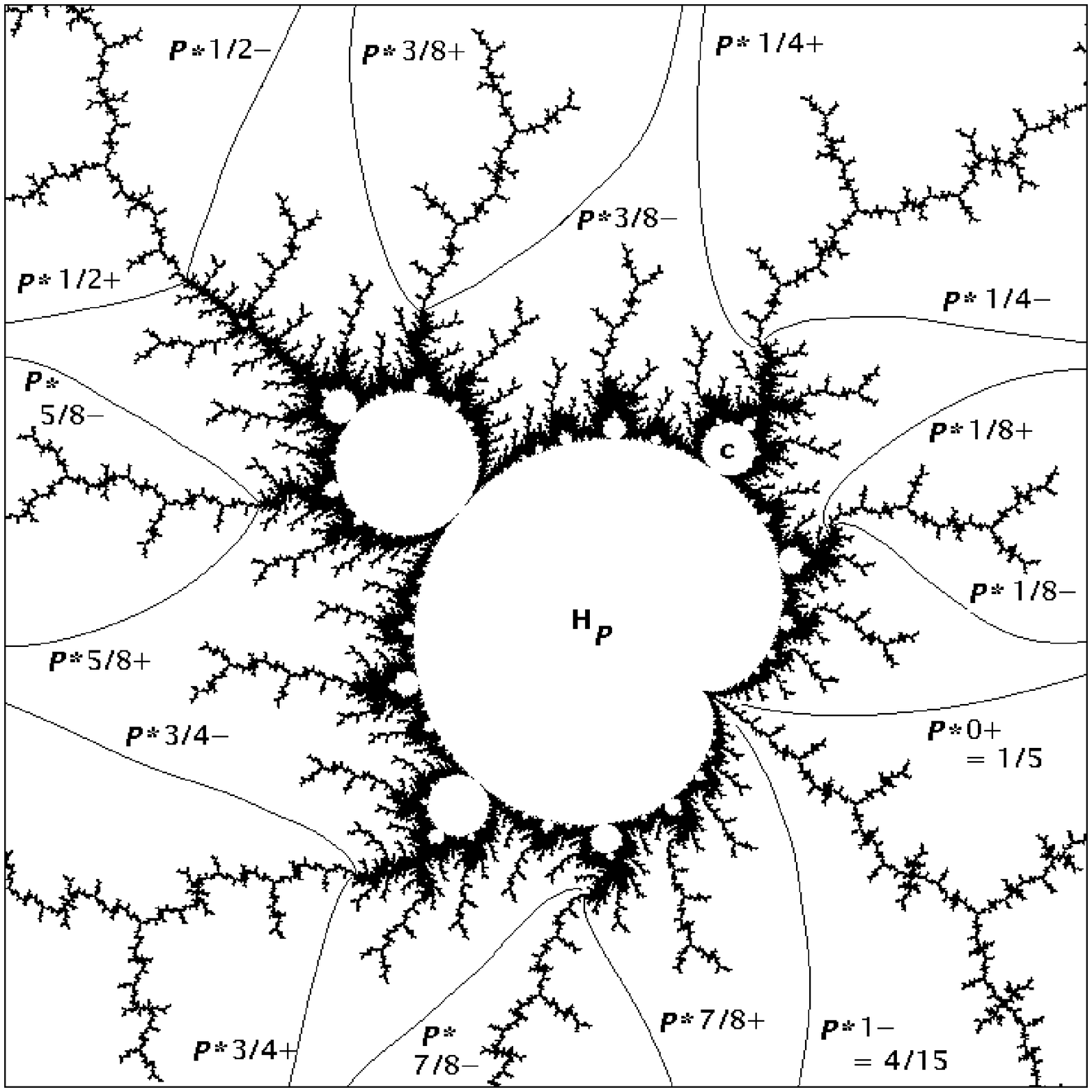,height=3.5in}}\ssk
\begin{quote}{\bit Figure 15. Detail near the period 4 hyperbolic component
$H_\po$ of Figure 12, where $\po=\po(1/5,4/15)$,
showing the first eight of the parameter
sectors which must be pruned away from $M$
to leave the small Mandelbrot set consisting of\break $\po$-renormalizable
parameter values.}\end{quote}
\end{figure}

To begin the construction, consider any orbit portrait $\po$ of ray period
$n\ge 2$ and valence $v\ge 2$. Let $c$ be a parameter value
in $W_\po\cup\{r_\po\}$, so that $f=f_c$ has a periodic orbit $\Or$
with portrait $\po$, and let $S=S(f)$ be the critical value
sector for this orbit (so that $\overline S$
is the critical value puzzle piece). To a first
approximation, we could try to say that $f$ is ``$\po$-renormalizable'' if
the orbit of $c$ under $f^{\circ n}$ is completely
contained in $S$. In fact this is a
necessary and sufficient condition whenever the map $f^{\circ n-1}|_S$
is univalent. However, in examples such as that of Figures 1, 2 one needs a
slightly sharper condition.

Let $\I_\po=(t_-,t_+)$ be the characteristic arc for this portrait,
so that $\partial S$ consists of the dynamic rays of angle $t_-$ and $t_+$
together with their common landing point $z_1$,
and let $\ell=t_+-t_-$ be the length of this arc.


\begin{figure}[htb]
\vskip .2in
\cl{\psfig{figure=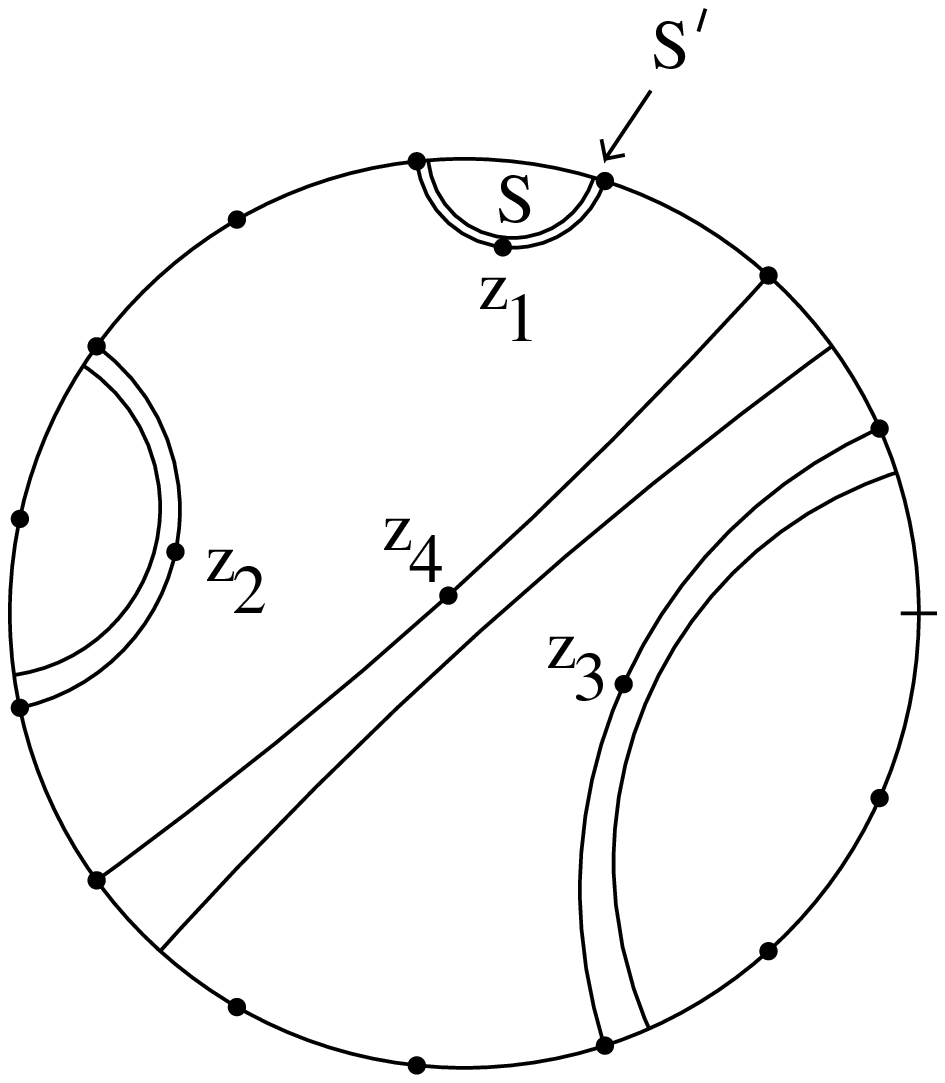,height=2.5in}}
\ssk\begin{quote}
{\bit Figure 16. The $n$-fold pull-back of the critical sector $S$
along the orbit $\Or$,
illustrated schematically for the orbit diagram of period $n=4$ which has
characteristic arc $(1/5,4/15)$. Compare Figures 5 (top), 12, 16.}
\end{quote}\end{figure}

{\ssk\nin{\tf Lemma 8.1. A (Nearly) Quadratic-Like Map.}
\it The dynamic rays of angle $t'_1=t_-+\ell/2^n$ and $t'_2=
t_+-\ell/2^n$ land at a common point $z'\ne z_1$ in $S\cap f^{-n}(z_1)$.
Let $S'\subset S$ be the region
bounded by $\partial S$ together with these two rays and their common
landing point. Then the map $f^{\circ n}$ carries $S'$ onto $S$ by
a proper map of degree two, with critical value equal to the critical
value $f(0)=c$.\ssk}

This region $S'$ can be described as the $n$-fold ``pull-back'' of $S$
along the orbit $\Or$. (Compare Figure 16, which
also shows the first three forward images of $S'$.)

{\tf Proof of 8.1.} First suppose that $c\in W_\po$ is outside the Mandelbrot
set. Then, following Appendix A, we can bisect the complex
plane by the two rays leading from infinity to the critical point.
(Compare the proof of 2.9.)
In order to check that the two rays of angle $t'_1$ and $t'_2$
have a common landing point, we need only show that they have
the same symbol sequence with respect to the resulting partition.
In other words, we must show, for every $k\ge 0$, that the $2^kt'_1$ and
$2^kt'_2$ rays
lie on the same side of the bisecting critical ray pair. For $k\ge n$
this is clear since $2^nt'_1\equiv t_+$ and $2^nt'_2\equiv t_-$ modulo $\Z$.

Now consider the critical puzzle piece $\Pi_0$ of \S7. Evidently $\Pi_0$
is a neighborhood, of angular radius $\ell/4$, of the bisecting critical
ray pair. For $k<n-1$ the dynamic rays with angle $2^kt_-$ and $2^kt_+$ both
lie in the same component of $\C\ssm\Pi_0$. Since $2^kt'_j$ differs from
$2^kt_j$ by at most $\ell/4$, it follows that the $2^kt'_1$ and $2^kt'_2$
rays have the same symbol. Finally, for $k=n-1$, it is not difficult to
check that the $2^kt'_1$ and $2^kt'_2$ rays both land at the same point $-z_0
\ne z_0$. This proves that the $t'_1$ and $t'_2$ rays land at the same point,
different from $z_1$, when $c\not
\in M$. A straightforward continuity argument now proves the same
statement for all $c\in W_\po$.

Thus we obtain the required region $S'\subset S$.
As in \S2, it will be convenient to complete the complex plane by adjoining
a circle of points at infinity. Note that the boundary of $S'$ within this
circled plane $\cc$ consists of two arcs of length $\ell/2^n$ at infinity,
together with two ray pairs and their common landing points. As we traverse this
boundary once in the positive direction, the image under $f^{\circ n}$
evidently traverses the boundary of $S$ twice in the positive direction.
Using the Argument Principle, it follows that the image of $S'$ is contained
in $S$, and covers every point of $S$ twice, as required. Thus
$f^{\circ n}|_{S'}$ must have exactly one critical point, which can only be
$c$.\QED\ssk

Thus we have an object somewhat like a quadratic-like map, as studied in
[DH3]. Note however that 
$S'$ is not compactly contained in $S$.\ssk

{\tf Definition.} We will say that $f$ is {\bit $\po$-renormalizable\/}
if $f(0)=c$ is contained in the closure
$\overline S'$, and furthermore the entire
forward orbit of $c$ under the map $f^{\circ n}$ is contained in
$\overline S'$. If this condition is satisfied, and the orbit of $c$ is also
bounded so that $c\in M$, then we will say that $c$ belongs to the
``{\bit small copy\/}'' $\po*M$ of the Mandelbrot set which is
associated with $\po$. (This terminology will be justified in 8.2. If the orbit
is unbounded, then we may say that $c$ belongs to a $\po$-{\bit renormalizable
external ray}.)


\begin{figure}[htb]
\cl{\psfig{figure=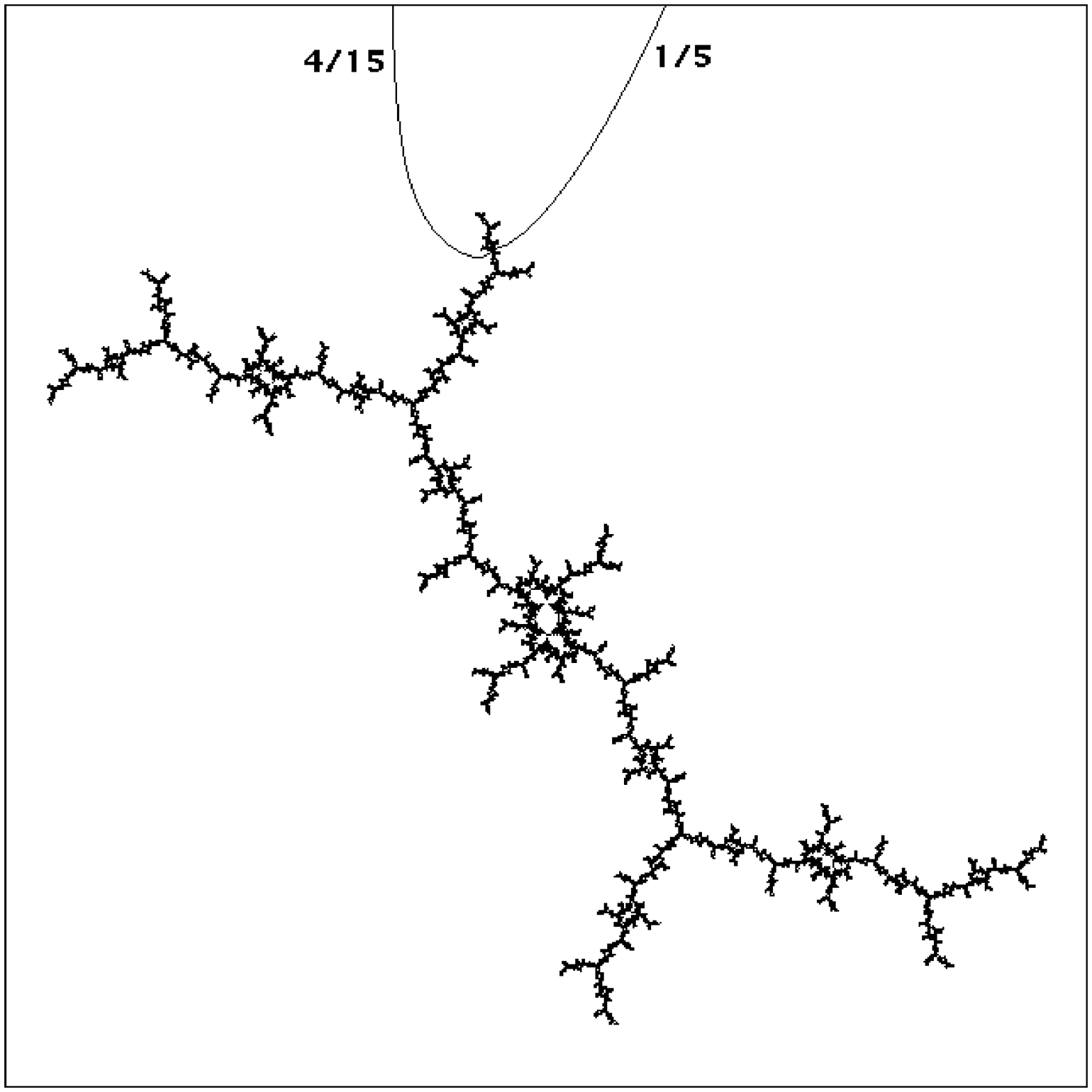,height=2.5in}\qquad
\psfig{figure=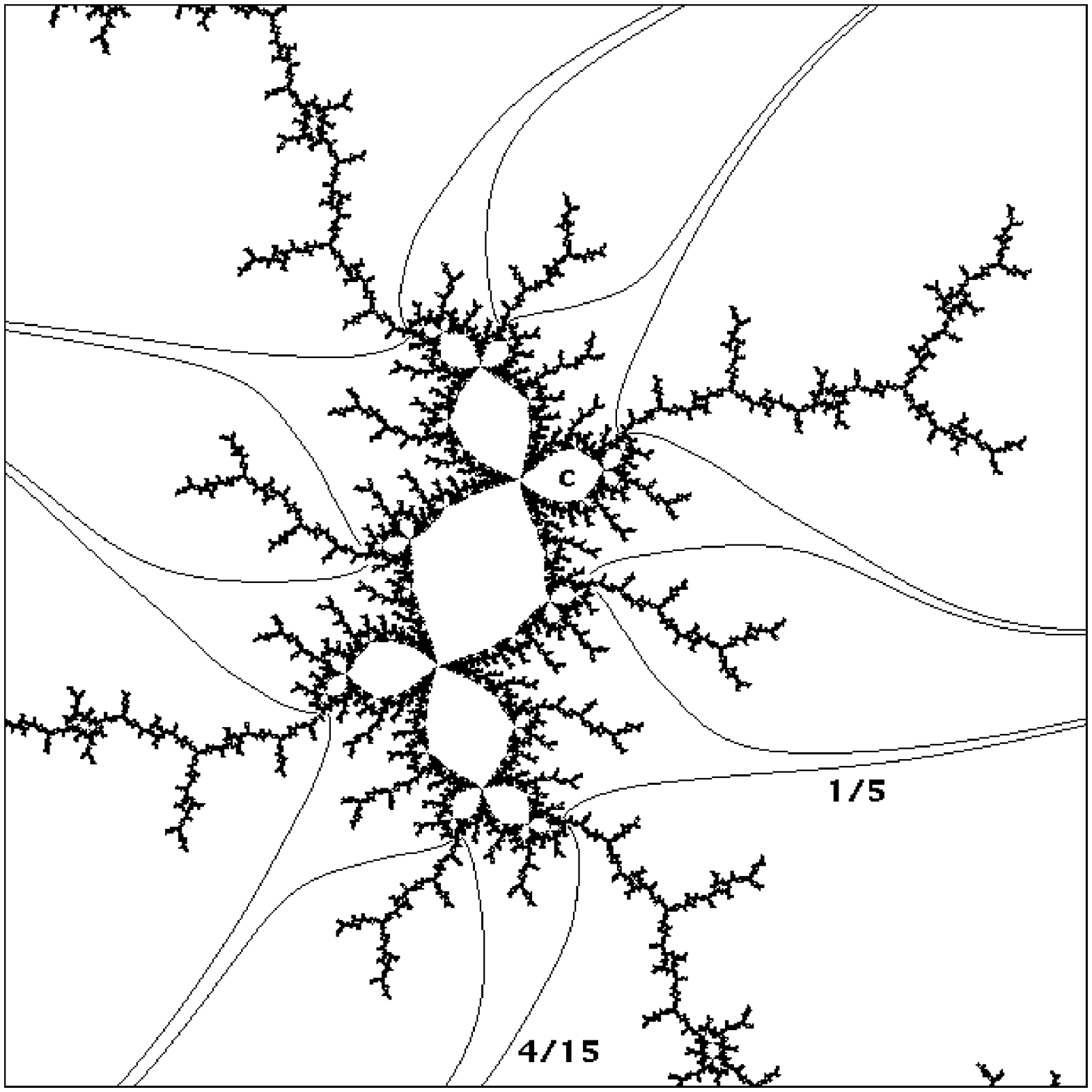,height=2.5in}}\ssk
\begin{quote}{\bit Figure 17. Julia set for the center point of the
period 12 satellite component of Figure 12 (the point $c$ of Figure 15), and a
detail near the critical value $c$, showing the first
eight of the sectors of the dynamic plane which must
be pruned away to leave the small Julia set associated with
$\po$-renormalization, with $\po$ as in Figures 12, 15.
(Here the right hand
figure has been magnified by a factor of 75.) This can be described as the
Julia set of Figure 13 (left) tuned by a ``Douady rabbit'' Julia set.
\ssk}\end{quote}
\end{figure}

Closely associated is the ``{\bit small filled Julia set\/}''
$K'=K(f^{\circ n}|S')$
consisting of all $z\in\overline S'$ such that the entire
forward orbit of $z$ under $f^{\circ n}$ is bounded and
contained in $\overline S'$.
(Compare Figure 17.) Thus the critical value $f(0)=c$ belongs to $K'$ if
and only if $f$ is\break $\po$-renormalizable, with $c\in M$.
As in the classical Fatou-Julia theory,
$c$ belongs to $K'$ if and only if $K'$ is connected.

\begin{figure}[htb]
\cl{\psfig{figure=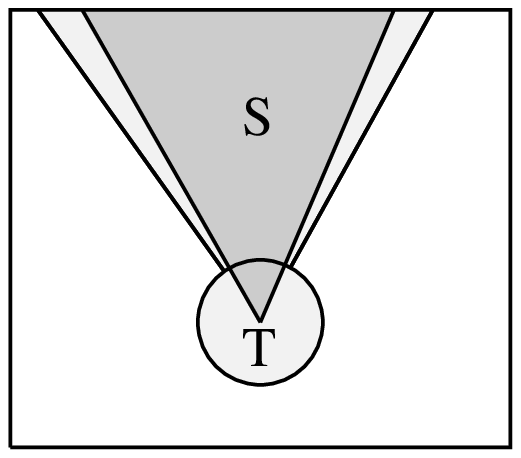,height=1.2in}}\ssk

\cl{\bit Figure 18. The sector $S$ and the thickened sector $T$.}
\end{figure}

In order to tie this construction up with Douady and Hubbard's
theory of polynomial-like
mappings, we need to thicken the sector $S$, and then cut it down to a
bounded set. (Compare [M3].) We exclude the exceptional
special case where $c$ is the root point $r_\po$. Thus we will suppose that
the periodic point $z_1\in\partial S$ is repelling. Choose a small disk
$D_\epsilon$ about $z_1$ which is mapped univalently by $f_c^{\circ n}$
and is compactly contained in $f_c^{\circ n}(D_\epsilon)$. Choose also
a very small $\eta>0$, and consider the dynamic rays with angle $t_--\eta$
and $t_++\eta$. Following these rays until they first meet $D_\epsilon$,
they delineate an open region $T\supset S\cup D_\epsilon$ in $\C$.
(Compare Figure 18.)
Now let $T'$ be the connected component of $f_c^{-n}(T)$ which contains
$S'$. It is not difficult to check that $\overline T'\subset T$, and that
$f_c^{\circ n}$ carries $T'$ onto $T$ by a proper map of degree two.

To obtain a bounded region, we let $U$ be the intersection of $T$ with
the set\break $\{z\in\C~;~G^K(z)< 1\}$, where $G^K$ is the Green's function
for $K=K(f_c)$. Similarly, let $U'$ be the intersection $T'$ with
$\{z~;~G^K(z)<1/2^n\}$. Then $U'$ is compactly contained in $U$, and
$f_c^{\circ n}$ carries $U'$ onto $U$ by a proper map of degree two.
{\it In other words, $f_c^{\circ n}|U'$ is a quadratic-like map.}

Evidently the forward orbit of a point $z\in U'$ under $f_c^{\circ n}$
is contained in $U'$ if and only if $z$ belongs to the small filled
Julia set $K'$. In particular, for $c\in M$, the map $f_c$ is
$\po$-renormalizable if and only if $c\in K'$, or if and only if $K'$
is connected.

If these conditions are satisfied, then according to [DH3]
the map $f_c^{\circ n}$ restricted to a neighborhood of $K'$ is
``hybrid equivalent'' to some uniquely defined quadratic map $f_{c'}$,
with $c'\in M$. Briefly, we will write $c=\po*c'$, or say that $c$
equals $\po$ {\bit tuned by\/} $c'$. Douady and
Hubbard show also that this correspondence
$$	c'~\mapsto \po*c' $$
is a well defined continuous embedding of $M\ssm\{1/4\}$ onto a
proper subset of itself. As an example, as $c'$ varies over the hyperbolic
component $H_{\{\{0\}\}}$ which is bounded by the cardioid, they show
that $\po*c'$ varies over the hyperbolic component $H_\po$.

It is convenient to supplement this construction, by defining
the operation\break
$\po\,,\,c'\mapsto \po*c'$ in two further special cases. If $c'$ is the
root point $1/4=r_{\{\{0\}\}}$ of $M$, then we define
$$	\po*(1/4)~=~r_\po $$
to be the root point of the $\po$-wake.
Furthermore, if $\po=\{\{0\}\}$ is the zero orbit portrait, then we define
$\{\{0\}\}*$ to be the identity map,
$$	\{\{0\}\}*c'~=~c' $$
for all $c'\in M$. With these definitions, we have the following basic
result of Douady and Hubbard.\ssk

{\ssk\nin{\tf
Theorem 8.2. Tuning.} \it For each non-trivial orbit portrait $\po$, the
correspondence $c\mapsto \po*c$ defines a continuous embedding of the
Mandelbrot set $M$ into itself. The image of this embedding is just the
``small Mandelbrot
set'' $\po*M\subset M$ described earlier. Furthermore, there is a unique
composition operation $\po\,,\,\q\mapsto\po*\q$ between non-trivial orbit
portraits so that the associative law is valid,
$$	(\po*\q)*c~=~\po*(\q*c) $$
for all $\po\,,\,\q$ and $c'$. Under this $*$ composition operation,
the collection of all non-trivial orbit portraits forms a free
(associative but noncommutative) monoid, with
the zero orbit portrait as identity element.\ms}

The proof is beyond the scope of this note.

We can better understand this construction by introducing a nested
sequence of open sets
$$	S=S^{(0)}~\supset ~S'=S^{(1)}~\supset~ S^{(2)}~\supset~\cdots $$
in the dynamic plane for $f$,
where $S^{(k+1)}$ is defined inductively as $S^{(k)}\cap f^{- n}(S^{(k)})$
for $k\ge 1$. Thus $S=S^{(0)}$ is bounded by the dynamic rays of angle
$t_-$ and $t_+$, together with their common landing point $z_1$.
Similarly, $S^{(1)}$ is bounded by $\partial S^{(0)}$ together with
the rays of angle $t_-+\ell/2^n$ and $t_+-\ell/2^n$,
together with their common landing point, which is an
$n$-fold pre-image of $z_1$. If $c\in S^{(1)}$, so that $S^{(2)}$
is a 2-fold branched covering of $S^{(1)}$, then
$S^{(2)}$ has two further boundary
components, namely the rays of angle $t_-+\ell/2^{2n}$ and
$t_-+\ell/2^n-\ell/2^{2n}$ and their common landing point, together
with the rays
of angle $t_+-\ell/2^n+\ell/2^{2n}$ and $t_+-\ell/2^{2n}$ and their common landing
point, for a total of 4 boundary components. Similarly, if $c\in S^{(2)}$,
then $S^{(3)}$ has 8 boundary components, as illustrated in Figure 17.

The angles which are left, after we have cut away the angles in
all of these (open) sectors, form a standard middle fraction
Cantor set ${\mathcal K}$, which can be described as follows. Let
$\phi$ be the fraction $1-2/2^n$. Start with the closure $[t_-\,,\,t_+]$ of the
characteristic arc for $\po$, with length $\ell$. First remove the
open middle segment of length $\phi\,\ell$, leaving two arcs of
length $\ell/2^n$. Then, from each of these two
remaining closed arcs, remove the middle segment of length
$\phi\,\ell/2^n$, leaving four segments of length $\ell/2^{2n}$, and
continue inductively. The intersection of all of the
sets obtained in this way is the required
Cantor set ${\mathcal K}\subset [t_-\,,\,t_+]$
of angles. These are precisely the angles of the dynamic rays which land
on the small Julia set $\partial K'$ (at least if we assume that these
Julia sets are locally connected).

There is a completely analogous construction in parameter space, as illustrated
in Figure 15. As noted earlier, parameter rays of angle $t_-$ and $t_+$ land
on a common point $r_\po$, and together form the boundary of the $\po$-wake.
Similarly, the parameter rays of angle $t_-+ \ell/2^n$ and $t_+-\ell/2^n$
must land at a common
point. These rays, together with their landing point, cut $W_\po$ into two
halves. For $c$ in the inner half, with boundary point $r_\po$, the
critical value of $f_c$ lies in $S'=S^{(1)}$, while for $c$ in the
outer half, this is not true. Similarly, for each pair of dynamic rays with
a common landing point in $\partial K$, forming part of
the boundary of $S^{(k)}$, there is a pair of parameter rays
with the same angles which have a common landing point in $\partial M$ and
form part of the boundary of a corresponding region $W_\po^{(k)}$
in parameter space. The basic property is that $c\in W_p^{(k)}$ if and only
if $c$ belongs to the corresponding region $S^{(k)}$ in the dynamic
plane for $f_c$.

Dynamically, the Cantor set ${\mathcal K}\subset\R/\Z$ can be described as the
set of angles in
$$	[t_-\;,\;t_-+\ell/2^n]\;\cup\;[t_+-\ell/2^n\;,\;t_+] $$
such that
the entire forward orbit under multiplication by $2^n$ is contained in
this set. Evidently the resulting dynamical system is topologically
isomorphic to the one-sided two-shift. Thus each element $t\in{\mathcal K}$
can be coded by an infinite sequence $(b_0\,,\,b_1\,,\,\ldots)$ of bits,
where each $b_k$ is zero or one according as $2^{nk}t$ belongs to the left
or right subarc. We will write $t=\po*(b_0b_1b_2\cdots)$. Intuitively, we can
identify this sequence of bits $b_i$ with the angle
$~.b_0b_1b_2\cdots=\sum b_k/2^{k+1}$. However, some care is needed since
the correspondence $.b_0b_1b_2\cdots\mapsto
\po*(b_0b_1\cdots)$has a jump discontinuity
at every dyadic rational angle, i.e., at those angles corresponding
to gaps in the Cantor set $\mathcal K$. Thus
we must distinguish between the left hand limit $\po*\alpha-$
and the right hand limit $\po*\alpha+$ when $\alpha$ is a dyadic rational.

With this notation, the angles of the bounding rays
for the various open sets $S^{(k)}$, or for the corresponding sets
$W_\po^{(k)}$ in parameter space, are just these left and right hand limits
$\po*\alpha\pm$, where $\alpha$ varies over the dyadic rationals; and
the composition operation between non-trivial orbit portraits
can be described as follows: {\it If $\q$ has characteristic arc
$(t_-,t_+)$, then $\po*\q$ has characteristic arc $(\po*t_-\,,\,
\po*t_+)$.\/} For further details, see [D3]. 

\section{Limbs and the Satellite Orbit.}
Let $\po$ be a non-trivial orbit
portrait with period $p\ge 1$ and ray period $rp\ge p$. (Thus $\po$ may be
either a primitive or a satellite portrait.) Recall that the {\bit limb\/}
$M_\po$ consists of all points which belong both to the Mandelbrot set $M$
and to the closure
$\overline W_\po$ of the $\po$-wake. By definition, a limb $M_\q$ with
$\q\ne\po$ is a {\bit satellite\/} of $M_\po$ if its root point ${\bf r}_\q$
belongs to the boundary of the associated hyperbolic component $H_\po$.
(See 6.4.)
We will prove the following two statements. (Compare [Hu], [S\o], [S3].)

{\ssk\nin{\tf Theorem 9.1. Limb Structure.} \it
Every point in the limb $M_\po$ 
either belongs to the closure $\overline
H_\po$ of the associated hyperbolic component, or else belongs
to some satellite limb $M_\q$.}\ssk

(For a typical example, see Figure 12.) For any parameter value $c$
in the wake $W_\po$, let $\Or(c)=\Or_\po(c)$ be the repelling orbit
for $f_c$ which has period $p$ and portrait $\po$. Clearly this
orbit $\Or(c)$ varies holomorphically with the parameter value $c$.

{\ssk\nin{\tf Corollary 9.2. The Satellite Orbit.}
\it To any $c\in W_\po$ there is associated another orbit
$\Or^\star(c)=\Or^\star_\po(c)$, distinct from $\Or(c)$, which has
period $n=rp$ and which also varies holomorphically
with the parameter value $c$. As $c$ tends to the root point
${\bf r}_\po$, the two orbits $\Or(c)$ and $\Or^\star(c)$ converge
towards a common parabolic orbit of portrait $\po$. $($Compare 4.1.$)$
This associated orbit $\Or^\star(c)$ is attracting if $c$ belongs to the
hyperbolic component $H_\po\subset W_\po$, indifferent for $c\in\partial
H_\po$, and is repelling
for $c\in W_\po\ssm\overline H_\po$, with portrait equal to $\q$
if $c$ belongs to the satellite wake $W_\q$. \ssk}

As an example, both statements apply to the zero portrait,
with $M_{\{\{0\}\}}$ equal to the entire Mandelbrot set, with
$W_{\{\{0\}\}}=\C\ssm(1/4,+\infty)$, and with
$H_{\{\{0\}\}}$ bounded by the cardioid. In this case,
for any $c\in W_{\{\{0\}\}}$, the orbit
$\Or(c)$ consists of the {\bit beta fixed point\/}
$(1+\sqrt{1-4c})/2$ while $\Or^\star(c)$ consists of the {\bit alpha fixed
point\/} $(1-\sqrt{1-4c})/2$,
taking that branch of the square root function with $\sqrt 1=1$.\ssk

{\tf Proof of 9.1.} For each $c\in H_\po$ let $\Or^\star(c)$ be the unique
attracting periodic orbit.
By the discussion in \S6, this orbit extends analytically as we vary $c$ over
some neighborhood of the closure $\overline H_\po$,
provided that we stay within the wake $W_\po$. Furthermore, this orbit
becomes strictly repelling as we cross out of $\overline H_\po$. Therefore
we can choose a neighborhood $N$ of $\overline H_\po$ which is small
enough so that this
analytically continued orbit $\Or^\star(c)$ will be strictly repelling for all
$c\in N\cap W_\po\ssm\overline H_\po$. If $c$ also belongs to the
Mandelbrot set, so that $c\in N\cap M_\po\ssm\overline H_\po$, it follows
that at least one
rational dynamic ray lands on the orbit $\Or^\star(c)$; hence there is an
orbit portrait $\q=\q(c)$ of period $n$ associated with $\Or^\star(c)$. Choosing the
neighborhood $N$ even smaller if necessary,
we will show that the rotation number of $\q(c)$
is non-zero, and hence that this portrait $\q(c)$
is non-trivial. In other words, we will prove that $c$ belongs to a
limb $M_\q$ which is associated to the orbit $\Or^\star(c)$.

First consider a point $\hat c$ which belongs to the boundary $\partial H_\po$.
Then $\Or^\star(\hat c)$ is an indifferent periodic orbit, with multiplier on the
unit circle.  Consider some dynamic ray $\ra_t^K$
which has period $n$, but does not participate in the portrait
$\po$, and hence does not land on the original
orbit $\Or(\hat c)$. Such a ray certainly cannot land on
$\Or^\star(\hat c)$, for that would imply that $\Or^\star(\hat c)$
was a repelling or parabolic orbit of rotation  number zero.
However, for $\hat c$
in the boundary of $H_\po$ the orbit $\Or^\star(\hat c)$ is never repelling, and is
parabolic of rotation number zero only when $\hat c$ is the root point of
$H_\po$, so that $\Or^\star(\hat c)=\Or(\hat c)$. Since we have assumed that the
ray $\ra_t^K$ does not land on $\Or(\hat c)$, it must land on some repelling
or parabolic periodic point which is disjoint
from $\Or^\star(\hat c)$. In fact it must land on a repelling orbit, since
a quadratic map cannot have a parabolic orbit and also a disjoint indifferent
orbit. (Compare \S6.) Now as we perturb $c$ throughout some neighborhood
of $\hat c$ it follows that the corresponding ray still lands on a repelling
periodic point disjoint from $\Or^\star(c)$. Since $\Or^\star(c)$ has period $n$,
but no ray of period $n$ can land on it, this proves that the rotation number
of the associated portrait $\q(c)$ is non-zero, as asserted.

Let $X$ be any connected component of $M_\po\ssm\overline H_\po$.
Since the Mandelbrot set is connected, $X$ must have some limit point in
$\partial H_\po$. Therefore, by the argument above, some point $c\in X$
must belong to a wake $W_\q$ associated with the orbit $\Or^\star(c)$.
Since the portrait $\q$ has period $n$,
the root point ${\bf r}_{\q}$ of its wake must lie on the boundary of some
hyperbolic component $H'$ which has period $n$ and is contained in $W_\po$.
In fact, for suitable choice of $c$,
we claim that $H'$ can only be $H_\po$ itself.
There are finitely many other components of period $n$, but these others
are all bounded away from $H_\po$, while the point $c\in X$ can be chosen
arbitrarily close to $\overline H_\po$. Thus we may assume that $W_\q$
is rooted at a point of $\partial H_\po$, and hence
is a satellite wake. Since the connected set X cannot
cross the boundary of $W_\q$, it follows that $X$ is completely contained
within $W_\q$, which completes the proof of 9.1\QED\ssk

{\tf Proof of 9.2.} As in the argument above, the orbit $\Or^\star(c)$ is well
defined for $c$ in some neighborhood of $W_\po\cap\overline H_\po$, and we
can try to extend
analytically throughout the simply connected region $W_\po$. There is a
potential obstruction if we ever reach a point in $W_\po$ where the multiplier
$\lambda_n$ of this analytically extended orbit is equal to $+1$.
However, this can never happen. In fact such a point would have to belong
to the Mandelbrot set, and hence to some satellite limb $M_\q$. But we
can extend analytically throughout the associated
wake $W_\q$, taking $\Or^\star(c)$ to be
the repelling orbit $\Or_\q(c)$ for every $c\in W_\q$. Thus there is no
obstruction. It follows similarly that
the analytically extended orbit must be repelling everywhere
in $W_\po\ssm\overline H_\po$. For if it became non-repelling at some
point $c$, then again $c$ would have to belong to some satellite limb
$M_\q$, but $\Or^\star(c)$ is repelling throughout the wake $W_\q$.\QED

\ssk\nin{\tf Corollary 9.3. Limb Connectedness.} {\it Each limb $M_\po=
M\cap\overline W_\po$ is connected, even if we remove its root point $r_\po$.}
\ssk

{\tf Proof.} The entire Mandelbrot set is connected by [DH1]. It follows
that each $M_\po$ is connected. For if some limb
$M_\po$ could be expressed as the union of two disjoint non-vacuous compact
subsets, then only one of these two could contain the root point $r_\po$.
The other would be a non-trivial open-and-closed subset of $M$, which is
impossible.

Now consider the open subset $M_\po\ssm\{r_\po\}$. This is a union of the
connected set $\overline H_\po\ssm\{r_\po\}$, together with the various
satellite limbs $M_{\mathcal Q}$, where each $M_{\mathcal Q}$ has root
point $r_{\mathcal Q}$ belonging to $\overline H_\po\ssm\{r_\po\}$. Since
each $M_{\mathcal Q}$ is connected, the conclusion follows.\QED\ssk

\nin{\tf Remark 9.4.}
It follows easily that every satellite root point separates the Mandelbrot set
into exactly two connected components, and hence that
exactly two parameter rays land at every such point. For a proof of
the corresponding statement for
a primitive root (other than $1/4$) see [Ta] or [S3].

\section*{Appendix A. Totally Disconnected Julia Sets and the Mandelbrot set.}
This appendix will be a brief 
review of well known material. For any parameter value $c$, let
$K=K(f_c)$ be the filled Julia set for the map $f_c(z)=z^2+c$, and let
$$	G(z)~=~G^K(z)~=~\lim_{n\to\infty}\;{1\over 2^n}\,\log\big|f^{\circ n}
(z)\big|$$ be the canonical potential function or {\bit Green's function\/},
which vanishes
only on $K$, and satisfies $G(f(z))=2G(z)$. The level sets
$\{z\;;\;G(z)=G_0\}$ are called {\bit equipotential curves\/} for $K$,
and the orthogonal trajectories which extend to infinity
are called the {\bit dynamic rays\/}
$\ra_t^K$, where $t\in\R/\Z$ is the angle at infinity.

Now suppose that $K$ is totally disconnected (and hence coincides with the
Julia set $J=\partial K$). Then the value $G(0)=G(c)/2>0$
plays a special role.
In fact there is a canonical conformal isomorphism $\psi_c$ from the
open set $\{z\;;\;G(z)>G(0)\}$
to the region $\{w\;;\;\log|w|>G(0)\}$. The map $z\mapsto f(z)$
on this region is conjugate under $\psi_c$
to the map $w\mapsto w^2$, and the equipotentials and
dynamic rays in the $z$-plane correspond to concentric circles and
straight half-lines through the origin respectively in the $w$-plane.
In particular, if we choose a constant $G_0>G(0)$, then the locus
$\{z\;;\;G(z)=G_0\}$ is a simple closed curve, canonically parametrized
by the angle of the corresponding dynamic ray. In particular, the critical
value $c\in\C\ssm K$ has a well defined external angle, which we denote
by $t(c)\in\R/\Z$. Thus $\psi_c(c)/|\psi_c(c)|=e^{2\pi it(c)}$, and
$c$ belongs to the dynamic ray $\ra_{t(c)}=\ra_{t(c)}^K$.

\begin{figure}[htb]
\cl{\psfig{figure=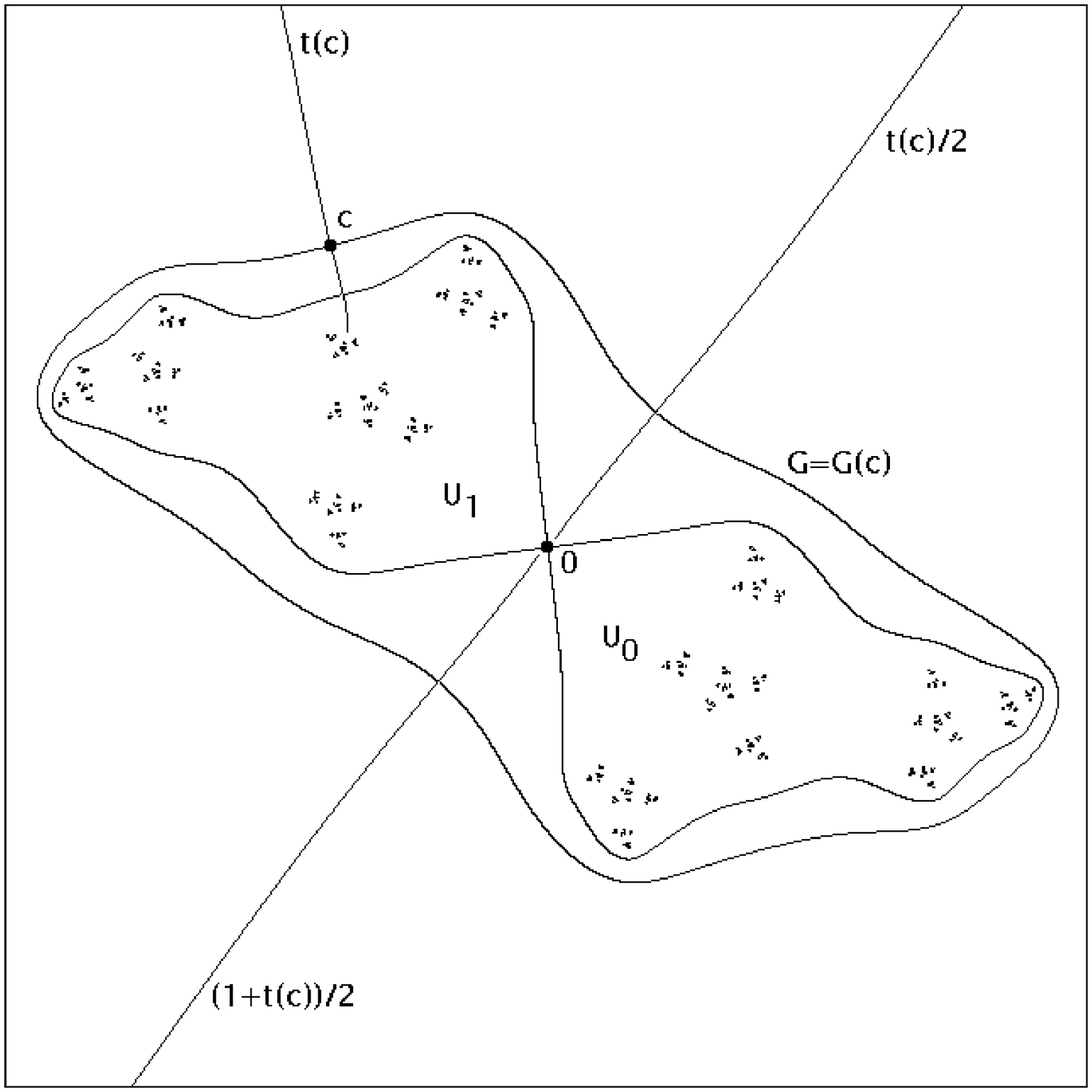,height=2.5in}\qquad
	\psfig{figure=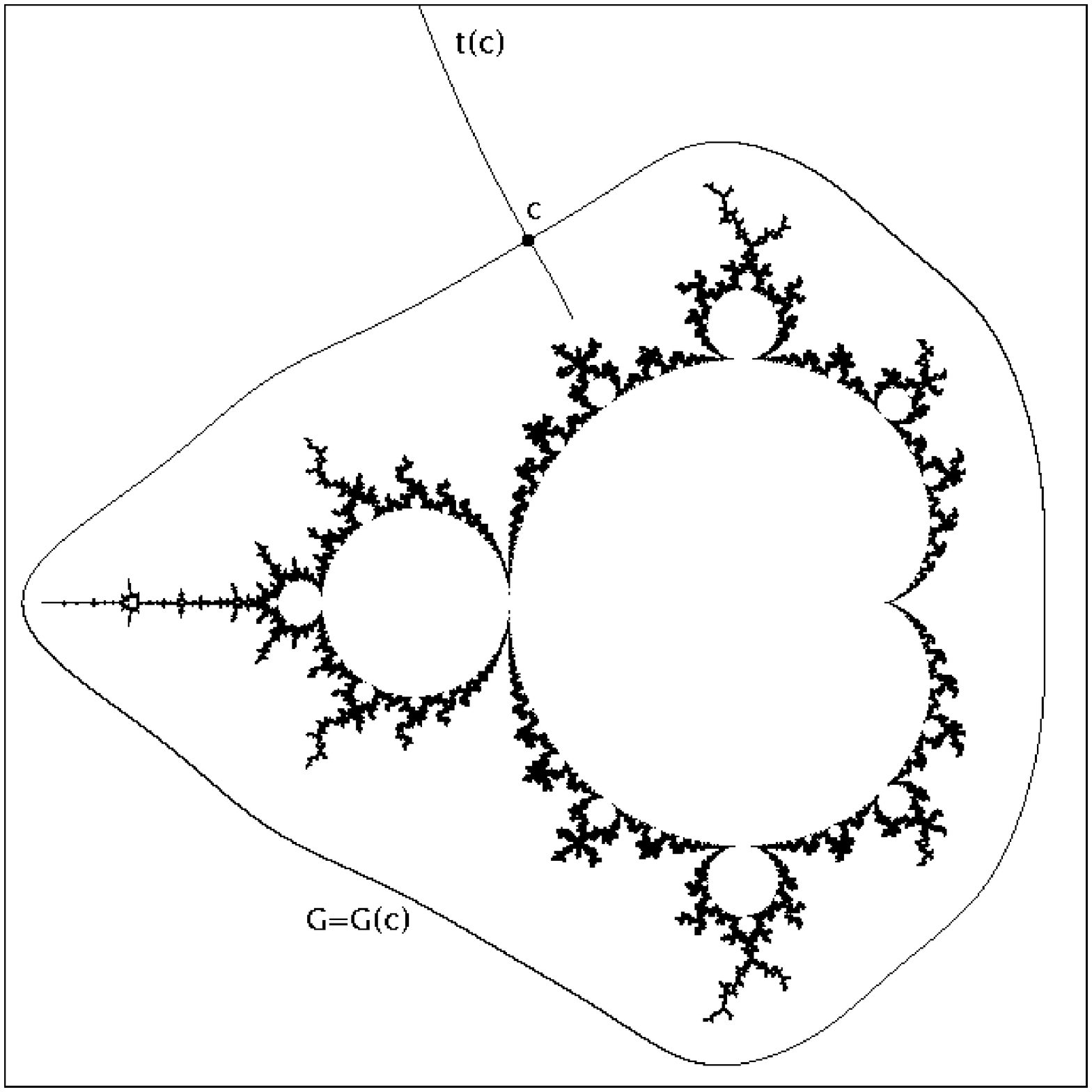,height=2.5in}}\ssk
\begin{quote}{\bit Figure 19. Picture in the dynamic plane
for a polynomial $f_c$ with $c\not\in M$, and a corresponding
picture in the parameter plane .\ssk}\end{quote}
\end{figure}

However,
for $G_0=G(0)$ this locus $\{z\;;\;G(z)=G(0)\}$ is a figure eight
curve. The open set $\{z\;;\;G(z)<G(0)\}$ splits as a disjoint union
$U_0\cup
U_1$, where the $U_b$ are the regions enclosed by the two lobes of
this figure eight. (We can express this splitting in terms of dynamic rays
as follows.
The ray $\ra_{t(c)}^K\subset\C\ssm K$ has two preimage rays under $f_c$,
with angles $t(c)/2$
and $(1+t(c))/2$ respectively. Each of these joins the critical point $0$
to the circle at infinity, and together they cut $\C$ into two open subsets,
say $V_0\supset U_0$ and $V_1\supset U_1$. If $c$ does
not belong to the positive real axis,
then we can choose the labels for these open sets so that the zero ray
is contained in $V_0$, and $c\in V_1$.)
We then cut the filled Julia set $K$ into two disjoint compact
subsets $K_b=K\cap U_b$. These constitute a {\bit Bernoulli
partition.\/} That is, for any one-sided-infinite sequence of bits
$b_0\,,\,b_1\,,\,\ldots\in\{0,1\}$, there is one and only one
point $z\in K$ with $f_c^{\circ k}(z)\in K_{b_k}$ for every $k\ge
0$. To prove this statement, let $U$ be the region $\{z\;;\;G(z)<G(c)\}$
and let $\phi_b:U\to U_b$ be the branch of $f^{-1}$ which maps
$U$ diffeomorphically onto $U_b$. Using the Poincar\'e metric
for $U$, we see that each $\phi_b$ shrinks distances by a factor bounded
away from one, and it follows easily that the diameter of the image
$$	\phi_{b_0}\circ\phi_{b_1}\circ\cdots\circ\phi_{b_n}(U)
 $$
shrinks to zero, so that this intersection shrinks to a single point
$z\in K$, as $n\to\infty$. {\it Thus each point of $J=K$ can be uniquely
characterized by an infinite sequence of symbols
$(b_0\,,\,b_1\,,\,\ldots)$ with $b_j\in\{0,1\}$.\/}
In particular,
$K$ is homeomorphic to the infinite cartesian product $\{0,1\}^{\bf N}$,
where the symbol $\bf N$ stands for the set
$\{0,1,2,\ldots\}$ of natural numbers. We say that the dynamical system
$(K,\,f_c|_K)$ is a {\bit one sided shift\/} on two symbols.

Similarly, given any angle $t\in\R/\Z$, if none of the successive images
$2^k\,t~({\rm mod~}\Z)$ under doubling is precisely equal to $t(c)/2$ or
$(1+t(c))/2$, then $t$ has an associated symbol sequence, called its
$t(c)$-itinerary, and
the ray $\ra_t^K$ lands precisely at that point of $K$ which has this
symbol sequence. For the special case $t=t(c)$, this symbol sequence
characterizes the point $c\in K$, and is called the {\bit kneading
sequence\/} for $c$ or for $t(c)$. (However, if $t(c)$ is periodic,
there is some ambiguity since
the symbols $b_{n-1}, b_{2n-1},\ldots$ of the kneading sequence are
not uniquely defined in the period $n$ case.)

If $t$ is periodic under doubling, then the itinerary is periodic (if uniquely
defined), and the ray $\ra_t^K$
lands at a periodic point of $K$.
For further discussion, see [LS], as well as Appendix B.

Here we have been thinking of $c=f(0)$ as a point in the dynamic plane
(the $z$-plane), but
we can also think of $c\in\C\ssm M$ as a point in the parameter plane
(the $c$-plane).
In fact Douady and Hubbard construct a conformal isomorphism from the
complement of $M$ onto the complement of the closed unit disk by mapping
$c\in\C\ssm M$ to the point $\psi_c(c)=
\exp({G^K(c)+2\pi it(c)})\in\C\ssm\overline{\bf D}$. {\it Thus they show that the
value of the Green's function on $c$ and the external angle $t(c)$ of $c$
are the same whether $c$ is considered as a point of $\C\ssm K(f_c)$ or as
a point of $\C\ssm M$.} In particular, the point $c\in\C\ssm M$ lies on
the external ray $\ra_{t(c)}^M$ for the Mandelbrot set. 

\section*{Appendix B. Computing Rotation Numbers.}
This appendix will outline how to actually compute the
rotation number $q/r$ of a periodic point for a map $f_c$ with
$c\not\in M$. Let $\tau=t(c)\in\R/\Z$ be the angle of the external ray
which passes through $c$. We may identify this critical value
angle with a number in the interval
$0<\tau\le 1$. The two preimages of $\tau$ under the angle doubling map
$m_2:\R/\Z\to\R/\Z$ separate the circle $\R/\Z$ into the two open arcs
$$	I(0)~=~I_\tau(0)~=~ \left(\frac{\tau-1}{2}\;,\;\frac{\tau}{2}\right)
\quad{\rm and}\quad I(1) ~=~I_\tau(1)~=~\left(\frac{\tau}{2}\; ,\;
\frac{\tau+1}{2}\right)~. $$
(We will write $I_\tau$ instead of $I$ whenever we want to emphasize
dependence on the critical value angle $\tau$.)
For any finite sequence $b_0\,,\,b_1\,,\,\ldots\,,\, b_k$
of zeros and ones, let $\overline I(b_0\,,\,b_1\,,\,\ldots\,,\, b_k)$
be the closure of the open set
$$	I(b_0\,,\,b_1\,,\,\ldots\,,\, b_k)~
  =~I(b_0)\cap m_2^{-1}I(b_1)\cap\cdots\cap m_2^{-k} I(b_k) $$
consisting of all $t\in\R/\Z$ with $m_2^{\circ i}(t)\in I(b_i)$
for $0\le i\le k$. (Caution: This is not the same as the intersection of the
corresponding closures $m_2^{-i}\overline I(b_i)$,
which may contain additional isolated points.) An easy induction shows that
$\overline I(b_0\,,\,b_1\,,\,\ldots\,,\, b_k)$
is a finite union of closed arcs with total length
$1/2^{k+1}$. {\it If $\sigma=(b_0\,,\,b_1\,,\,\ldots)$ is any
infinite sequence of zeros and ones, it follows that the intersection
$$  \overline I(\sigma)~=~\bigcap_k\,
\overline I(b_0\,,\,b_1\,,\,\ldots\,,\, b_k)$$
is a compact non-vacuous set of measure zero.}
For each angle $t\in\R/\Z$ there are two possibilities:\ssk

{\bf Precritical Case.} If
$t$ satisfies $m_2^{\circ i}(t)\equiv\tau$ for some $i>0$, then there
will be two distinct infinite symbol sequences with $t\in\overline
I(b_0,b_1,b_2\cdots)$. In this case, the associated dynamic ray $\ra_t^K$ does
not land, but rather bounces off some precritical point for the map $f_c$.
(Compare [GM].)\ssk

{\bf Generic Case.} Otherwise there will be a unique infinite symbol sequence
with $t\in\overline I(b_0, b_1,\cdots)$. The corresponding ray $\ra_t^K$
will land at the unique point of the Julia set for $f_c$ which has this
same symbol sequence, as described in Appendix A. In particular, if $t$
is periodic under doubling, then $\ra_t^K$ must land at a periodic point
of the Julia set, possibly with smaller period.\ssk

\ssk{\nin{\tf Lemma B.1. Symbol Sequences and Rotation Numbers.} \it
For any symbol sequence $\sigma=(b_0,b_1,\ldots)
\in\{0,1\}^{\bf N}$ which is periodic of period $p$, the map $m_2^{\circ p}$
on the compact set $\overline I_\tau(\sigma)\subset\R/\Z$
has a well defined rotation number
$\rot(b_0,\ldots, b_{p-1}\,;\,\tau)\in\R/\Z$ which is invariant under
cyclic permutation of the bits $b_i$. This number
increases monotonically with $\tau$, and winds $b_0+\cdots+b_{p-1}$
times around the circle as $\tau$ increases from $0$ to $1$.}\ssk

To see this, we introduce an auxiliary monotone degree one map which is
defined on the entire circle and agrees with $m_2^{\circ p}$ on
$\overline I_\tau(\sigma)$. (Compare [GM].)
By definition, a {\bit monotone degree one circle map\/} $\psi:\R/\Z\to\R/\Z$
is the reduction modulo $\Z$ of a map $\Psi:\R\to\R$ which is monotone
increasing and satisfies the identity $\Psi(u+1)=\Psi(u)+1$. Such a $\Psi$,
called a {\bit lift\/} of $\psi$, is unique up to addition of an integer
constant. The {\bit translation number\/}
of such a map $\Psi$ is defined to be the real number
$$	\Rot(\Psi)~=~\lim_{k\to\infty} \big(\Psi^{\circ k}(u)-u\big)/k~.$$
This always exists, and is independent of $u$.
The {\bit rotation number\/}
$\rot(\psi)$ of the associated circle map is now defined to be the image
of this real number $\Rot(\Psi)$ under the projection $\R\to\R/\Z$.
This is well defined, since $\Rot(\Psi+1)=\Rot(\Psi)+1$. One important
property is the identity
$$	\Rot(\Psi_1\circ\Psi_2)~=~\Rot(\Psi_2\circ\Psi_1)~, \eqno (9)$$
where $\Psi_1$ and $\Psi_2$ are the lifts of two different monotone degree
one circle maps.
If $\Psi_1$ is a homeomorphism, this is just invariance under a suitable
change of coordinates, and the general case follows by continuity.

Given any $b\in\{0,1\}$, and given a critical value angle
$\tau$, define an auxiliary monotone map $\Phi_{b,\tau}$
by the formula
$$	\Phi_{b,\tau}(u)~=~\begin{cases} \min(2u\,,\,\tau) & \text{if}~~ b=0,\cr
	 \max(2u\,,\,\tau) & \text{if}~~ b=1,\cr\end{cases}$$
for $u$ between $(\tau-1)/2$ and $(\tau+1)/2$, extending by the
identity $\Phi(u+1)=\Phi(u)+1$ for $u$ outside this interval. (See
Figure 20.) Note that
$\overline I(b)$ is just the set of points on the circle where the associated
circle map $\phi_{b,\tau}$ is not locally constant, and that
$\phi_{b,\tau}(u)\equiv 2u~({\rm mod~}\Z)$ whenever $u\in \overline I(b)$.

\begin{figure}[htb]
\cl{\psfig{figure=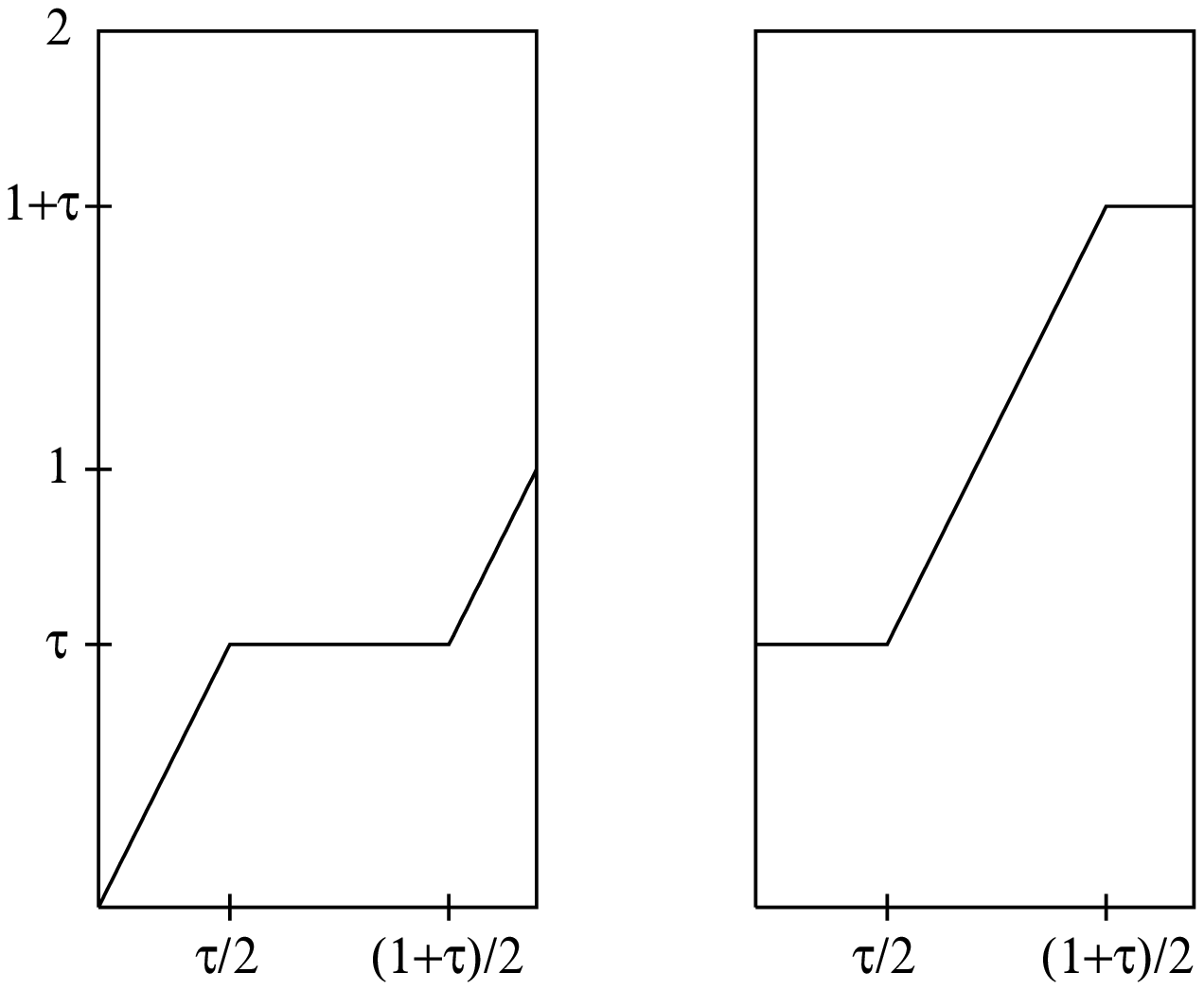,height=2.4in}}\ssk

\cl{\bit Figure 20. Graphs of $\Phi_{0,\tau}$ and $\Phi_{1,\tau}$ (with
$\tau=0.6$.)}
\end{figure}

For any symbol sequence $\sigma$ which is periodic of period $p$,
we set $\Phi_{\sigma,\tau}$ equal to the $p$-fold composition
$\Phi_{b_{p-1},\tau}\circ\cdots\circ
\Phi_{b_{0},\tau}$. (Note that $\overline I_\tau(\sigma)$ is just the set of
all points $t\in\R/\Z$
such that the orbit of $t$ under the associated circle map
$\phi_{\sigma,\tau}$ coincides
with the orbit of $t$ under $m_2^{\circ p}$.) This composition
is also monotone, with $\Phi(t+1)=\Phi(t)+1$, and
therefore has a well defined translation number, which we denote by
$$	\Rot(b_0,\cdots,b_{p-1}\;;\;\tau)~=~\Rot(\Phi_{\sigma,\tau})\;\in\;
  \R~.$$
It follows from property (9) that
this translation number is invariant under cyclic permutation of the bits
$b_0\,,\,\ldots\,,\,b_{p-1}$.
Since each $\Phi_{b,\tau}(u)$ increases monotonically with $\tau$, with
$\Phi_{b,0}(0)=0$ and $\Phi_{b,1}(0)=b$, it follows easily that
$\Rot(\Phi_{\sigma,\tau})$ depends monotonically on
$\tau$, increasing from $0$ to
$b_0+\cdots+b_{p-1}$ as $\tau$ increases from $0$ to $1$. In other
words its image in $\R/\Z$ wraps $b_0+\cdots+b_{p-1}$ times around the
circle as $\tau$ varies from $0$ to $1$. By definition, the rotation number
$\rot(b_0,\ldots ,b_{p-1}\,;\,\tau)$ of
$m_2^{\circ p}$ on the compact set $\overline I_\tau(\sigma)$ is equal
to the image of the real number $\Rot(\Phi_{\sigma,\tau})$ in the
circle $\R/\Z$. \QED\ssk

If a map $f_c$ has critical value angle $t(c)=\tau$, then it is not hard to
see that $\rot(b_0,\ldots ,b_{p-1}\,;\,\tau)$ coincides with the rotation
number as defined in 2.12 for the orbit with periodic symbol sequence
$\overline{b_0,\ldots ,b_{p-1}}=(b_0,\ldots ,b_{p-1},b_0,\ldots
 ,b_{p-1},\ldots)$, so long as at least one rational ray lands on this orbit.
(Compare [GM, Appendix C].)\ssk

\begin{figure}[htb]
\bs
\cl{\psfig{figure=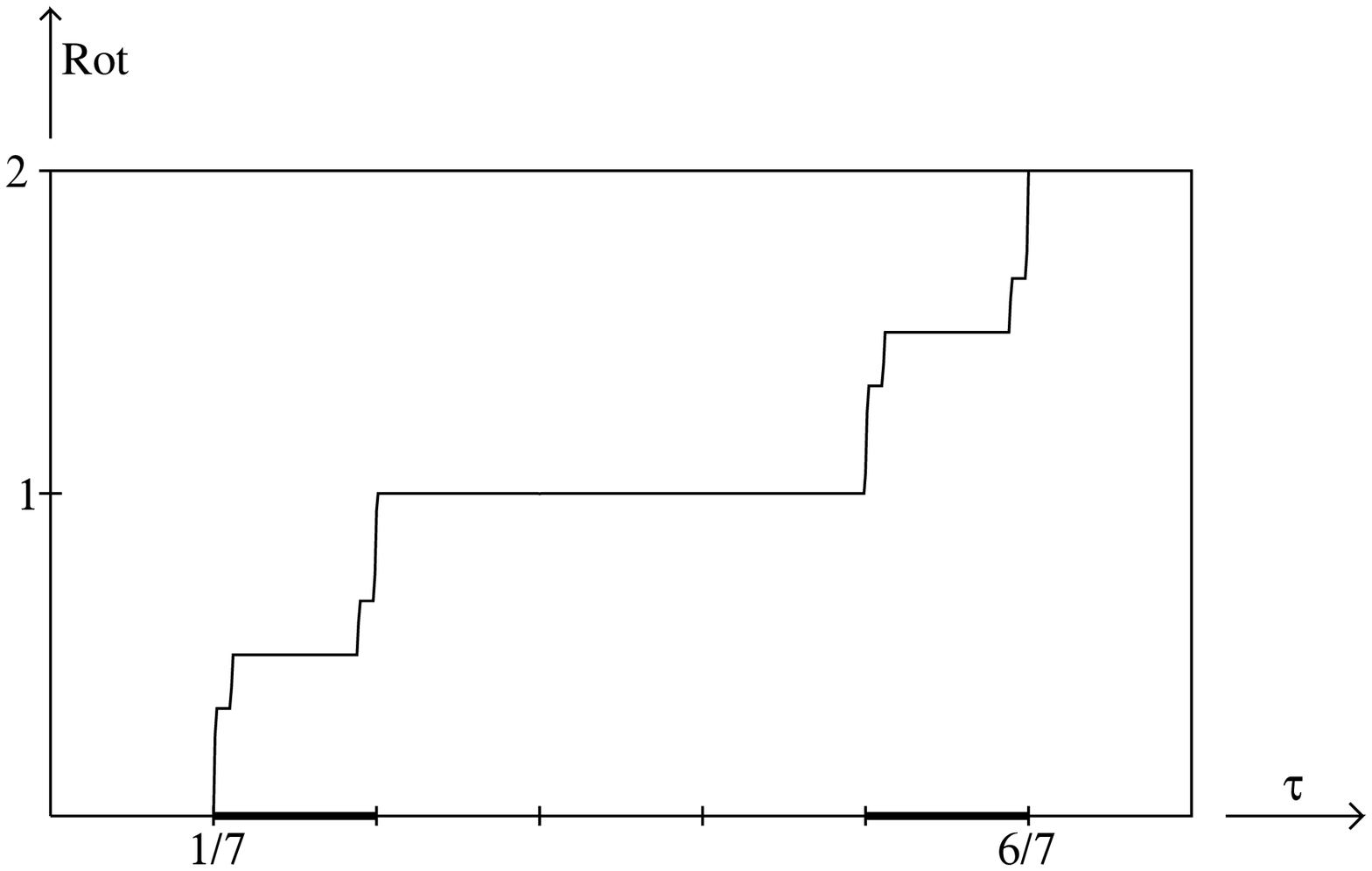,height=1.8in}}\ms
\begin{quote}
{\bit Figure 21. The translation number as a function of the critical value
exterior angle $\tau$ for the period 3 point with symbol sequence
$\overline{110}=( 1, 1, 0,1,1,0,\ldots)$.
}\end{quote}
\end{figure}

We will use the notation ${\mathcal S}(q/r)$ for the orbit
portrait with orbit period $p=1$ and rotation number $q/r$,
associated with the $q/r$-satellite of the
main cardioid. (Compare [G].) If $\po$ is an arbitrary orbit portrait, then
$\po*{\mathcal S}(q/r)$ can be described as its $(q/r)$-{\bit satellite
portrait\/}. (See 6.4, 8.2.)

To any orbit portrait $\po$ with period $p\ge 1$ and ray period $n=rp\ge p$
we can associate a symbol sequence $\sigma=\sigma(\po)$ of period
$p$ as follows. Choose any $c\not\in M$ in the wake $W_\po$, and
number the points of the $f_c$-orbit with portrait $\po$ as
$z_0\mapsto z_1\mapsto\cdots$, where $z_0$ is on the boundary of the critical
puzzle piece and $z_1$ is on the boundary of the critical value puzzle piece.
Now let $\sigma(\po)$ be the
symbol sequence for $z_0$, as described in Appendix A. This is independent of
the choice of $c\in W_\po\ssm M_\po$.

There is an associated {\bit satellite symbol sequence\/} $\sigma^\star=
\sigma^\star(\po)$ of period $n=rp$, constructed as follows. (Compare 9.2.)
By definition, the $k$-th bit of $\sigma^\star$ is identical to the $k$-th
bit of $\sigma$ for $k\not\equiv 0~({\rm mod}~n)$, but is reversed,
so that $0\leftrightarrow 1$, when $k\equiv 0~({\rm mod}~ n)$.

\ssk{\nin{\tf Lemma B.2. Satellite Symbol Sequences.}
\it For every satellite $\po*{\mathcal S}(q'/r')$
of $\po$, the symbol sequence $\sigma(\po*{\mathcal S}
(q'/r'))$ coincides with the satellite sequence
$\sigma^\star(\po)$. The translation number $\Rot(\sigma(\po),\tau)$ is
constant for $\tau$ in the characteristic arc $\I_\po$, while
$\Rot(\sigma^\star(\po),\tau)$ increases by $+1$ as $\tau$ increases through
$\I_\po$, taking the value $q'/r'~({\rm mod}~ \Z)$
on the characteristic arc of $\po*{\mathcal S}(q'/r')$.\ssk}

Intuitively, if we tune a map in $H_\po$ by a map in $H_{{\mathcal S}(q'/r')}$
then we must replace the Fatou component containing the critical point
for the first map by a small copy of the filled Julia set for a
$(q'/r')$-rabbit. Here the period $p$ point $z_0$ for $\po$ corresponds to
the $\beta$-fixed point of this small rabbit, while the period $n$ point
$z_0$ for $\po*{\mathcal S}(q'/r')$ corresponds to the $\alpha$ fixed
point for this rabbit. Perturbing out of the connectedness locus $M$,
these two points will be separated by the ray pair terminating at the
critical point. Further details will be omitted.\QED\ssk

For example, starting with $\sigma\big(\{\{0\}\}\big)=\overline 0$,
where the overline indicates infinite repetition, we find that
$$	~\sigma({\mathcal S}(q/r))~=~
\sigma^\star\big(\{\{0\}\}\big)~=\overline 1~,$$
while
$$ \sigma^\star({\mathcal S}(1/2))~=~\overline{0 1}~,\quad
 \sigma^\star({\mathcal S}(q/3))~=~\overline{0 1 1}~,\quad
 \sigma^\star({\mathcal S}(q/4))~=~\overline{0 1 1 1}~,\quad \ldots \,.$$

We can use this discussion to provide a different insight on the counting
argument of \S5. Since $\Rot(\sigma^\star(\po)\,;\,\tau)$ increases by $+1$
on the characteristic arc $\I_\po$, we see that the total number
of portraits (or the total number of
characteristic arcs) with ray period $rp=n$ is equal to the
sum of $b_0+\cdots+b_{n-1}$ taken over all cyclic equivalence classes
of symbol sequences of period exactly $n$. But the number of such symbol
sequences, up to cyclic permutation, is $\nu_2(n)/n$, and the average value
of $b_0+\cdots+b_{n-1}$ is equal to $n/2$, since each symbol sequence
with sum different from $n/2$ has
an opposite with zero and one interchanged. Therefore, this sum is equal to
$\nu_2(n)/2$, as in \S5.\ssk

\ssk{\tf Examples}  (Compare Figure 4).
Here is a list for all cyclic equivalence classes
of symbol sequences of period at most four:\ssk

\Ref $\Rot(0\;;\;\tau)$ is identically zero.

\Ref $\Rot(1\;;\;\tau)$ increases from $0$ to $1$ for $0\le\tau\le 1$,
taking the value $q/r$ in the characteristic arc for ${\mathcal S}(q/r)$.

\Ref $\Rot(1,0\;;\;\tau)$ increases from $0$ to $1$ as $\tau$ passes through
$(1/3\,,\,2/3)$, the characteristic arc for ${\mathcal S}(1/2)$.

\Ref$\Rot(1,0,0\;;\;\tau)$ increases from $0$ to $1$ as $\tau$ passes through
the characteristic arc $(3/7\,,\,4/7)$ for the period
3 portrait with root point $c=-1.75$.

\Ref$\Rot(1,1,0\;;\;\tau)$ increases by one in the arc $(1/7\,,\,2/7)$ for
${\mathcal S}(1/3)$, and by one more in the arc $(5/7\,,\,6/7)$ for
${\mathcal S}(2/3)$. (Compare Figure 21.)

\Ref$\Rot(1,0,0,0\;;\;\tau)$ increases by one in the arc $(7/15\,,\,8/15)$,
corresponding to the leftmost period 4 component on the real
axis.

\Ref $\Rot(1,1,0,0\;;\;\tau)$ increases by one in the arcs $(1/5\,,\,4/15)$
and $(11/15\,,\,4/5)$ associated with the
period 4 components on the $1/3^{\rm rd}$ and $2/3^{\rm rd}$ limbs.
(Figure 12.)

\Ref
$\Rot(1,1,1,0\;;\;\tau)$ increases by one in the arcs $(1/15\,,\,2/15)$ and
$(13/15\,,\,14/15)$ for ${\mathcal S}(1/4)$ and ${\mathcal S}(3/4)$,
and also in the arc
$(2/5\,,\,3/5)$  for the portrait ${\mathcal S}(1/2)*{\mathcal S}(1/2)$
with root point $-1.25$.


\vskip .5in

\hskip 2.7in{John Milnor}

\hskip 2.7in{Institute for Mathematical Sciences}

\hskip 2.7in State University of New York

\hskip 2.7in Stony Brook, NY 11794-3660\ssk

\hskip 2.7in{jack@math.sunysb.edu}

\hskip 2.7in{http://www.math.sunysb.edu/$\sim$jack}

\vfil

\end{document}